\pdfoutput=1
\documentclass[final,leqno]{siamltex}

\usepackage{graphicx}
\usepackage{bm}
\usepackage{bbm}
\usepackage{amsmath}
\usepackage{amssymb}
\usepackage{color}
\usepackage{enumitem}
\usepackage{units}

\newcommand{\IGNORE}[1]{}
\newcommand{\ignore}[1]{}
\newcommand{\veps}{\varepsilon}
\newcommand{\opn}{\operatorname}
\newcommand{\re}{\operatorname{Re}}
\newcommand{\im}{\operatorname{Im}}

\newcommand{\mbb}[1]{\mathbb{#1}}

\newcommand{\mc}[1]{\mathcal{#1}}

\newcommand{\jd}{\displaystyle}

\newcommand{\jt}{\textstyle}

\newcommand{\der}[2]{\frac{\partial #1}{\partial #2}}

\newcommand{\erf}{\opn{erf}}
\renewcommand{\H}{\mathcal{H}}

\newcommand{\rlam}{\tau}

\newcommand{\la}{\langle}
\newcommand{\ra}{\rangle}

\newcommand{\seglen}{\ell}
\newcommand{\nn}{n}
\newcommand{\bern}{\alpha}

\newcommand{\sgn}{\beta}

\title{A Spectral Transform Method for Singular Sturm--Liouville Problems
       with Applications to Energy Diffusion in Plasma Physics}

\author{
  Jon Wilkening
  \thanks{Department of Mathematics and Lawrence
    Berkeley National Laboratory, University of California, Berkeley,
    CA 94721 ({\tt wilkening@berkeley.edu}). This work was supported
    in part by the US Department of Energy, Office of Science, Applied
    Scientific Computing Research, under award number
    DE-AC02-05CH11231, and by the National Science Foundation under
    award number DMS-0955078.}
  \and
  Antoine Cerfon
  \thanks{Courant
    Institute of Mathematical Sciences, New York, NY 10012 ({\tt
      cerfon@cims.nyu.edu}). This work was supported in part by the
    U.S. Department of Energy, Office of Science, Fusion Energy
    Sciences under award number DE-FG02-86ER53223.}
}

\begin{document}

\maketitle

\begin{abstract}
  We develop a spectrally accurate numerical method to compute
  solutions of a model partial differential equation used in plasma
  physics to describe diffusion in velocity space due to Fokker-Planck
  collisions.  The solution is represented as a discrete and
  continuous superposition of normalizable and non-normalizable
  eigenfunctions via the spectral transform associated with a singular
  Sturm-Liouville operator.  We present a new algorithm for computing
  the spectral density function of the operator that uses Chebyshev
  polynomials to extrapolate the value of the Titchmarsh-Weyl
  $m$-function from the complex upper half-plane to the real axis.
  The eigenfunctions and density function are rescaled and a new
  formula for the limiting value of the $m$-function is derived to
  avoid amplification of roundoff errors when the solution is
  reconstructed.  The complexity of the algorithm is also analyzed,
  showing that the cost of computing the spectral density function at
  a point grows less rapidly than any fractional inverse power of the
  desired accuracy. A WKB analysis is used to prove that the spectral
  density function is real analytic.  Using this new algorithm, we
  highlight key properties of the partial differential equation and
  its solution that have strong implications on the optimal choice of
  discretization method in large-scale plasma physics computations.
\end{abstract}

\begin{keywords}
Sturm-Liouville theory,
continuous spectrum,
Titchmarsh-Weyl $m$-function,
spectral density function,
Fokker-Planck collisions,
WKB approximation
\end{keywords}

\section{Introduction}

Partial differential equations involving singular Sturm-Liou\-ville
operators with continuous spectra arise frequently in computational
physics.  Common approaches to solving them include domain truncation,
which often regularizes the operator and makes the spectrum discrete,
or projection onto finite dimensional orthogonal polynomial or finite
element subspaces, which also leads to discrete spectra.  Here we
develop an alternative approach in which the continuous spectrum is
treated analytically via a spectral transform, and the numerical
challenge is in accurately representing and evaluating the integrals
giving the exact solution.

While the methods developed in this paper to diagonalize singular
Sturm-Liouville operators are quite general, we will describe them
in the context of velocity-space diffusion in one dimension,
\begin{equation}\label{eq:U:intro}
  \der{h_a}{t} = \frac{1}{v^2}\der{}{v}\left[
    \tilde\Psi(v)\left(2v^2 h_a + v \der{h_a}{v}\right)\right], \qquad
  (v>0,\;t>0),
\end{equation}
where $\tilde\Psi(v) = \big[\erf(v) -
v\erf'(v)\big]/(2v^2)$ is the Chandrasekhar function and
$\erf(v)=2\pi^{-1/2}\int_0^v e^{-r^2}dr$ is the error function. The
diffusion operator on the right-hand side of (\ref{eq:U:intro}) plays
an important role in numerical simulations of systems governed by the
Fokker-Planck kinetic equation
\cite{hazeltine_meiss,hazeltine_waelbroeck}
\begin{equation}\label{eq:full:FP}
  \der{f_a}{t} + \mathbf{v}\cdot\nabla f_a + \frac{q_a}{m_a}(\mathbf{E} +
  \mathbf{v}\times
  \mathbf{B})\cdot \nabla_\mathbf{v}f_a = \sum_b C(f_a,f_b).
\end{equation}
Here $f_a(\mathbf{x},\mathbf{v},t)$ is the distribution function for
particles of species $a$ of charge $q_a$ and mass $m_a$, and
$C(f_a,f_b)$ is the Fokker-Planck collision operator \cite{rosenbluth,
  helander, pataki} describing collisions between species $a$ and
other species, including itself. In many problems of
interest, the distribution function $f_{a}$ is close to a
Maxwell-Boltzmann distribution, in which case one often writes
$f_a=f_{Ma}+h_a$, where $f_{Ma}$ is the Maxwell-Boltzmann distribution
and $h_a\ll f_{Ma}$ \cite{candy,barnes1,landreman2}. The collision
operator $C(f_a,f_b)$ may then be linearized about $f_{Ma}$
\cite{helander,barnes2,landreman1}. When this is done, the operator on
the right-hand side of (\ref{eq:U:intro}) is the term in the
linearized version of $C(f_a,f_b)$ that represents energy diffusion
resulting from the collision of $h_a$ with a Maxwellian background
$f_{Ma}$ \cite{helander}. In this context, the variable $v$ in
Equation (\ref{eq:U:intro}) is the speed coordinate
$|\mathbf{v}|$.

Numerical solution of (\ref{eq:full:FP}) is expensive due to the
high-dimensional phase space in which the distribution functions
evolve \cite{candy,barnes1,barnes2}. To reduce the computational time
without sacrificing accuracy, it is important to develop optimized
discretization techniques \cite{bratanov,landreman1}. For the speed
coordinate, non-classical orthogonal polynomials
\cite{shizgal,landreman1,ghiroldi} are emerging as a promising
alternative to finite difference methods \cite{candy,barnes3}. Equation
(\ref{eq:U:intro}) is well-suited to assess the merits of these
discretization schemes for the speed coordinate in (\ref{eq:full:FP})
without the computational overhead of a high-dimensional phase
space. It is physically relevant since the right-hand side
of (\ref{eq:U:intro}) can be identified as the energy-diffusion term
in the ``test-particle'' part of the linearized Landau collision
operator \cite{helander,abel}, as already discussed. The solution also
relaxes to a Maxwellian distribution as $t\rightarrow\infty$, as one
would physically expect from a collision operator, and the equation is
``mass-conserving,'' i.e.~the integral $\int_0^\infty 4\pi v^2
h_a(v,t)\,dv$ remains constant in time. As such, (\ref{eq:U:intro})
and close variants of it are often used as standard tests of the
accuracy and conservation properties of numerical solvers, as
illustrated by Example IV.A in \cite{barnes2} for
instance.

In the present paper, we develop a spectral transform method to study
the dynamics of (\ref{eq:U:intro}) in detail.  In subsequent work
\cite{vsck2}, joint with Landreman, we will study the projected
dynamics of this equation in finite-dimensional spaces of orthogonal
polynomials.  Roughly speaking, we show in this article how to efficiently
evaluate the exact solution by discretizing a continuous transform,
while in \cite{vsck2} we discretize the PDE before evolving the
solution. The latter approach is faster and better suited to large
scale computations of the full Fokker-Planck equation, while the
current approach clarifies the role of the continuous spectrum in the
dynamics and provides an independent means of validating the
orthogonal polynomial approach. Validation is especially important in
singular cases where the true solution of (\ref{eq:U:intro}) leaves
the finite-dimensional subspace and later returns to a point that may
or may not agree closely with the solution of the projected
dynamics. These singular cases are not only of academic interest, but
in fact correspond to situations of practical interest in plasma
physics, such as the calculation of the resistivity of a homogeneous,
unmagnetized plasma \cite{landreman1}. The construction of an accurate
discretization of the exact solution of (\ref{eq:U:intro}) with
singular initial conditions is therefore a useful tool to evaluate
the performance and accuracy of numerical solvers for the
Fokker-Planck equation, and to understand their behavior. It is also
of intrinsic theoretical interest, with applications beyond plasma
physics, to be able to diagonalize differential operators with
continuous spectra.

Existing algorithms for computing spectral density functions of
singular Sturm-Liouville problems employ either a domain truncation
technique (accelerated by Rich\-ard\-son extrapolation)
\cite{pruess1993,fulton1994,fulton1998}, or use formulas for the
density function \cite{titchmarsh,fulton2005,fulton2008b,fulton2008c}
that avoid domain truncation but only apply to Sturm-Liouville
problems in standard form, $Lu=-u''+q(x)u=\lambda u$.  Further details
about both approaches are given in Sections~\ref{sec:comp}
and~\ref{sec:compare}, and in \cite{wilkening:irk}.  While it is
always possible to reduce to standard form via the Liouville
transformation \cite{milson1997}, it is often preferable to work with
the most physically relevant variables.  Our idea is to turn Weyl's
original construction \cite{weyl} into a numerical algorithm by
complexifying the spectral parameter $\lambda=\tau+i\veps$ and
extrapolating the Titchmarsh-Weyl $m$-function to the real axis using
Chebyshev interpolating polynomials. In more detail, when $\lambda$ is
complex, $m(\lambda)=-\lim_{x\rightarrow\infty}
[u_0(x;\lambda)/u_1(x;\lambda)]$ is the limiting ratio of two
solutions of $Lu=\lambda u$.  We show that this limit converges
exponentially fast relative to the work required to compute the
solutions $u_0$ and $u_1$, with a decay rate proportional to
$\im\{\sqrt{\lambda}\}/|\lambda|^{1/2}$ (as a function of
work). Making use of arbitrary-order ODE solvers and interpolation
methods, the cost of the new algorithm for computing the spectral
density function $\rho'(\tau)=(1/\pi)\lim_{\veps\rightarrow0^+}
\im\{m(\tau+i\veps)\}$ to a tolerance $\delta$ grows slower than
$\delta^{-\gamma}$ (for any $\gamma>0$) as $\delta\rightarrow0$.
Proving this requires precise information about the asymptotic
behavior of solutions of $Lu=\lambda u$ for large $x$. We present
a WKB analysis in Appendix~\ref{sec:asym} and a proof of analyticity
of $\rho'(\lambda)$ in Appendix~\ref{sec:analytic}
for the Sturm-Liouville operator associated with (\ref{eq:U:intro}).

In addition to developing a new algorithm for computing spectral
density functions of singular Sturm-Liouville problems of the form $Lu=\lambda u$,
we show how to use them to evaluate the solution of $u_t=-Lu$ at any
later time.  Through appropriate changes of variables, the spectral
transform of the solution, $\hat u(\lambda,t)$, can be represented
accurately and concisely as a trigonometric polynomial. The solution
$u(x,t)$ in physical space is then expressed as an oscillatory
integral in $\lambda$.  For some initial conditions, the spectral
transform $\hat u(\lambda,t)$ decays exponentially in $\lambda$ for
$t>0$ but only algebraically at $t=0$. Thus, with limited
computational resources, the solution of $u_t=-Lu$ often cannot be
resolved to the desired level of accuracy until $t$ surpasses a
critical value, $t^*$, where the decay rate of $\hat u(\lambda,t^*)$
becomes fast enough.

Remarkably, the same is true of the projected dynamics in some spaces
of orthogonal polynomials \cite{vsck2}. For singular initial
conditions, the projected dynamics is a poor approximation of the true
solution initially, regardless of which space of polynomials is used
to represent the solution.  However, the true solution will generally
return (very nearly) to the space once $t$ exceeds some $t^*$. For the
class of orthogonal polynomials introduced by Shizgal \cite{shizgal}
and Landreman and Ernst \cite{landreman1}, the true solution agrees
with the projected dynamics to 29 digits of accuracy for $t>t^*$ (in
quadruple-precision arithmetic). By contrast, for classical Hermite
polynomials, it only agrees to 2--3 digits of accuracy for similar
computational work.  Thus, in one case the projected dynamics evolves
to the correct state when $t$ reaches $t^*$, while in the other case
it does not.  The methods of the current paper were developed in order
to quantify these errors and understand these results.

\section{Preliminaries}
\label{sec:problem}

Our goal is to develop a spectral representation for solutions of the
PDE (\ref{eq:U:intro}). For notational convenience, we will use the variable $x$
instead of $v$ for the speed coordinate. The equation may then be
written
\begin{align}
  \label{eq:U}
  \der{h_a}{t} = \frac{1}{x^2}\der{}{x}\left[
    \Psi(x)x^2 e^{-x^2}\der{}{x}\Big(e^{x^2}h_a\Big)\right],
  \qquad (x>0,\;t>0),
\end{align}
where
\begin{align}\label{eq:psi}
  \Psi(x) = \frac{1}{2x^3}\left[\erf(x) - \frac{2}{\sqrt{\pi}}xe^{-x^2}\right], \qquad
  \erf(x) = \frac{2}{\sqrt{\pi}}\int_0^x e^{-s^2}\,ds.
\end{align}
Several properties of $\Psi(x)$, which differs from the Chandrasekhar
function $\tilde\Psi(x)$ in the introduction by a factor of $x$,
are established in Lemma~\ref{lem:psi} of Appendix~\ref{sec:lem}.
In particular,
$\Psi(0)=2/(3\sqrt\pi)\approx0.3761$, $\Psi'(0)=0$, $\Psi(x)$ is
monotonically decreasing for $x\ge0$, and $\Psi(x)\sim(2x^3)^{-1}$ as
$x\rightarrow\infty$. Furthermore, it is clear that as the ratio of two odd functions, $\Psi(x)$ is even.

We begin by transforming (\ref{eq:U}) to a self-adjoint system.  Let
\begin{equation}
  u(x,t) = h_a(x,t)e^{x^2}.
\end{equation}
Then $u$ satisfies
\begin{equation}\label{eq:PDE}
  u_t = -Lu, \qquad
  Lu = -\frac{(\Psi w u')'}{w}, \qquad
  w(x) = x^2e^{-x^2},
\end{equation}
where $u_t:=\partial u/\partial t$ and the prime symbol stands for the
derivative with respect to $x$.  The domain $D$ of $L$ can be
characterized precisely (see \cite{coddington} and \S\ref{sec:green}),
but it is difficult to show that $L$ is symmetric on all of $D$
directly.  Thus, initially, we will work with the set $D_1$ of
bounded, $C^2$ functions on $(0,\infty)$ with two bounded derivatives.
Such functions extend continuously to $x=0$ with finite limiting value
and slope.  Since
\begin{equation}
  \langle Lu,v\rangle = \langle u,Lv\rangle, \qquad\qquad
  (u,v\in D_1),
\end{equation}
where $\langle u,v\rangle = \int_0^\infty u(x)\overline{v(x)}\,w(x)\,dx$,
we see that the (densely defined) restriction operator
$L_1=L\restriction D_1$ is symmetric on the Hilbert space
\begin{equation}\label{eq:H}
  \H = L^2(\mathbb{R}_+;w\,dx) =
  \big\{u\;:\;\int_0^\infty |u(x)|^2w(x)\,dx < \infty\big\}.
\end{equation}
$L$ is defined as the graph closure of $L_1$, which exists since $L_1$
is symmetric. A well-known theorem \cite{reed:simon} asserts that $L$
is self-adjoint iff $\opn{Ran}(L_1+i)$ and $\opn{Ran}(L_1-i)$ are
dense in $\H$, which may be proved by construction of a Green's
function (see \S\ref{sec:green}).  We note that $D_1$ must be small
enough that $L_1$ is symmetric but large enough that $D_1$ and
$\opn{Ran}(L_1\pm i)$ are dense in $\H$.  For singular Sturm-Liouville
operators, this boils down to imposing the correct boundary conditions
at the endpoints.

\subsection{Classification of the endpoints}
\label{sec:classify}

The operator $L$ in (\ref{eq:PDE}) is singular at $x=0$ since
$w(0)=0$, and at $x=\infty$ since the domain is unbounded.
We now show that $L$ is of limit circle type at $x=0$ and limit point
type at $x=\infty$ \cite{coddington, stakgold, krall}.
To classify the endpoints, we study the behavior of solutions of
\begin{equation}\label{eq:eval}
  -(\Psi w u')' = \lambda w u, \qquad\quad (\lambda\in\mathbb{C})
\end{equation}
as $x\rightarrow 0$ and $x\rightarrow\infty$.  When $\lambda=0$,
the general solution is
\begin{equation}\label{eq:u:lam0}
  u(x) = \alpha_1 - \alpha_0 \int_1^x \frac{e^{s^2}}{s^2\Psi(s)}\,ds, \qquad
  (0<x<\infty).
\end{equation}
The integrand may be expanded in a Laurent series about $s=0$ to obtain
\begin{equation}
  \frac{e^{s^2}}{s^2\Psi(s)} = \frac{3\sqrt{\pi}}{2s^2} + \frac{12\sqrt{\pi}}{5}
  + O(s^2).
\end{equation}
Thus,
\begin{equation*}
  u(x) = \frac{3\sqrt{\pi}}{2x} \alpha_0 + O(1), \qquad\quad
  (x\ll 1).
\end{equation*}
Since $1/x$ belongs to $\H$ in (\ref{eq:H}), all solutions of
(\ref{eq:eval}) are square-integrable on $(0,1)$ with weight function
$w(x)$ when $\lambda=0$.  Weyl's theorem \cite{coddington,
  stakgold, krall} states that this is true for all
$\lambda\in\mathbb{C}$ if it is true for one $\lambda$.  Thus, the
limit circle case prevails at $x=0$.  

The situation is different at $x=\infty$. Since
$\lim_{s\rightarrow\infty}s^3\Psi(s)=1/2$, there is an
$x_0\in(0,\infty)$ such that
\begin{equation}
  \frac{e^{s^2}}{s^2\Psi(s)}\ge s e^{s^2}, \qquad\quad (s\ge x_0).
\end{equation}
It follows that $u(x)$ in (\ref{eq:u:lam0}) with $\alpha_1=0$ and $\alpha_0=-1$
satisfies
\begin{equation*}
  u(x) = \bigg(u_0 + \int_{x_0}^x \frac{e^{s^2}}{s^2\Psi(s)}\,ds\bigg)
  \ge \bigg(u_0 + \int_{x_0}^x s e^{s^2}\,ds\bigg) =
  \bigg(u_0 + \frac{e^{x^2} - e^{x_0^2}}{2}\bigg)
\end{equation*}
for $x\ge x_0$.  The function on the right is not square integrable on
$(x_0,\infty)$ with weight function $w(x)$, so neither is $u$ and the
limit point case prevails at $x=\infty$.

\subsection{Rescaled variables}
\label{sec:rescaled}

The limit circle case requires a boundary condition. It suffices for this
to require that solutions of $Lu=\lambda u$ remain bounded at $x=0$.
However, a linearly independent solution (that blows up at the origin)
must also be computed to evaluate the Titchmarsh-Weyl $m$-function,
and both of these solutions grow rapidly as $x\rightarrow\infty$.
Thus, it is convenient to rescale $u$ and its derivative to avoid
overflow in numerical computations.  We define
\begin{equation}\label{eq:yz:def}
  y(x) = xe^{-x^2/2}u(x), \qquad z(x) = \Psi(x) x^2 e^{-x^2/2} u'(x)
\end{equation}
and note that $u$ belongs to $\H$ iff $y\in L^2(0,\infty)$.  In terms
of $y$ and $z$, the ODE (\ref{eq:eval}) can be rewritten as
\begin{equation}\label{eq:y:ode}
  \frac{d\vec r}{dx} = A(x)\vec r, \qquad
  A(x) := \frac{1}{x}\begin{pmatrix} 1-x^2 & \Psi(x)^{-1} \\
    -\lambda x^2 & x^2 \end{pmatrix}, \qquad
  \vec r = \begin{pmatrix} y \\ z \end{pmatrix},
\end{equation}
which has a singularity of the first kind \cite{coddington} at $x=0$.
Formal series solutions of (\ref{eq:y:ode}) are therefore convergent,
yielding actual solutions.  In the present case, the solution may be
expanded in a Taylor series, though other problems may require the use
of more general Frobenius series or logarithmic sums to obtain a formal
solution \cite{bender,coddington}.  Since
\begin{equation*}
  \Psi(x)^{-1} = \frac{3\sqrt{\pi}}{2} + \frac{9\sqrt{\pi}}{10}x^2 +
  \frac{153\sqrt\pi}{700}x^4 + \cdots
\end{equation*}
is even, we see that $A(x)=x^{-1}[A_0 + A_2x^2 + A_4x^4+\cdots]$
where $A_0=\begin{pmatrix} 1 & 3\sqrt\pi/2 \\ 0 & 0 \end{pmatrix}$.
It is therefore natural to construct a fundamental matrix with one
column even and the other odd:
\begin{equation}\label{eq:fund:mat}
  \Phi(x) = \Big( \vec r_0(x) \,,\, \vec r_1(x) \Big) =
\begin{pmatrix} y_0(x) & y_1(x) \\ z_0(x) & z_1(x) \end{pmatrix},
  \qquad
  \left(\begin{aligned}
    \vec r_0 &= \vec c_0+\vec c_2x^2+\cdots, \\
    \vec r_1 &= \vec c_1x+\vec c_3x^3+\cdots
  \end{aligned}\right).
\end{equation}
Matching terms yields
\begin{align*}
  A_0 \vec c_0 &= 0, &
  (2kI-A_0)\vec c_{2k} &= \jt \sum_{j=1}^k A_{2j}\vec c_{2k-2j}, & \quad (k&\ge1), \\
  (I-A_0)\vec c_1 &= 0, & \quad
  ((2k+1)I-A_0)\vec c_{2k+1} &= \jt \sum_{j=1}^k A_{2j}\vec c_{2k+1-2j}, &  (k&\ge1),
\end{align*}
where $I$ is the $2\times2$ identity matrix.  Since the eigenvalues of
$A_0$ are 0 and 1, nontrivial vectors $\vec c_0$ and $\vec c_1$ exist,
which are defined up to multiplicative factors.  Once these factors
are chosen, the higher order coefficients $\vec c_{2k}$ and $\vec
c_{2k+1}$ are uniquely determined from the recursion relationships
given above.  The leading terms are
\begin{equation}\label{eq:taylor}
  \bigg( \vec r_0\,,\, \frac{\vec r_1}{x} \bigg)
  = \begin{pmatrix}
    \frac{3\sqrt{\pi}}{2} & 1 \\
    -1 & 0 \end{pmatrix}
  - \begin{pmatrix}
    \frac{9}{40}(14\sqrt{\pi} + 5\pi\lambda) & \frac{1}{4}(2+\sqrt{\pi}\lambda) \\
    \frac{1}{4}(2+3\sqrt{\pi}\lambda) & \lambda/3
    \end{pmatrix}x^2 + O(x^4).
\end{equation}
The arbitrary constants were chosen so that
\begin{equation}\label{eq:W}
  u_1(0)=1, \qquad
  W[u_0,u_1] = \Psi w(u_0u_1' - u_1u_0') = \frac{1}{x}\det\Phi(x) \equiv 1,
\end{equation}
where $u_j=w^{-1/2}y_j$, $u_{j}'=x^{-1}w^{-1/2}\Psi^{-1}z_j$,
($j=0,1$), and $W$ is the Wronskian determinant.  The general solution
$u = \alpha_0 u_0 + \alpha_1 u_1$ reduces to (\ref{eq:u:lam0}) when
$\lambda=0$.

Note that $\Phi(x)$ is analytic in a complex neighborhood of $x=0$,
i.e.~(\ref{eq:y:ode}) has only an apparent singularity
\cite{coddington} at $x=0$.  The determinant of a fundamental matrix
is always zero at an apparent singularity, which is true in our case
as $\det\Phi(x)=x$.  An alternative first-order system using $u(x)$
and $u'(x)$ as components would yield a fundamental matrix with a pole
at $x=0$.  Another alternative in which $z(x)$ is replaced by
$\Psi(x)y'(x)$ yields the equation
\begin{equation}\label{eq:y:ode:intro}
  -(\Psi y')' + V(x)y = \lambda y, \qquad
  V(x) = (1-x^2)\frac{\Psi'(x)}{x} + (x^2-3)\Psi(x),
\end{equation}
which is self-adjoint and regular at the origin, alleviating the need
for initialization with series solutions.  This advantage comes at the
cost of $V(x)$ being more expensive to evaluate than $A(x)$ in
(\ref{eq:y:ode}) due to the additional $\Psi'(x)/x$ term. We also note
that the formulas for $\Psi(x)$ and $\Psi'(x)/x$ in
(\ref{eq:y:ode:intro}) are numerically unstable near the origin, and
have to be computed with a series for small $x$ anyway.

\subsection{Green's function and the Titchmarsh-Weyl $m$-function}
\label{sec:green}

For any $\lambda\in\mathbb{C}$ with $\im\{\lambda\}\ne0$, we can
construct a Green's function for $L-\lambda$.  We seek an operator
$[\mathcal{G}(\lambda)f](x)= \int_0^\infty
G(x,\xi;\lambda)f(\xi)w(\xi)\,d\xi$ that satisfies
\begin{equation}\label{eq:GL:inv}
  \mathcal{G}(\lambda)(L-\lambda)u=u, \qquad
  (L-\lambda)\mathcal{G}(\lambda)f=f
\end{equation}
for a wide class of functions $u$ and $f$, which we characterize in detail below.  The key
to the construction is to identify the complex number $m(\lambda)$,
unique in the limit point case, for which
\begin{equation}\label{eq:weyl:m}
  \chi(x;\lambda):= 
  u_0(x;\lambda) + m(\lambda)u_1(x;\lambda) \;\;
  \text{belongs to } \H.
\end{equation}
Although $u_0(x;\lambda)=x^{-1}e^{x^2/2}y_0(x;\lambda)$ has a simple
pole at $x=0$ for all $\lambda$, it is square integrable on $(0,1)$
because $w(x)=x^2e^{-x^2}$. Thus, $m(\lambda)$ is determined by the
behavior of $u_0$ and $u_1$ as $x\rightarrow\infty$.
Indeed, when $\im\{\lambda\}>0$, one may show \cite{coddington} that
the set of complex numbers $m_b(\lambda)$ for which $\chi(x) =
u_0(x;\lambda) + m_b(\lambda)u_1(x;\lambda)$ satisfies
$\cos\beta\,\chi(b) + \sin\beta\,p(b)\,\chi'(b)=0$ forms a nested
family of circles (parametrized by $\beta$) in the upper half-plane
that converge to $m(\lambda)$ as $b\rightarrow\infty$.  The parameter
$\beta$ is independent of $b$, $m$ and $\lambda$, and represents a
general self-adjoint boundary condition that could be imposed at the
right endpoint if the domain were truncated to the finite interval
$[0,b]$. Taking $\beta=0$, we have
\begin{equation}\label{eq:m:lim}
  m(\lambda)
  = -\lim_{b\rightarrow\infty}\frac{u_0(b;\lambda)}{u_1(b;\lambda)}
  = -\lim_{b\rightarrow\infty}\frac{y_0(b;\lambda)}{y_1(b;\lambda)}.
\end{equation}
Note that this also holds for $\im\{\lambda\}<0$, the only
difference being that the nested family of circles are then in the
lower half-plane.

Once $m(\lambda)$ and $\chi(x;\lambda)$ are known, the Green's
function may be written
\begin{equation}\label{eq:green}
  G(x,\xi;\lambda) = \begin{cases}
    u_1(x;\lambda)\chi(\xi;\lambda), & x<\xi, \\
    u_1(\xi;\lambda)\chi(x;\lambda), & x>\xi.
  \end{cases}
\end{equation}
Suppressing $\lambda$ to simplify the notation, we have
\begin{equation}\label{eq:Gf}
  \mathcal{G}f(x) = \chi(x)\int_0^x u_1(\xi)f(\xi)w(\xi)\,d\xi +
  u_1(x)\int_x^\infty \chi(\xi)f(\xi)w(\xi)\,d\xi.
\end{equation}
If $f$ is continuous on $(0,\infty)$ and belongs to $\H$, it follows
from (\ref{eq:Gf}) and the Wronskian identity that $u=\mc Gf$ has two
continuous derivatives on $(0,\infty)$ and satisfies $(L-\lambda)u=f$,
where $L-\lambda$ is applied pointwise as a differential operator.  It
follows that $(L_1-\lambda)$ from \S\ref{sec:problem} has dense range
for $\im\{\lambda\}\ne0$, i.e.~$L_1$ is essentially self-adjoint.  In
more detail, let $S_1$ be the set of continuous functions $f$ with
compact support that satisfy $\int_0^\infty u_1fw\,d\xi=0$.  If $f\in
S_1$, then $\mc G(\lambda)f\in D_1$ (defined in \S\ref{sec:problem})
and $(L_1-\lambda)\mc G(\lambda)f=f$ (already shown).
Moreover, $S_1$ is dense in $\H$. In particular, the ``orthogonality'' condition
$\la u_1,f\ra=0$ does not preclude density as $u_1\not\in\H$. In fact,
if $f\in\H$ is any function for which $\la u_1,f\ra$ is finite,
then $\la u_1,f\ra$ can be adjusted to zero with a continuous,
compactly supported perturbation to $f$ of arbitrarily small
$\H$-norm. We leave the details to the reader.

We now characterize the domain $D$ of $L$.  The formula (\ref{eq:Gf})
is well-defined for $f\in\H$ (still assuming
$\im\{\lambda\}\ne0$). Its derivative exists almost everywhere, where
it equals $\chi'(x)\int_0^x u_1fw\,d\xi + u_1'(x)\int_x^\infty\chi
fw\,d\xi$, which is locally absolutely continuous. Thus, $\mc
G(\lambda)f(x)$ is actually differentiable everywhere, and
$(L-\lambda)\mc{G}(\lambda)f=f$ almost everywhere.  We claim
$\|\mc{G}(\lambda)f\|\le\|f\|/|\im\lambda|$, i.e.~$\mc G(\lambda)$ is
bounded.  This can be seen as follows. Let $f_n\rightarrow f$ in $\H$
with $f_n\in S_1$.  A standard argument \cite{richtmyer} using the
symmetry of $L_1$ shows that
$\|\mc{G}(\lambda)f_n\|\le\|f_n\|/|\im\lambda|$.  Since
$\mc{G}(\lambda)f_n$ is a Cauchy sequence, it converges to some
$g\in\H$ with $\|g\|\le\|f\|/|\im\lambda|$. It also converges
pointwise to $\mc G(\lambda)f(x)$, by (\ref{eq:Gf}). As a result, $\mc
G(\lambda)f = g$ a.e., as required.  In addition to belonging to $\H$,
$\mc G(\lambda)f(x)$ remains finite as $x\rightarrow0$.  Indeed,
the first term on the right-hand side of (\ref{eq:Gf})
approaches zero since $\chi(x)=O(x^{-1})$ and the integral is bounded
by $\big(\int_0^x
|u_1(\xi)|^2w(\xi)\,d\xi\big)^{1/2}\|f\|=O(x^{3/2})$.  The second term
approaches $u_1(0)\la f,\chi\ra$, which is finite.  By uniqueness of
the solutions of $(L-\lambda)u=f$ with $u\in\H$ and $u$ bounded near
the origin, the range $D$ of $\mc G(\lambda)$ consists precisely of
those functions $u\in\H$ with a locally absolutely continuous
derivative such that $Lu\in\H$. Moreover, (\ref{eq:GL:inv}) holds for
all $u\in D$, $f\in\H$, $\lambda\not\in\mathbb{R}$.  This set $D$,
which is independent of $\lambda$, is the domain of $L$.

\subsection{Spectral transform}
\label{sec:spec:tr}

From the general theory of singular self-adjoint eigenvalue problems
\cite{coddington, titchmarsh, stakgold, krall}, we know that
$m(\lambda)$ (and hence the Green's function) is analytic on
$\mbb{C}\setminus\mbb{R}$, with simple poles at the eigenvalues of $L$
and a branch cut across the continuous spectrum of $L$.  The imaginary
part of $m(\lambda)$ is positive for $\im\lambda>0$, so
$m(\lambda)$ is a Pick-Nevanlinna function \cite{donoghue}. Any such
function can be represented in the canonical form \cite{donoghue}
\begin{equation}\label{eq:m:diff}
  m(\lambda) = A\lambda + B + \int_{-\infty}^\infty \left( \frac{1}{s-\lambda} -
    \frac{s}{s^2+1}\right) \,d\rho(s),
  \qquad (\lambda\in\mathbb{C}\setminus\mbb{R}),
\end{equation}
where $A\ge0$, $B$ is real, and $\rho(s)$ is a non-decreasing,
real-valued function such that $\int d\rho(s)/(1+s^2)<\infty$.  Using
the fact from \S\ref{sec:green} that
$m(\lambda)=\lim_{b\rightarrow\infty} m_b(\lambda)$, one may show
\cite{bennewitz} that $A=0$.  Setting $\lambda=\tau+i\veps$, the imaginary part of
(\ref{eq:m:diff}) yields
\begin{equation}\label{eq:poisson}
  \jd\im\{m(\rlam+i\veps)\}=
  \int_{-\infty}^\infty \frac{\veps\,d\rho(s)}{
    (s-\rlam)^2+\veps^2}.
\end{equation}
An expression for $\rho$ as a function of $\tau$ can be derived from (\ref{eq:poisson}).
Integrating (\ref{eq:poisson}) from $\tau_0$ to $\tau_1>\tau_0$, one finds $\int_{-\infty}^\infty
\left[\tan^{-1}\left(\frac{s-\tau_0}{\veps}\right)-
  \tan^{-1}\left(\frac{s-\tau_1}{\veps}\right)\right]d\rho(s)$ on the right-hand side of the equation.  Taking the limit as
$\veps\rightarrow0$, this integrand approaches $\pi$ for
$\tau_0<s<\tau_1$ and 0 outside this range. It is also positive and
bounded by $12\pi(\tau_1-\tau_0)^2/[(2s-(\tau_0+\tau_1))^2+
3(\tau_1-\tau_0)^2]$ for $s\in\mathbb{R}$ and
$0<\veps<\tau_1-\tau_0$. Thus, by the dominated convergence theorem,
\begin{equation}\label{eq:rho:diff}
  \rho(\rlam_1) - \rho(\rlam_0) = \lim_{\veps\rightarrow0^+}
  \frac{1}{\pi}\int_{\rlam_0}^{\rlam_1}\im\{m(\rlam+i\veps)\}\,d\rlam
\end{equation}
at points of continuity $\tau_0$, $\tau_1$ of $\rho$.
Moreover, there is a 1-1 norm-preserving correspondence between $f\in \H$ and $\hat f\in
L^2(\mathbb{R};d\rho)$, where one function may be obtained from the other by the relations \cite{coddington}
\begin{equation}\label{eq:fhat}
  f(x) = \int_{-\infty}^\infty \hat{f}(\lambda)u_1(x;\lambda)
  \,d\rho(\lambda), \qquad
  \hat{f}(\lambda) = \int_0^\infty f(x)u_1(x;\lambda)w(x)\,dx.
\end{equation}
The second integral defines $\hat f(\lambda)$, the first
gives the inversion formula, and $\|\hat f\|=\|f\|$.
As with the Fourier transform, if $f$ or $\hat f$ belongs to
$L^2$ but not $L^1$, then the integrals in (\ref{eq:fhat}) must be
defined through a limiting procedure, e.g.~by multiplying $f$ or $\hat
f$ by the characteristic function supported on $[-k,k]$ and letting
$k\rightarrow\infty$. Equivalently, the integrals may be regarded as
improper integrals rather than Lebesgue integrals over the whole
integration domain; see \cite{coddington} for precise statements.

Since the transform pair (\ref{eq:fhat}) preserves the $L^2$ norm, the
jump discontinuities of $\rho$ are precisely the eigenvalues of $L$,
and the size of the jump at the $k$th eigenvalue is
$\|u_1(\cdot\,;\lambda_k)\|_\H^{-2}$. We note that all eigenvalues are
simple since any solution of $Lu=\lambda_ku$ that is not a multiple of
$u_1(\cdot\,;\lambda_k)$ will blow up as $x\rightarrow0^+$.
Since $e^{im(\lambda)}$ is bounded in the upper
half-plane, standard arguments \cite{simon:opuc,rudin:big} show that
\begin{equation*}
  m(\tau^+) := \lim_{\veps\rightarrow0^+} m(\tau + i\veps)
\end{equation*}
exists and is finite for a.e.~$\tau\in\mathbb{R}$. Decomposing
$d\rho = \rho'(\tau)\,d\tau + d\rho_s$ into absolutely continuous
and singular components, it follows from (\ref{eq:poisson}) that
for a.e.~$\tau\in\mathbb{R}$,
\begin{equation}\label{eq:rho:m}
  \rho'(\rlam) = \frac{1}{\pi}\im\{m(\rlam^+)\},
\end{equation}
and $d\rho_s$ is supported on the set
$\{\tau\,:\,\lim_{\veps\rightarrow0^+}|m(\tau+i\veps)|=\infty\}$; see
\cite{simon:opuc, rudin:big}.  While examples can be
constructed in which $d\rho_s$ has a singular continuous component
\cite{pearson:singular},
the usual situation \cite{weidmann,lavine} is
that $d\rho_s$ has a pure point spectrum and $\rho(s)$ is
absolutely continuous between eigenvalues.  In that case, if there are
no eigenvalues between $\tau-a$ and $\tau+a$, (\ref{eq:m:diff}) may be
written as a Cauchy integral
\begin{equation}\label{eq:m:cauchy}
  \frac{m(\lambda)}{2\pi i} = \varphi(\lambda) +
  \frac{1}{2\pi i}\int_{\tau-a}^{\tau+a}
  \frac{\rho'(s)}{s-\lambda}\,ds.
\end{equation}
The function $\varphi(\lambda)$, which includes $B$ and the remaining
portions of the integral in (\ref{eq:m:diff}), is analytic in the
upper and lower half-planes as well as in the disk $|\lambda-\tau|<a$.
We will show in Appendix~\ref{sec:analytic} that for $L$ in
(\ref{eq:PDE}), $\rho'(\lambda)$ is real-analytic (for $\lambda>0$).
Thus, if the radius of convergence of $\rho'(\lambda)$ at
$\lambda=\tau$ exceeds $a$,
then the contour from $\tau-a$ to $\tau+a$ along the real axis in
(\ref{eq:m:cauchy}) can be deformed to the semicircle $S_{\tau,a}$ in
the lower half-plane, and $m(\lambda)$ has an analytic continuation
from the upper half-plane to
$S_{\tau,a}$.  This leads to rapid convergence of polynomial
extrapolation methods from the upper half-plane to the real axis, as
we will see in Section~\ref{sec:length}.

Note that for $\lambda\in\mathbb{R}$,
$\rho'(\lambda)u_1(x;\lambda)u_1(\xi;\lambda)=\frac{1}{2\pi
  i}[G](x,\xi;\lambda)$, where $[G]$ is the jump in the Green's
function across the real $\lambda$ axis. Some authors
\cite{titchmarsh,stakgold} make use of this in deriving the transform
pair (\ref{eq:fhat}).  It is also useful to know that $\rho(\lambda)=
\lim_{b\rightarrow\infty}\rho_b(\lambda)$ at points of continuity of
$\rho$, where $\rho_b$ is the spectral density function for the
eigenvalue problem $Lu=\lambda u$ over the finite interval $(0,b)$
with appropriate boundary conditions. Each $\rho_b$ is a
right-continuous step function with arbitrary additive constant chosen
(in our case) so that $\rho_b(\lambda)=0$ for $\lambda<0$.

\subsection{The spectrum of $L$ and behavior of solutions of $Lu=\lambda u$}
\label{sec:point:spec}

The operator $L$ in (\ref{eq:PDE}) has only one eigenvalue,
$\lambda_0=0$.  The corresponding eigenvector is $u_1(x;0)=1$, which
has norm $\frac{1}{2}\pi^{1/4}$.  There are no negative eigenvalues
since $\langle Lu,u\rangle\ge0$.  In Appendix~\ref{sec:asym}, we show
that the asymptotic behavior of the general solution of
(\ref{eq:y:ode}) for $\lambda>0$ and $x\gg1$ has the form
\begin{equation}\label{eq:asym}
  \begin{aligned}
  & \jt y(x;\lambda) = C x^{3/4}\left(1+\frac{1}{8x\lambda}
  + \frac{5}{128x^2\lambda^2} + \frac{15}{1024x^3\lambda^3}\right) 
  \cos\Big\{
    \frac{\sqrt{2x\lambda}}{\lambda^2}\left[\frac{2}{5}x^2\lambda^2
      - \frac{x\lambda}{6} \right. \\
      & \jt \quad \left. - \frac{1}{16} + \frac{1}{64 x\lambda}
    + \frac{5}{3072x^2\lambda^2} + \left(\frac{7}{20480\lambda^3}
  - \frac{9\lambda}{160}\right)x^{-3}\right]
    -\theta \Big\}
  + O(x^{-11/4}),
  \end{aligned}
\end{equation}
where $C$ and $\theta$ are constants determined by the initial
conditions.  Thus, $y_1(x;\lambda)$ does not belong to
$L^2(\mathbb{R}_+;dx)$, $u_1(x;\lambda)$ does not belong to $\H$, and
there are no positive eigenvalues.  Moreover, a Green's function does
not exist for $\lambda>0$ since there is no $m(\lambda)$ for which
(\ref{eq:weyl:m}) holds; thus, the continuous spectrum includes
$(0,\infty)$. It actually equals $(0,\infty)$ since (\ref{eq:m:diff})
shows that the Green's function is analytic across the real axis in
regions where $\rho(\lambda)$ is constant, and
$\rho(\lambda)=\lim_{b\rightarrow\infty}\rho_b(\lambda)=0$ for
$\lambda<0$. Alternatively, when $\lambda<0$, an $m(\lambda)$ for
which (\ref{eq:weyl:m}) holds can be constructed explicitly since
one of the asymptotic solutions $u_\pm(x)$ in (\ref{eq:asym:u})
belongs to $\H$. Thus, $\mc{G}(\lambda)$ exists and $L-\lambda$
has a bounded inverse when $\lambda<0$.

\begin{figure}
  \begin{center}
\includegraphics[width=.495\linewidth]{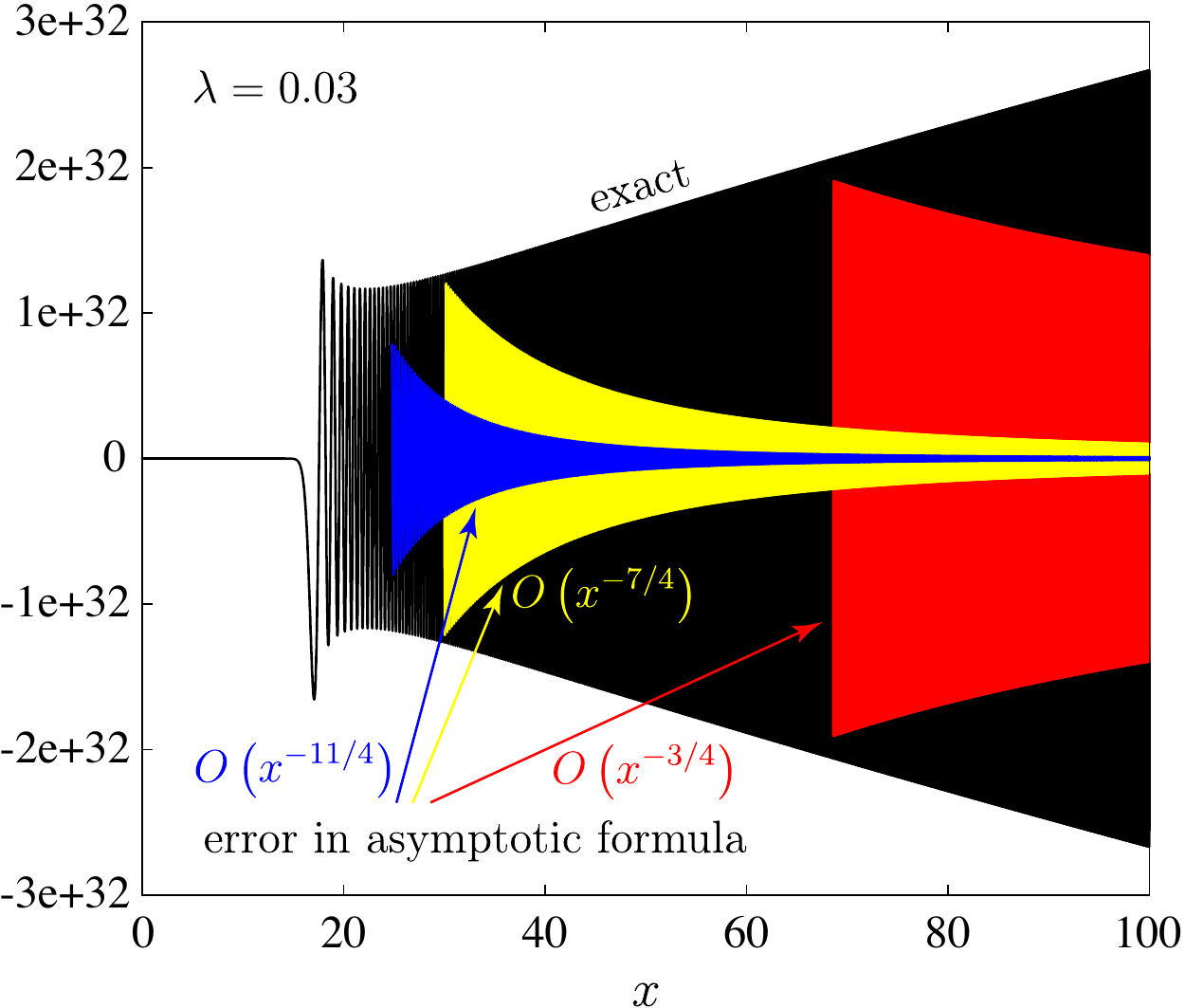}\hfill
\includegraphics[width=.495\linewidth]{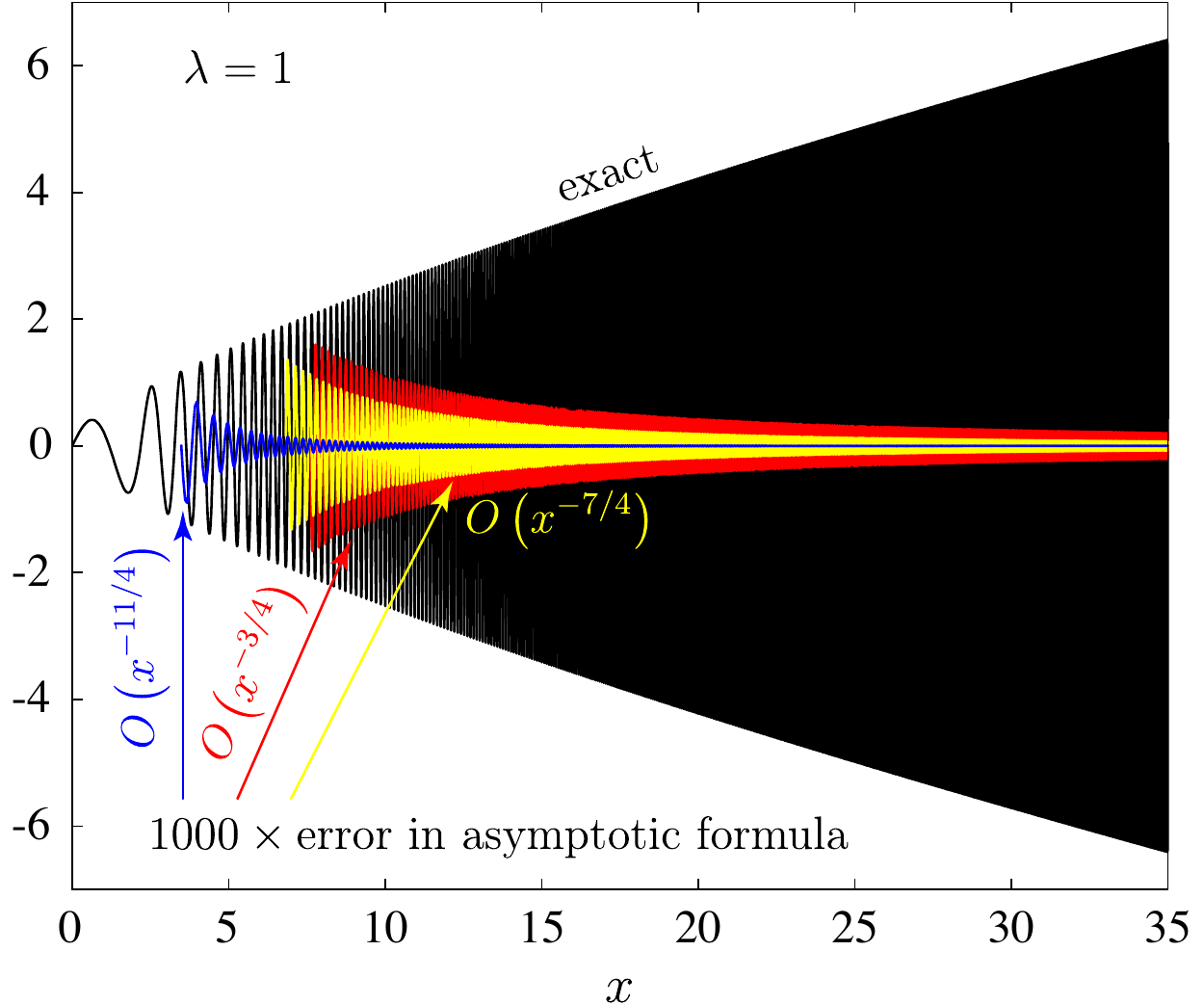}
\end{center}
\caption{\label{fig:asym} Solution $y_1(x;\lambda)$ for two values of
  $\lambda$ and the error in using the asymptotic formula
  (\ref{eq:asym}) and its lower order variants.  The ``exact''
  solution was computed using a 50th order fully implicit Runge-Kutta
  collocation (IRK) method in quadruple precision arithmetic. }
\end{figure}

Plots of $y_1(x;\lambda)$ and the error in the asymptotic
approximation (\ref{eq:asym}) are given in Figure~\ref{fig:asym} for
$\lambda=0.03$ and $\lambda=1$. For small $\lambda$, the solution
exhibits a rapid growth phase before becoming oscillatory.  For larger
values of $\lambda$, the asymptotic formula (\ref{eq:asym}) is
accurate even for small values of $x$. The three error curves
correspond to the difference between the exact solution (from solving
the ODE) and the asymptotic formula (\ref{eq:asym}), truncated at
orders $x^{-11/4}$, $x^{-7/4}$ and $x^{-3/4}$, respectively.  The
amplitude $C$ and phase $\theta$ in (\ref{eq:asym}) were obtained in
two stages.  First we computed $C_k$, $\theta_k$ by fitting the
solution through 201 data points near $x_k=175+25k$ for $\lambda=1$
and $x_k=500+80k$ for $\lambda=0.03$, with $0\le k\le 5$. Then we
extrapolated to $x=\infty$ assuming $C_k \approx C + C_\infty/x_k^{4}$,
$\theta_k \approx \theta + \theta_\infty/x_k^{3.5}$.  These values of $C$
and $\theta$ were also used for the lower order truncations.

\begin{figure}
  \begin{center}
    \includegraphics[width=.92\linewidth]{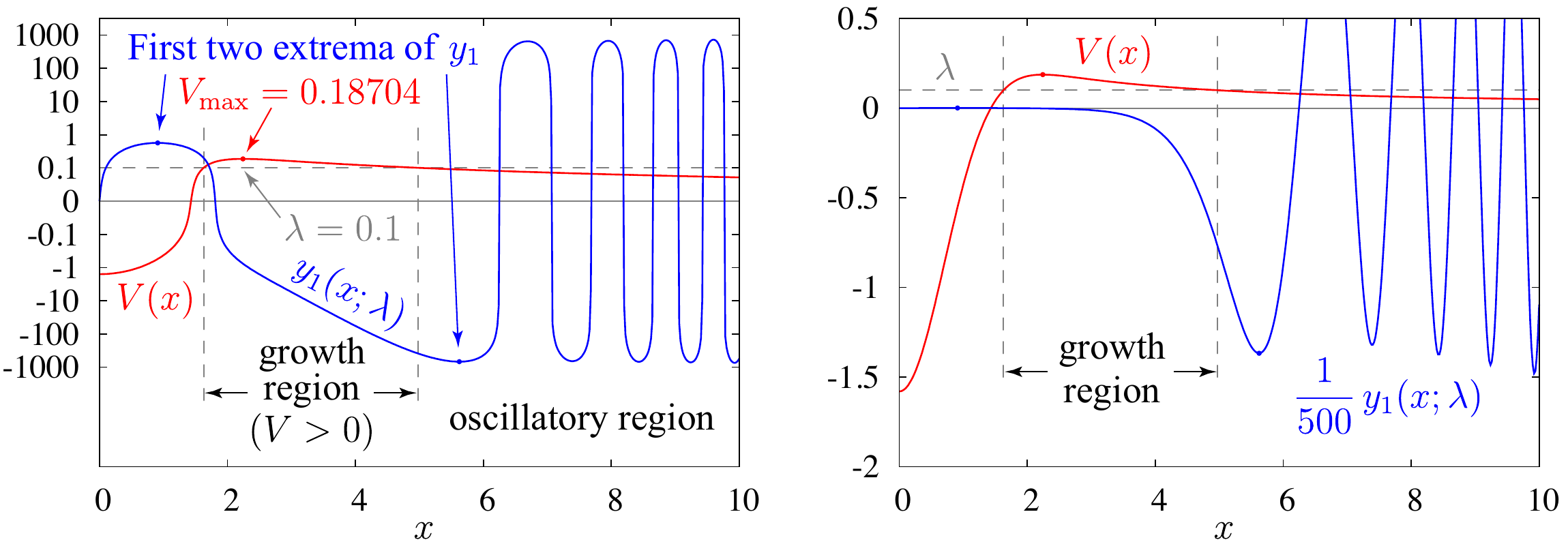}
  \end{center}
  \caption{ \label{fig:qplot} Plots of $V(x)$ and $y_1(x;\lambda)$,
    $\lambda=0.1$, on a signed log (i.e.~arcsinh) scale
    (left) and linear scale (right).  Since $\lambda<V_\text{max}$,
    there is a band of $x$ values where $V(x)>\lambda$.  In this
    region, $y_1(x;\lambda)$ grows rapidly as soon as $y_1$ and $y_1'$
    are both negative, since $y_1''=\Psi^{-1}\big[ (V-\lambda)y_1 -
    \Psi'y_1'\big]$ remains negative, accelerating the growth rate,
    until $V-\lambda$ changes sign. Decreasing $\lambda$ increases the
    size of the growth region. }
\end{figure}

The growth phase observed in Fig.~\ref{fig:asym} occurs only for $0<\lambda<0.18704$, and always
begins after the first extremum of $y_1(x;\lambda)$.  This is because
$V(x)$ in (\ref{eq:y:ode:intro}) is negative near the origin, causing
$y_1(x;\lambda)$ to execute a small half-oscillation before rapid
growth begins (see Fig.~\ref{fig:qplot}).  In more detail,
$u_1=w^{-1/2}y_1$ satisfies
\begin{equation}\label{eq:u:le:1}
  -\Psi(x) w(x) u_1'(x) = \lambda\int_0^x u_1(s) w(s)\,ds,
  \qquad u_1(0)=1.
\end{equation}
For as long as $u_1(x)$ is positive, the integral is positive and
increasing, and $u_1'(x)<0$ (assuming $\lambda>0$). Thus, by the mean
value theorem, $u_1(x)<1$ for $0<x<x_1$, where $x_1$ is the first zero
of $u_1(x)$.  This zero exists since we can use (\ref{eq:u:le:1}) to
bound $u'(x)$ away from zero for $x\in[\veps,x_1]$, where $\veps>0$ is
chosen so $u_1(x)\ge1/2$ on $[0,\veps]$.
As a result, $0<y_1(x)<xe^{-x^2/2}$ on $(0,x_1)$, so its first
extremum is of modest size. The second extremum can be much larger
in magnitude if $0<\lambda<V_\text{max}$.

\subsection{Spectral representation of the solution}
\label{sec:spec:rep}

Since the point spectrum of $L$ is $\{0\}$, the singular continuous
spectrum is absent (see Appendix~\ref{sec:analytic}), and the
absolutely continuous spectrum is $(0,\infty)$, the transform pair
(\ref{eq:fhat}) simplifies slightly, and the solution of the PDE
(\ref{eq:PDE}) with initial condition $u(x,0)=f(x)$, $f\in \H$, may be
written
\begin{align}
\label{eq:fhat2}
u(x,t) &= \frac{4}{\sqrt{\pi}}\hat{f}(0) + e^{x^2/2}
  \int_0^\infty \big[\hat{f}(\lambda)e^{-\lambda t}\big]
  \frac{y_1(x;\lambda)}{xY(\lambda)}[Y(\lambda)\rho'(\lambda)]\,d\lambda, \\
\label{eq:fhat3}
  \hat{f}(\lambda) &= \int_0^\infty \big[ xe^{-x^2/2}f(x) \big]
  y_1(x;\lambda)\,dx, \qquad (\lambda\ge0),
\end{align}
where we used $\hat u(\lambda,t)=\hat f(\lambda)e^{-\lambda t}$
and $u_1(x;\lambda)=w(x)^{-1/2}y_1(x;\lambda)$ in (\ref{eq:fhat}). Here
$Y(\lambda)$ is an arbitrary scale factor that will be used
later to account for the rapid initial growth of $y_1(x;\lambda)$ when
$\lambda$ is small.  Evaluation of $\hat{f}(\lambda)$ is
straightforward if $xe^{-x^2/2}f(x)$ decays rapidly (e.g.~if $f$ is a
polynomial). For a given $\lambda$, we just add a third component $F$
to $\vec y$ in (\ref{eq:y:ode}) and evolve $dF/dx =
xe^{-x^2/2}f(x)y_1(x;\lambda)$ along with $y_1$ and $z_1$ until
$dF/dx$ is negligible.

\section{New algorithm for computing spectral density functions}
\label{sec:comp}

The most difficult part of computing the solution $u(x,t)$
via (\ref{eq:fhat2}) and (\ref{eq:fhat3}) is determining
$\rho'(\lambda)$.  A popular approach, implemented in the software
package SLEDGE \cite{pruess1993,fulton1994,fulton1998}, is to compute
the step function $\rho_b(\lambda)$ for the regular problem over
$(0,b)$ and let $b\rightarrow\infty$. For problems of the form
$-u''+V(x)u=\lambda u$, more sophisticated methods have recently been
developed \cite{fulton2005, fulton2008b,fulton2008c} that compute
$\rho'(\lambda)$ directly, without computing auxiliary functions
$\rho_b$ and taking a limit.  Further details are given in
\S\ref{sec:compare} and \cite{wilkening:irk}.  See also
\cite{fulton2008a,fulton2013}, which focus on solving $-u''+qu=\lambda
u$ when both endpoints are singular. In addition, \cite{fulton2008a}
contains a wealth of information on the history of the problem and
provides an insightful review of the literature.  There are also many
papers on computing eigenvalues and eigenfunctions (rather than
continuous spectra and density functions) for singular Sturm-Liouville
problems; see e.g.~\cite{hammerling2010,makarov2012}.

\subsection{Chebyshev extrapolation and convergence rate}
\label{sec:cheb}

Rather than approximating $\rho(\lambda)$ by $\rho_b(\lambda)$ on a
truncated domain or searching for specialized formulas for
$\rho'(\lambda)$ that avoid solving ODEs with complex $\lambda$, we
have developed a simpler approach based directly on (\ref{eq:m:lim})
and (\ref{eq:rho:m}).  For a given $\lambda>0$, we choose a complex
line segment $\Gamma$ of length $\seglen>0$,
\begin{equation}\label{eq:Gamma:def}
  \Gamma = \lambda+i\seglen\theta, \qquad (0<\theta<1)
\end{equation}
and choose $\nn$ collocation points on $\Gamma$ consisting of nodes of a
Chebyshev-Lobatto quadrature scheme, omitting the node at $\theta=0$:
\begin{equation}\label{eq:cheb:pts}
  \lambda_k = \lambda + i\seglen\theta_k, \qquad
  \theta_k = \frac{1}{2}\left[1-\cos\left(\frac{k\pi}{\nn}\right)\right],
  \qquad 1\le k\le \nn.
\end{equation}
For each $\lambda_k$, we evaluate $m(\lambda_k)$ by computing the
limit (\ref{eq:m:lim}).  We do this by evolving $\vec
r_0(x;\lambda_k)$ and $\vec r_1(x;\lambda_k)$ simultaneously using an
arbitrary (e.g.~50th) order fully implicit Runge-Kutta collocation
(IRK) method \cite{hairer} in double or quadruple-precision
arithmetic.  As shown in Appendix~\ref{sec:asym}, there exist two
solutions of (\ref{eq:y:ode}) of the form
\begin{equation}\label{eq:y:asym:cx}
  y_\pm(x;\lambda_k) = x^{3/4} P_0(x;\lambda_k) \exp\left\{\pm i
    \frac{2\sqrt{2\lambda_k}}{5} P_1(x;\lambda_k)x^{5/2}\right\}
  \big[1 + O\big(x^{-7/2}\big)\big]
\end{equation}
for $x\gg1$, where $P_0(x;\lambda)$ and $P_1(x;\lambda)$ are defined in
(\ref{eq:pq:def}) and approach 1 as $x\rightarrow\infty$.  Since
$\im\sqrt{\lambda_k}$ is positive, $y_+(x;\lambda_k)$ decays
super-exponentially as $x\rightarrow\infty$ while all other linearly
independent solutions grow.  The function $y_1(x;\lambda_k)$ is
guaranteed to grow super-exponentially since decaying would cause the
corresponding $u_1$ to be an eigenfunction of $L$.  The same is true
of $y_0(x;\lambda_k)$ since there is a related self-adjoint boundary
value problem in which $u_0$ would then be an eigenfunction.  The
limit (\ref{eq:m:lim}) emerges when the decaying mode in
$y_0(x;\lambda_k)$ and $y_1(x;\lambda_k)$ becomes negligible in
comparison to the growing mode.

Quantitatively, (\ref{eq:y:asym:cx}) gives an asymptotic estimate for the
rate of convergence of $-y_0/y_1$ to $m$.  Expressing
$y_j(x;\lambda_k) = c_j^- y_-(x;\lambda_k) + c_j^+ y_+(x;\lambda_k)$ for
$j=0,1$, we have
\begin{equation}\label{eq:rate}
  \frac{y_0(x;\lambda_k)}{y_1(x;\lambda_k)} =
  \frac{c_0^-}{c_1^-}
    \cdot\frac{1 + (c_0^+y_+)/(c_0^-y_-)}{1 + (c_1^+y_+)/(c_1^-y_-)} =
    -m(\lambda_k)\Big[1 + o\Big(e^{-\im\{4\sqrt{2\lambda_k}/5\}x^{5/2}}\Big)\Big]
\end{equation}
for $x\gg1$.  Here we
used the fact that $c_0^-\ne0$ and $c_1^-\ne0$, as explained above,
as well as the estimate $|(1+\delta_1)/(1+\delta_2)-1|\le 2|\delta_1| + 2|\delta_2|$
when $|\delta_2|\le1/2$, and
\begin{equation}\label{eq:yp:ym:rate}
  |y_+/y_-|e^{\im\{4\sqrt{2\lambda_k}/5\}x^{5/2}}=
  \big[1+O(x^{-7/2})\big]e^{
    -\im\big\{
    (4/5)\sqrt{2\lambda_k}\big[P_1(x;\lambda_k) - 1\big]\big\}x^{5/2}},
\end{equation}
which converges to zero as $x\rightarrow\infty$ since
$\sqrt{\lambda_k}[P_1(x;\lambda_k)-1]=
\frac{5x^{3/2}}{12}\big[-\lambda_k^{-1/2}+O(x^{-1})\big]$, so its imaginary
part approaches $+\infty$ as $x\rightarrow\infty$ (since $\im\lambda_k>0$).
We conclude that the relative error in approximating $m(\lambda_k)$ by
$-y_0(x;\lambda_k)/y_1(x;\lambda_k)$ decays extremely rapidly, faster
than $\exp\big(-\im\{4\sqrt{2\lambda_k}/5\}x^{5/2}\big)$, as
$x\rightarrow\infty$.

Since this asymptotic estimate for the convergence rate only applies
for large $x$, it is useful to develop an estimate for the error that
can be monitored as the solution is evolved numerically.  We find that
the number of digits in $y_0(x;\lambda_k)/y_1(x;\lambda_k)$ that
remain frozen as $x$ increases is roughly the same as the number of
correct digits in the numerically computed Wronskian.  Thus, we use
the stopping criterion that $|W[u_0,u_1]-1|$ in (\ref{eq:W}) exceeds
1.  At this point, $\vec r_0$ and $\vec r_1$ in (\ref{eq:fund:mat})
are linearly dependent to machine precision and continuing further in
the evaluation of the limit (\ref{eq:m:lim}) does more harm than good
due to additional roundoff errors.  Note that loss of accuracy in the
Wronskian does not mean $y_0$ and $y_1$ are inaccurate; the
catastrophic cancellation of digits occurs when the determinant of
$\Phi$ is computed.

Once $m(\lambda_k)$ is known for each $\lambda_k$ in
(\ref{eq:cheb:pts}), we compute the interpolating polynomial
$q(\theta)$ satisfying
\begin{equation}\label{eq:q:interp}
  q(\theta_k) = m(\lambda_k), \qquad 1\le k\le \nn
\end{equation}
and evaluate $\frac{1}{\pi}\im\{q(0)\}$ to approximate
$\rho'(\lambda)=\frac{1}{\pi}\im\{m(\lambda^+)\}$.  The results of
this ``naive'' algorithm, and the improved version described in \S
\ref{sec:split} below, are shown in Figure~\ref{fig:rho}, where we
computed $\rho'(\lambda)$ at 768 values of $\lambda$ of the form
\begin{equation}\label{eq:grid}
  \lambda_j = e^{\sigma_j}, \qquad \sigma_j = -4 + \frac{3j}{128},
  \qquad 0\le j< 768.
\end{equation}
When $\lambda$ decreases below $0.2$, $\rho'(\lambda)$ begins to
decrease rapidly and typical values of $y_1(x;\lambda)$ grow very
large. To account for this, we introduced a scale factor $Y(\lambda)$
in (\ref{eq:fhat2}), which we define as
\begin{equation}\label{eq:Y:def}
  Y(\lambda) = \sqrt{1+y_\text{max}^2(\lambda)},
\end{equation}
where $y_\text{max}(\lambda)$ is the first negative extremum of
$y_1(x;\lambda)$, computed using Newton's method to solve
$y_1'(x;\lambda)=0$ for $x$, which occurs when
$z_1(x;\lambda)=(x^2-1)\Psi(x)y_1(x;\lambda)$.  For example, in
Figure~\ref{fig:asym}, $y_\text{max}(0.03) = -1.656\times10^{32}$
while $y_\text{max}(1) = -0.7377$.  The reason we use the first
negative extremum (rather than the positive one that precedes it) was
explained in Figure~\ref{fig:qplot}.

\begin{figure}
  \begin{center}
    \includegraphics[width=.99\linewidth]{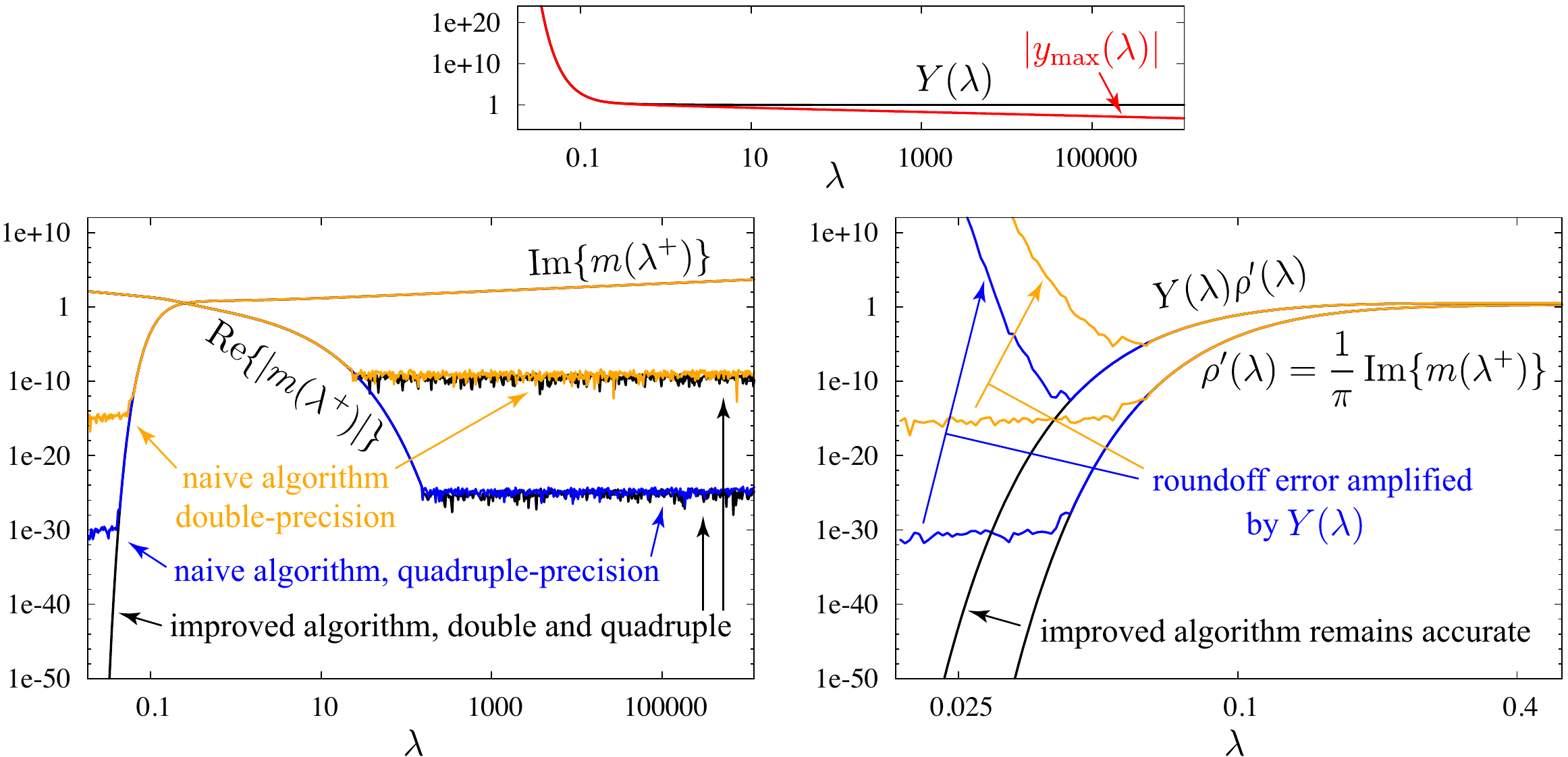}
  \end{center}
  \caption{\label{fig:rho} Plots of the real and imaginary parts of
    $m(\lambda^+)$, the spectral density function $\rho'(\lambda)$,
    and its re-scaled version, $Y(\lambda)\rho'(\lambda)$, which is
    used in (\ref{eq:fhat2}) to represent $u(x,t)$. Roundoff errors in
    the naive algorithm are amplified to unacceptable levels when
    re-scaled.}
\end{figure}

\subsection{Avoiding amplification of roundoff error}
\label{sec:split}

The poor scaling of $y_1(x;\lambda)$ poses a problem in the
reconstruction of the solution via (\ref{eq:fhat2}). Roundoff errors
in $\rho'(\lambda)$ near $10^{-15}$ in double-precision and $10^{-30}$
in quadruple-precision are amplified to large values by $Y(\lambda)$
when $\lambda$ is small enough. Fortunately, we are able to improve the
algorithm to achieve small \emph{relative} errors in $\rho'(\lambda)$.

Let us complexify $\lambda$ again and write $\lambda=\tau+i\veps$ with
$\tau>0$ fixed.  We wish to compute $\rho'(\tau)$.  First, we write
the fundamental matrix in (\ref{eq:fund:mat}) as a product,
\begin{equation}\label{eq:xstar}
  \Phi(x;\lambda) = \tilde\Phi(x;\lambda)\Phi(x^*;\lambda), \qquad
  \tilde\Phi(x^*;\lambda) = I,
\end{equation}
where $x^*$ is the location of the first negative extremum of
$y_1(x;\tau)$, which is independent of $\veps$, and
$\tilde\Phi(x;\lambda)$ is another fundamental matrix for
(\ref{eq:y:ode}), evolved from the identity at $x=x^*$.  By analogy
with (\ref{eq:fund:mat}), we denote the entries of $\tilde\Phi$ by
$\tilde y_0$, $\tilde y_1$, $\tilde z_0$, $\tilde z_1$.  When
$\veps\ne0$, the limit (\ref{eq:m:lim}) exists and we find
from (\ref{eq:xstar}) that
\begin{equation}\label{eq:mtil}
\begin{aligned}
  \tilde m(\lambda) := -\lim_{x\rightarrow\infty}
  \frac{\tilde y_0(x;\lambda)}{\tilde y_1(x;\lambda)} &=
  \lim_{x\rightarrow\infty}\frac{
    y_0(x;\lambda)z_1(x^*;\lambda) - y_1(x;\lambda)z_0(x^*;\lambda)}{
    y_0(x;\lambda)y_1(x^*;\lambda) - y_1(x;\lambda)y_0(x^*;\lambda) } \\
  &= \frac{
    m(\lambda) z_1(x^*;\lambda) + z_0(x^*;\lambda) }{
    m(\lambda) y_1(x^*;\lambda) + y_0(x^*;\lambda) }.
\end{aligned}
\end{equation}
Since $\Phi(x^*;\lambda)$ is invertible, the numerator and denominator
cannot both be zero, so the limit exists in the extended complex
plane. It is in fact finite, for $\tilde m(\lambda)$ is the
$m$-function (multiplied by $\exp[(x^*)^2]/x^*$)
 for the eigenvalue
problem $Lu=\lambda u$ on $[x^*,\infty)$ with Dirichlet boundary
conditions at the left endpoint.
Moreover, $\tilde m(\tau^+)$ exists (since $m(\tau^+)$ does),
and is given by the final formula of (\ref{eq:mtil}) with $m(\lambda)$
replaced by $m(\tau^+)$ and $\lambda$ replaced by $\tau$ in $y_0$,
$y_1$ $z_0$, $z_1$.  Next, since $\det\Phi(x^*;\tau)=x^*$, the
singular value decomposition of $\Phi(x^*;\tau)$ has the form
\begin{equation}\label{eq:Phi:star}
  \begin{pmatrix} y_0(x^*;\tau) & y_1(x^*;\tau) \\
    z_0(x^*;\tau) & z_1(x^*;\tau) \end{pmatrix} =
  \begin{pmatrix} a & -b \\ b & a \end{pmatrix}
  \begin{pmatrix} \alpha\sqrt{x^*} & \\ & \alpha^{-1}\sqrt{x^*}
  \end{pmatrix} \begin{pmatrix} c & d \\ -d & c \end{pmatrix},
\end{equation}
where $a^2+b^2=1$, $c^2+d^2=1$ and $\alpha>0$. We know that if the
first singular value is $\alpha\sqrt{x^*}$ the second singular value
must be $\alpha^{-1}\sqrt{x^*}$ because the Wronskian is equal to
1. This is convenient because the larger singular value
$\alpha\sqrt{x^*}$ can be computed accurately while the smaller one
may be severely corrupted by roundoff error. In our algorithm, we
ignore the computed version of the smaller singular value and assume
it equals $\alpha^{-1}\sqrt{x^*}$ instead. Combining
(\ref{eq:Phi:star}) with (\ref{eq:mtil}), we obtain
\begin{equation}\label{eq:m:clever}
  \begin{aligned}
    m &= -\frac{\tilde m y_0 - z_0}{\tilde m y_1 - z_1} =
    -\frac{\alpha(b-a\tilde m)c + \alpha^{-1}(a+b\tilde m)(-d)}{
      \alpha (b-a\tilde m)d + \alpha^{-1}(a+b\tilde m)(c)} \\
    & = -\beta + \frac{1+\beta^2}{\alpha^{2}\mu^{-1}+\beta}
    = \beta^{-1} - \frac{1+\beta^{-2}}{\alpha^{-2}\mu + \beta^{-1}},
    \qquad
    \left(\beta=\frac{c}{d}, \quad \mu = \frac{a+b\tilde m}{b - a\tilde m}\right),
  \end{aligned}
\end{equation}
where $m=m(\tau^+)$, $\tilde m=\tilde m(\tau^+)$, and
$y_0$, $z_0$, $y_1$ and $z_1$ are evaluated at $(x^*;\tau)$.
The resulting formula
\begin{equation}\label{eq:rho:clever}
  \rho'(\tau) = \begin{cases}
    \frac{1}{\pi}\im\left\{
    \frac{1+\beta^2}{\alpha^{2}\mu^{-1}+\beta} \right\}, & |\beta|\le 1, \\[5pt]
    -\frac{1}{\pi}\im\left\{
    \frac{1+\beta^{-2}}{\alpha^{-2}\mu+\beta^{-1}} \right\}, & |\beta|>1,
  \end{cases}
\end{equation}
avoids the cancellation of digits that occurs if $\alpha$ is large and
$m$ is not simplified in (\ref{eq:m:clever}) to separate out the
$-\beta$ (or $\beta^{-1}$) term.

The choice of $x^*$ as the first negative extremum of $y_1(x;\tau)$
ensures that $\Phi(x^*;\tau)$ captures the growth phase of
$y_1(x;\tau)$ observed in Figures~\ref{fig:asym} and~\ref{fig:qplot}.
As a result, $\tilde\Phi(x;\lambda)$ begins as the identity matrix in
the oscillatory phase of the ODE, where growth is due to $\lambda$
being complex rather than all solutions growing. This causes $\tilde
m$ in the improved algorithm to be more accurately computed than $m$
in the naive algorithm. Another advantage of the split
(\ref{eq:xstar}) is that $\lambda$ can be set to $\tau$ when computing
$\Phi(x^*;\lambda)$ since it does not depend on $x$.  Of course,
$\lambda$ must remain complex in $\tilde\Phi(x;\lambda)$ since the
limits $x\rightarrow\infty$ and $\veps\rightarrow0$ do not commute.

We remark that if $z(x)$ is replaced by $\Psi(x)y'(x)$ and $\Phi(x)$
is evolved according to
\begin{equation}\label{eq:1st:order:alt}
  \Phi'=\begin{pmatrix} 0 & 1/\Psi \\ V-\lambda & 0
  \end{pmatrix}\Phi,
\end{equation}
with $V$ as in (\ref{eq:y:ode:intro}), then $\det\Phi(x^*;\tau)=1$,
rather than $x^*$. Formulas (\ref{eq:Phi:star})--(\ref{eq:rho:clever})
remain unchanged, except that $\sqrt{x^*}$ should be omitted from the
diagonal in (\ref{eq:Phi:star}).  This technique of computing
$\rho'(\lambda)$ with high relative accuracy works generally, and is
not tied to $L$ in (\ref{eq:PDE}). Any convenient choice of $x^*$ can
be used as long it is in the ``oscillatory'' region where the solution
grows due to $\lambda$ being complex. To compute $a$, $b$, $c$, $d$,
$\alpha$ and $\tilde m$, one may use any first order system\,
$d\vec{r}/dx=A(x)\vec r$\, that is equivalent to $Lu=\lambda u$,
provided the first component of $\vec r$ is a multiple of $u$.  (In
our case, $\vec r=(y;z)$ with $y=xe^{-x^2/2}u$.)  The Abel-Liouville
formula\, $(d/dx)\det\Phi=(\opn{tr}A)\det\Phi$\, will determine the
factor to include with $\alpha$ and $\alpha^{-1}$ in
(\ref{eq:Phi:star}). Typically, as in (\ref{eq:1st:order:alt}),
$\opn{tr}A=0$ and the factor is 1.

\subsection{Error bounds, complexity, and optimal parameters}
\label{sec:length}

In the improved algorithm, we need to compute $\tilde m(\lambda^+)$
for $\lambda>0$.  This is done with an interpolating polynomial
$q(\theta)$, just as in (\ref{eq:q:interp}), but matching $\tilde
m(\lambda_k)$ instead of $m(\lambda_k)$ at $\theta_k$.  As in
(\ref{eq:cheb:pts}), $\lambda$ is now real while $\lambda_k=\lambda +
i\seglen\theta_k$, $1\le k\le \nn$, are complex.  We will show in this
section that the cost of computing $\tilde m(\lambda_k)$ for
$k=1,\dots,\nn$ is dominated by the $k=1$ term, regardless of $\nn$,
and determine the parameters $\seglen$ and $\nn$ to best take
advantage of the super-exponential rate of convergence of $\tilde
y_0(x;\lambda_k)/\tilde y_1(x;\lambda_k)$ to $-\tilde m(\lambda_k)$.

\subsubsection{Convergence rate}
Writing $\tilde y_j(x;\lambda_k) = \tilde c_j^{\,-} y_-(x;\lambda_k) +
\tilde c_j^{\,+} y_+(x;\lambda_k)$ for $j=0,1$, we obtain
(\ref{eq:rate}) again, but with tildes placed over $y_0$, $y_1$,
$c_0^\pm$, $c_1^\pm$ and $m$. The coefficients $\tilde
c_j^{\,\pm}$ are determined from the initial condition
$\tilde\Phi(x^*)=I$ by solving
\begin{equation}\label{eq:match:ypm}
  E(x^*)C =
  \begin{pmatrix}
    y_+(x^*) & y_-(x^*) \\
    z_+(x^*) & z_-(x^*)
  \end{pmatrix}
  \begin{pmatrix}
    \tilde c_0^{\,+} & \tilde c_1^{\,+} \\
    \tilde c_0^{\,-} & \tilde c_1^{\,-}
  \end{pmatrix} =
  \begin{pmatrix}
    1 & 0 \\
    0 & 1
  \end{pmatrix},
\end{equation}
where $z_\pm(x) = \Psi(x)\big[ xy_\pm'(x) +(x^2-1)y_\pm(x)\big]$ and
we have dropped $\lambda_k$ to simplify the notation, since it is fixed in this
discussion. We emphasize that $(y_\pm;z_\pm)$ are exact solutions of
the ODE (\ref{eq:y:ode}) for which $y_\pm$ have the form
(\ref{eq:y:asym:cx}) for $x\gg x^*$.  As before, the coefficients
$\tilde c_j^{\,-}$ must be non-zero since there are self-adjoint
boundary-value problems, this time on $[x^*,\infty)$, for which
$\tilde u_0$ or $\tilde u_1$ would be an eigenfunction with complex
eigenvalue if it were a multiple of $u_+=w^{-1/2}y_+$, the solution
that decays as $x\rightarrow\infty$.  In other words, neither $y_+(x)$
nor $z_+(x)$ can ever vanish since $\lambda_k\not\in\mathbb{R}$.  The
same conclusions hold if we define $z_\pm(x)=\Psi(x)y_\pm'(x)$ in
(\ref{eq:match:ypm}) and use (\ref{eq:y:ode:intro}) and
(\ref{eq:1st:order:alt}) to evolve the solution instead of
(\ref{eq:y:ode}).

The growing solution, $y_-(x)$, is not uniquely determined. Adding any
multiple of $y_+(x)$ will not affect the asymptotics of $y_-(x)$ as
$x\rightarrow\infty$. Thus, we may assume that the columns of $E(x^*)$
above are orthogonal to each other. Also, $\tilde\Phi(x)=E(x)C$ is
unchanged if we replace $E$ by $ED$ and $C$ by $D^{-1}C$, where $D$ is
diagonal. Thus, we may assume $E(x^*)$ and $C$ are orthogonal matrices
in (\ref{eq:match:ypm}) if we allow $y_+$ and $y_-$ to be of the form
(\ref{eq:y:asym:cx}) up to constant factors. The signs can be arranged
(via $D$) so that $\tilde c_0^{\,+} = -\tilde c_1^{\,-}$, $\tilde
c_1^{\,+} = \tilde c_0^{\,-}$, $y_+(x^*) = -\tilde c_1^{\,-}$, and
$y_-(x^*) = \tilde c_0^{\,-}$.  The argument of (\ref{eq:rate}) then
gives the relative error estimate
\begin{equation}\label{eq:rate2}
  \frac{1}{|\tilde m|} \left|
    \tilde m + \frac{\tilde y_0(x)}{\tilde y_1(x)}
  \right| \le 2\kappa \epsilon(x) = o\Big(
  e^{-\im\{4\sqrt{2\lambda_k}/5\}[x^{5/2}-(x^*)^{5/2}]} \Big),
\end{equation}
where
\begin{equation}\label{eq:kappa:eps}
  \kappa = 1 + \left|\frac{\tilde c_1^{\,-}}{\tilde c_0^{\,-}}\right|^2 =
  1 + |\tilde m|^{-2}, \qquad
  \epsilon(x) = \left| \frac{y_+(x)/y_+(x^*)}{y_-(x)/y_-(x^*)}\right|.
\end{equation}
The first inequality in (\ref{eq:rate2}) is valid once
$\epsilon(x)\le\frac{1}{2}$.  Note that $\epsilon(x)$ is the ratio of
the magnitudes of the decaying and growing solutions if both are
scaled to equal 1 at $x=x^*$, and if the growing solution is chosen to
be orthogonal to the decaying solution at $x=x^*$.  The $(x^*)^{5/2}$
term in (\ref{eq:rate2}) accounts for the fact that the error is
$O(1)$ when $x=x^*$. It is not strictly necessary, but avoids hiding a
large constant in the $o(\cdots)$ notation that would delay its
convergence to zero.  As mentioned already, $\tilde c_0^{\,-}$ and
$\tilde c_1^{\,-}$ are guaranteed to be nonzero, so $\kappa$ is
finite. The results of Figure~\ref{fig:cheb} below show that
$\kappa\le2$ for $4\le\lambda\le10^6$ and remains less than 3000 if
$\lambda$ is decreased to $0.01875$, assuming $\im\lambda_k$ is small
enough that $|\tilde m(\lambda_k)|$ is similar in size to $|\tilde
m(\lambda^+)|$.

\subsubsection{Step count for high-order collocation methods}

In Appendix~\ref{sec:step:bound}, we show that the number of steps
required to evolve (\ref{eq:y:ode}) or (\ref{eq:y:ode:intro}) at
$\lambda_k$ from $x_1$ to $x_2$ using a $\nu$-stage Runge-Kutta
collocation method of order $2\nu$ while maintaining a bound of
$\delta$ on the local truncation error satisfies
\begin{equation}\label{eq:N:steps}
  N_{\text{steps},k} \le K_\nu(\delta) \int_{x_1}^{x_2}
  5|\lambda_k|^{1/2}\la x\ra^{3/2} + 3.5 |\lambda_k|^{-1/2}\la x\ra^{1/2}\,dx.
\end{equation}
Here $K_\nu(\delta)=\opn{max}\big(\frac{3}{2},
\frac{4}{9}\delta^{-1/(2\nu-1)}\big)$, $\la x\ra = \sqrt{1+x^2}$, and
$\nu\ge5$ was assumed to obtain these particular constants.  For the
complexity analysis, it is convenient to decouple $x^*$ from the
location where $y_1(x)$ achieves its second extremum. Instead, we set
\begin{equation}\label{eq:xstar2}
  x^* = \begin{cases}
    1, & \lambda\ge1, \\
    \lambda^{-1}, \quad & 0<\lambda<1.
  \end{cases}
\end{equation}
This change has little effect on the running time of the algorithm,
but is easier to analyze.  Recall that we evolve to $x^*$ with
$\lambda$ real to get past the growth phase, then complexify
$\lambda_k=\lambda+i\seglen\theta_k$ and evolve from $x^*$ until
$-y_0/y_1$ converges to $\tilde m(\lambda_k)$.  As shown in
Figure~\ref{fig:qplot} above, the growth phase occurs in the band of
$x$ values for which $V(x)\ge\lambda$.
This band only exists when $\lambda<V_\text{max}=0.18704$, and ends
well before $x$ reaches $1/\lambda$, since $V(x)<1/(2x)$ for
$x>0$. Thus, (\ref{eq:xstar2}) is sufficient to traverse the growth
phase, when it exists.

Combining (\ref{eq:N:steps}) and (\ref{eq:xstar2}) and making use of
$\la x\ra^\alpha\le 1+x^\alpha$ for $x>0$, $\alpha\le2$, we find
that the number of steps required to evolve from $x^*$ to
$x$ is bounded by
\begin{equation}\label{eq:N:steps2}
  \frac{N_{\text{steps},k}}{K_\nu(\delta)|\lambda_k|^{1/2}} \, \le \,
  \frac{34}{5}\big[x^{5/2}-(x^*)^{5/2}\big],
\end{equation}
where we used $|\lambda_k|^{-1}\le1$ when $\lambda=\re\lambda_k\ge1$
and $|\lambda_k|^{-1}\le x$ when $\lambda<1$ and $x\ge
x^*=\lambda^{-1}$.  An additional
$(77/6)\max(\lambda^{-2},\lambda^{1/2})K_\nu(\delta)$ steps are
required to evolve from 0 to $x^*$.

\subsubsection{Cost of evolving all $n$ solutions from $x^*$}
We see from (\ref{eq:rate2}) that for small enough $\delta$, the
relative error in estimating $\tilde m$ by $-\tilde y_0/\tilde y_1$
will be less than $\delta$ if
\begin{equation}\label{eq:evolve:to}
  b(x) = \ln(1/\delta), \qquad\quad
  \bigg( b(x) := \im\{4\sqrt{2\lambda_k}/5\}[x^{5/2} - (x^*)^{5/2}]\bigg).
\end{equation}
Here we rely on surpassing the point beyond which $|o(e^{-b(x)})|\le
e^{-b(x)}$ in (\ref{eq:rate2}). Equivalently, $\epsilon(x)$ and
$\kappa$ in (\ref{eq:kappa:eps}) must satisfy $\epsilon(x)\le
e^{-b(x)}/(2\kappa)$ by the time $b(x)$ reaches $\ln(1/\delta)$. This
happens rapidly since $\epsilon(x)e^{b(x)}$ converges to zero
super-exponentially in spite of the rapid growth of $e^{b(x)}$, as in
(\ref{eq:yp:ym:rate}).
Combining (\ref{eq:N:steps2}) and (\ref{eq:evolve:to}), the number of
steps required to achieve a relative error of $\delta$ satisfies
\begin{equation*}
  \frac{5}{34}\cdot\frac{4\sqrt{2}}{5}\cdot
  \frac{N_{\text{steps},k}}{K_\nu(\delta)\ln(1/\delta)}\le
  \frac{\sqrt{|\lambda+i\seglen\theta_k|}}{
    \im\big\{\sqrt{\lambda+i\seglen\theta_k}\big\}} \;=\;
    2\frac{\lambda}{\seglen\theta_k} +
  \frac{3}{4}\frac{\seglen\theta_k}{\lambda} +
  O\left(\frac{\seglen\theta_k}{\lambda}\right)^3.
\end{equation*}
We see that for fixed $\delta$, $N_{\text{steps},k}$ is a function of
$\seglen\theta_k/\lambda$, and, to leading order, scales inversely
with it. Retaining only this leading term, which is particularly
accurate for those $\theta_k$ closest to zero (which matter most), we
can estimate the total cost (in steps) of computing $\tilde
m(\lambda_1)$, \dots, $\tilde m(\lambda_n)$:
\begin{equation}\label{eq:cost}
  \frac{\sqrt{2}}{17}\cdot
  \frac{N_\text{tot}}{K_\nu(\delta)\ln(1/\delta)} \; \lesssim \;
  \frac{\lambda}{\ell}\sum_{k=1}^\nn \theta_k^{-1} \;=\; 
  \frac{\lambda}{\seglen}
  \left(\frac{2\nn^2+1}{3}\right).
\end{equation}
The final equality follows from observing that the polynomial
\begin{equation*}
  P(\theta)=\theta(\theta_1-\theta)(\theta_2-\theta)\cdots(\theta_\nn-\theta)
\end{equation*}
satisfies $P'(0)=\theta_1 \theta_2\cdots\theta_\nn$ and
$P''(0)=-2P'(0)(\theta_1^{-1}+\cdots +\theta_\nn^{-1})$.  Setting
$Q(\vartheta)=P\big((1-\cos\vartheta)/2\big)$, we find (from the choice of
the $\theta_k$ in (\ref{eq:cheb:pts}) as Chebyshev-Lobatto nodes) that
$Q(\vartheta)=\sin n\vartheta\,\sin \vartheta$.  This gives
\begin{equation*}
  P'(0) = \lim_{\vartheta\rightarrow0}\frac{2Q'(\vartheta)}{\sin\vartheta}, \qquad
  P''(0) = \lim_{\vartheta\rightarrow0}\left(
    \frac{4Q''(\vartheta)}{\sin^2\vartheta} - \frac{4Q'(\vartheta)}{\sin^3\vartheta}\cos\vartheta
    \right),
\end{equation*}
which evaluate to $4\nn$ and $-\frac{8}{3}\nn(2\nn^2+1)$, respectively.  The
result (\ref{eq:cost}) follows from $-P''(0)/[2P'(0)]=(2\nn^2+1)/3$.

It is worth noting that $\theta_1^{-1}=2/[1-\cos(\pi/\nn)]$ also grows
quadratically for large $\nn$. A series expansion shows that
$d/d\nn\big([1-\cos(\pi/\nn)](2\nn^2+1)\big)>0$ for $\nn\ge1$.  Thus, the
ratio of $\theta_1^{-1}$ to $\sum_1^\nn\theta_k^{-1}$ decreases
monotonically from 1 (at $\nn=1$) to $6/\pi^2\approx0.60793$ (as
$\nn\rightarrow\infty$). This shows that for $\nn\ge1$, the cost of
computing $\tilde m(\lambda_1)$ is more than $60\%$ of the total cost of
computing $\tilde m(\lambda_k)$ at all the collocation points along the
complex line segment $\Gamma$ --- the closest point to the real
$\lambda$ axis dominates the others.

\subsubsection{Bound on the extrapolation error}
As explained in Section~\ref{sec:spec:tr}, if
$\rho(s)$ is real analytic at $s=\lambda>0$, then $m(\zeta)$ can be
analytically continued from the upper half-plane to a disk of radius
$a$ centered at $\lambda$. By (\ref{eq:mtil}), the same is true of
$\tilde m(\zeta)$.  Thus, there exists $\bern>1$ such that $\tilde
m(\zeta)$ is analytic and bounded in the region $B$ bounded by
the Bernstein ellipse \cite{trefethen:book}
\begin{equation*}
  \zeta = \lambda+i\frac{\seglen}{2}\left(1-\frac{z+z^{-1}}{2}\right), \qquad
  z=\bern e^{i\theta}, \qquad
  0\le\theta\le 2\pi,
\end{equation*}
which encloses $\Gamma$, has foci at its endpoints, and has semi-major
and semi-minor axes equal to $\seglen(\bern+\bern^{-1})/4$ and
$\seglen(\bern-\bern^{-1})/4$, respectively.  In particular, it is not
hard to show that the choice
\begin{equation}\label{eq:rho:from:a}
  \bern=\sqrt{1+\frac{a^2}{\seglen^2}}+\frac{a}{\seglen}
  + \sqrt{2} \left(
    \frac{a}{\seglen}\sqrt{1+\frac{a^2}{\seglen^2}}+\frac{a^2}{\seglen^2}
  \right)^{1/2}
\end{equation}
causes $\partial B$ to cross the real $\zeta$-axis at $\lambda\pm a$
and to remain strictly above the semicircle
$\{\zeta\,:\,|\zeta-\lambda|=a\,,\, \im \zeta<0\}$. We regard
$\seglen$ as a free parameter of the method and $a$ as a property of
the operator $L$: it is the (effective) radius of convergence of
$\rho'(s)$ at $s=\lambda$. Note that $a$, $\bern$ and $z$ here are not
related to those in (\ref{eq:Phi:star}), (\ref{eq:y:ode}) and
(\ref{eq:evolve:to}).  From (\ref{eq:rho:from:a}), we see that
$\bern\sim 4a/\seglen$ for $\seglen\ll a$ and $\bern\sim
1+\sqrt{2a/\seglen}$ for $\seglen\gg a$.
The significance of $\bern$ is that the Chebyshev
coefficients of $\tilde m(\zeta)\vert_\Gamma$ decay like $\bern^{-j}$, i.e. we can write
\begin{equation}\label{eq:cheb:mtil}
  \tilde m(\lambda+i\theta\seglen) = \sum_{j=0}^\infty a_j
  T_j(1-2\theta), \qquad  0\le\theta\le 1, \qquad
  |a_j|\le 2M\bern^{-j},
\end{equation}
where the $T_{j}$ are Chebyshev polynomials and $M=\max_{\zeta\in\overline B} |\tilde m(\zeta)|$; (see
\cite{trefethen:book}).  Let $p(\theta)$ be the polynomial of degree
$\nn$ that interpolates the value of $\tilde m(\lambda+i\theta\seglen)$
at $\{\theta_k\}_{k=0}^\nn$ (including $k=0$). Then $p(0)=\tilde
m(\lambda^+)$ and the Chebyshev coefficients of
$p(\theta)=\sum_0^\nn\hat p_jT_j(1-2\theta)$ satisfy\, $|\hat p_j -
a_j|\le 2M\bern^{-\nn}/(\bern-1)$, due to aliasing
\cite{trefethen:book}. The polynomial $q(\theta)$ of degree $\nn-1$ that
interpolates $\tilde m$ at $\{\theta_k\}_{k=1}^\nn$ differs from
$p(\theta)$ by a multiple of the Lagrange polynomial
\begin{equation*}
  l_0(\theta) \;=\; \prod_{j=1}^\nn \frac{\theta-\theta_j}{0-\theta_j} \;=\;
  \frac{T_0(1-2\theta) + T_\nn(1-2\theta)}{2\nn} +
  \sum_{j=1}^{\nn-1} \frac{T_j(1-2\theta)}{\nn}.
\end{equation*}
The multiple is chosen to reduce the degree of $p(\theta)$ by one,
i.e.~$q(\theta)=p(\theta)-2\nn\hat p_\nn l_0(\theta)$. As a result,
\begin{equation}\label{eq:q0:err}
  |q(0)-\tilde m(\lambda^+)| \;=\; 2\nn|\hat p_\nn| \;\le\;
  2\nn\big(|a_\nn|+|\hat p_\nn - a_\nn|\big) \;\le\; \frac{4\nn M\bern^{1-\nn}}{\bern-1}.
\end{equation}
In addition to this bound on the extrapolation error, we find that the
Chebyshev coefficients of $q$ agree closely with those of $\tilde
m(\zeta)\vert_\Gamma$:
\begin{equation}\label{eq:qj:err}
  |\hat q_j - a_j| \;\le\; |\hat q_j - \hat p_j| + |\hat p_j - a_j| \;\le\;
  \frac{6M\bern^{1-\nn}}{\bern-1}, \qquad\quad (0\le j<\nn),
\end{equation}
where we used $|\hat q_j - \hat p_j| = |2\nn\hat p_\nn \hat{l_0}_j|$,
which is bounded by $2\hat p_\nn$.

\begin{figure}
  \begin{center}
    \includegraphics[width=.99\linewidth]{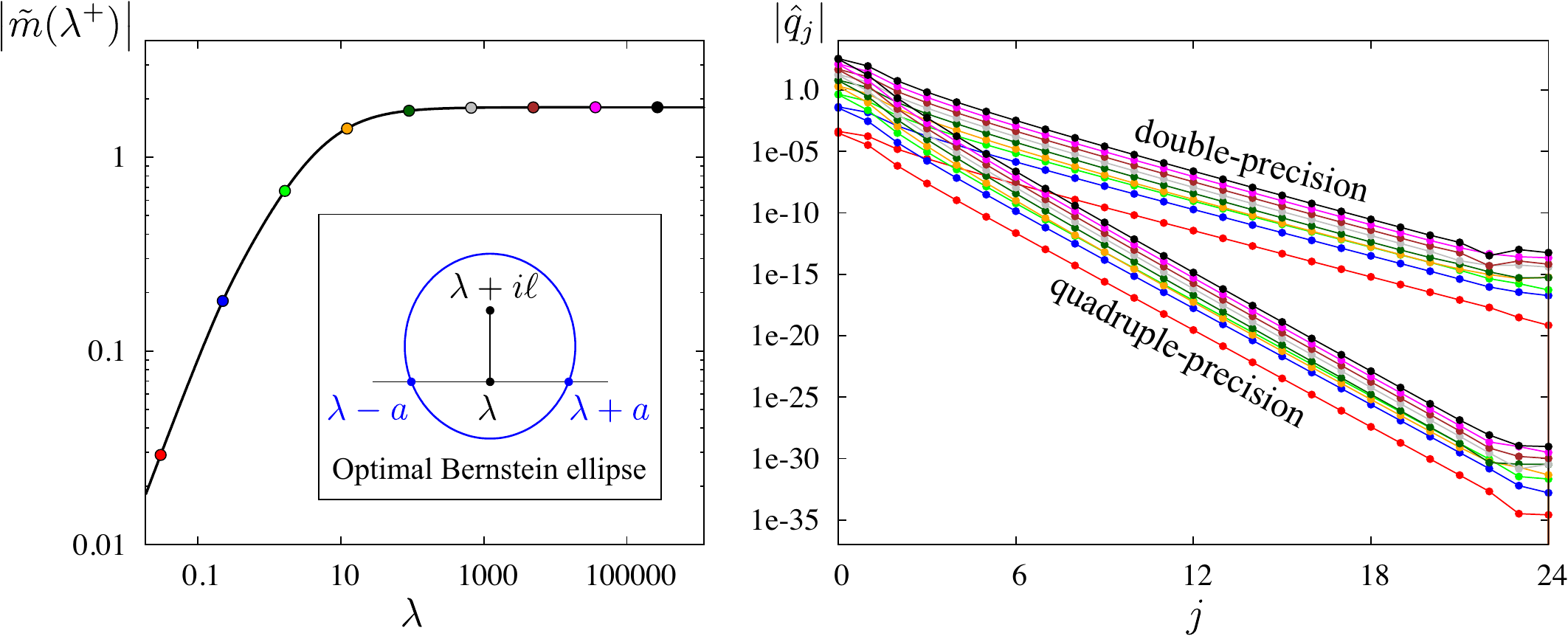}
  \end{center}
  \caption{\label{fig:cheb} Plots of $\left|\tilde
      m(\lambda^+)\right|$ and the magnitudes of the Chebyshev
    coefficients $\hat q_j$ used to compute $\tilde m(\lambda^+)$.
    The points at left are color coded to correspond to two curves at
    right, one for each type of machine precision used to compute it.
    The interval lengths $\seglen$ in (\ref{eq:sigma:opt}) were chosen
    so the Chebyshev coefficients decay just fast enough to reach
    machine precision at $j=24$. The inset at left shows the
  optimal Bernstein ellipse with $\alpha=4.9834$.}
\end{figure}

\subsubsection{Numerical computation of optimal segment length}
Figure~\ref{fig:cheb} confirms that the Chebyshev modes of $q(\theta)$
decay exponentially, as predicted from $|a_j|\le 2M\bern^{-j}$ and
(\ref{eq:qj:err}).  The function $\tilde m(\lambda^+)$ computed here
corresponds to $x^*$ as defined in \S\ref{sec:split} rather than the
simplified version in (\ref{eq:xstar2}).  We find empirically
that with $\nn=25$ collocation points, the segment lengths
\begin{equation}\label{eq:sigma:opt}
  \seglen = \frac{13\lambda}{11+\lambda^{-11/8}}
  \;\; (\text{double precision}), \quad
  \seglen = \frac{1.85\lambda}{9+\lambda^{-11/8}}
  \;\; (\text{quadruple})
\end{equation}
are close to optimal over a wide range of $\lambda$ values.  Indeed,
for each $\lambda$ in the range computed in Figure~\ref{fig:cheb}
($e^{-4}\le\lambda\le e^{14}$), the Chebyshev modes $\hat q_j$ reach
the roundoff threshold just before $j$ reaches $24=\nn-1$.  Assuming
the observed modes $|\hat q_j|$ are a good predictor of the bound
$|a_j|\le2M\bern^{-j}$, then we can estimate $2M\bern^{1-n}$ by $|\hat
q_{n-1}|$. If the latter is $O(\delta)$, then $|q(0)-\tilde
m(\lambda^+)|$ is at most $2n/(\bern-1)$ times larger, or
$O(n\delta)$, due to (\ref{eq:q0:err}). The value of $\bern$ is
substantially larger than $1$ since reducing $|a_j|$ by a factor of
$\delta$ when $j$ varies from 0 to $n-1$ requires
$\bern\gtrsim\delta^{-1/(n-1)}$.  When $n=25$, this gives
$\bern\gtrsim 4.2$ in double-precision and $\bern\gtrsim 18$ in
quadruple-precision arithmetic.  We note from (\ref{eq:rho:from:a})
that $\bern$ can be controlled easily by varying $\seglen$.  In
practice, we choose $\seglen$ to be just small enough that roundoff
errors in the highest-frequency Chebyshev modes of $\hat q_j$ begin to
be visible, hitting a plateau rather than continuing to decay
exponentially.

\subsubsection{Roundoff error} \label{sec:roundoff}
The effects of floating-point arithmetic
can be taken into account using the extrapolation formula
\begin{equation}\label{eq:extrap:formula}
  q(0)=(-1)^{\nn-1}q(\theta_\nn) - 2\sum_{k=1}^{\nn-1}(-1)^k q(\theta_k).
\end{equation}
In exact arithmetic, if $q(\theta)$ interpolates the exact values of
$\tilde m(\lambda+i\seglen\theta_k)$ for $1\le k\le \nn$, we have
shown how to choose $\seglen$ so that the right-hand side of
(\ref{eq:extrap:formula}) will equal $\tilde m(\lambda^+)$ to
$O(n\delta)$.  If each term of the right-hand side of
(\ref{eq:extrap:formula}) is perturbed by $O(\delta)$, the accumulated
effect remains $O(n\delta)$.  In practice, the errors really are this
small, as shown in the follow-up paper \cite{wilkening:irk} by comparing
double-precision results to ``exact'' solutions computed in quadruple-precision,
and indirectly in \cite{vsck2} by comparing this method of solving (\ref{eq:PDE})
to a projected dynamics approach using orthogonal polynomials.

\subsubsection{Complexity estimate}
\label{sec:complexity}

We conclude this section with a complexity estimate for the algorithm
when $\nn$ (the number of collocation points) and $\nu$ (the number of
fully implicit Runge-Kutta stages) are chosen optimally.

First, for fixed $\lambda$ and $\delta$, we wish to find $n$ and
$\seglen$ to minimize $N_\text{tot}$ in (\ref{eq:cost}) subject to the
constraint $\bern=\delta^{-1/(n-1)}$.  A nearly identical optimization
problem (with a simpler answer) is to minimize $\nn^2/\seglen$ subject
to $\bern=\delta^{-1/n}$, where $\bern$ is regarded as a function of
$\seglen$ (with $a$ fixed) in (\ref{eq:rho:from:a}). A routine
calculation shows that the optimal solution satisfies
\begin{equation*}
  \alpha(\seglen) \ln \alpha(\seglen) + 2\seglen \alpha'(\seglen) = 0, \qquad
  \nn = \frac{\ln(1/\delta)}{\ln\alpha(\seglen)},
\end{equation*}
which yields $\seglen/a = 0.9065$, $\bern = 4.9834$ and $\nn =
0.6226\ln(1/\delta)$. This gives $\nn=23$ for $\delta=10^{-16}$ and
$\nn=46$ for $\delta=10^{-32}$. Thus, our choice of $\nn=25$ above is
nearly optimal in double-precision but too small in
quadruple-precision, which is what we observe in practice as well (see
\S\ref{sec:compare} below). Note that the eccentricity of the optimal
Bernstein ellipse equals $2/(\bern+\bern^{-1})=0.3858$, independent of
$\lambda$ and $\delta$ (see Figure~\ref{fig:cheb}).  Its size is
determined by $a$, the effective radius of convergence of $\rho'(s)$ at
$s=\lambda$.  We can estimate $a$ as a function of $\lambda$
indirectly, using $\seglen/a = 0.9065$ and the left formula in
(\ref{eq:sigma:opt}), since $\nn=25$ is close to optimal in
double-precision.  The target accuracy $\delta$ affects the number of
grid points via $\nn=0.6226\ln(1/\delta)$, but not the size or shape
of the ellipse.  With these parameter choices, (still approximating
$2n^2+1$ by $2n^2$), we obtain
\begin{equation}\label{eq:cost2}
  N_\text{tot} \le 3.107\, K_\nu(\delta)
  \ln^3(1/\delta)\,\frac{\lambda}{\seglen},
\end{equation}
where $\lambda/\seglen\approx 11/13 + (1/13)\lambda^{-11/8}$ in our
case, and $K_\nu(\delta)=\opn{max}\big(\frac{3}{2},
\frac{4}{9}\delta^{-1/(2\nu-1)}\big)$.

The additional $(77/6)\max(\lambda^{-2},\lambda^{1/2})K_\nu(\delta)$
steps required to evolve from 0 to $x^*$ are normally a small fraction
of $N_\text{tot}$. These initial steps do not require complex
arithmetic, which further reduces their cost.  Nevertheless,
$\lambda^{-2}$ will eventually dominate $\lambda^{-11/8}$ when
$\lambda\rightarrow0$. Over the range of $\lambda$ and tolerances
considered here ($e^{-4}\le\lambda\le e^{14}$, $\delta\le 10^{-15}$),
$3.107\ln(1/\delta)^3$ is at least 9 times larger than
$(77/6)\lambda^{1/2}$. But for still larger values of $\lambda$, the
latter could become significant. Fortunately, this term is an artifact
of choosing $x^*=1$ in (\ref{eq:xstar2}), and is not present when
$x^*$ is defined as the first negative extremum of
$y_1(x;\lambda)$. The naive algorithm (in which $x^*$ is set to zero
and only one fundamental matrix is computed) can even be used since
there is no growth phase when $\lambda\ge V_\text{max}=0.18704$.

The difficulty in analyzing the case in which $\lambda$ is large and $x^*$
is allowed to drop below 1 is that the estimate (\ref{eq:rate2}) of
the convergence rate is based on the large $x$ asymptotics of
$\Psi(x)$, yielding errors of the form
$\exp(-\gamma[x^{5/2}-(x^*)^{5/2}])$, where $\gamma$ is a multiple of
$\im\sqrt{\lambda_k}$.  This formula (incorrectly) suggests that if
$x^*$ is close to zero, little progress will be made until $x$ reaches
1. By contrast, (\ref{eq:N:steps}) does not rely on asymptotics, and
shows that the number of steps needed to evolve from $x^*$ to $x$
scales linearly, like $|\lambda_k|^{1/2}[x-x^*]$, when $x^*$ is close
to zero and $x\le 1$.  We expect that a more refined WKB analysis
using (\ref{eq:xi:def}) instead of (\ref{eq:xi:asym}) would show that
for large $\lambda$, $|\tilde y_0/\tilde y_1+\tilde m|$ decays
initially like $\exp(-\tilde\gamma[x-x^*])$, where $\tilde \gamma$ is
also a multiple of $\im\sqrt{\lambda_k}$, before exhibiting the
$\exp(-\gamma[x^{5/2}-(x^*)^{5/2}])$ behavior. As a basic check, if we
approximate $\sqrt{F(x)/\Psi(x)}$ in (\ref{eq:xi:def}) by
$\sqrt{F(0)/\Psi(0)}$ and assume that $h_\pm(x)$ in (\ref{eq:vpm})
remain small down to $x=0$ when $\lambda\gg1$, then $\epsilon(x)$ in
(\ref{eq:kappa:eps}) becomes $\exp(-\tilde\gamma[x-x^*])$ with
$\tilde\gamma = \sqrt{6}\pi^{1/4} \im\{\sqrt{\lambda_k}\} +
O(|\lambda_k|^{-1/2})$.  Since the number of steps and the logarithm
of the inverse error grow similarly as functions of $x$, with
prefactors proportional to $|\lambda_k|^{1/2}$ and
$\im\{\sqrt{\lambda_k}\}$, respectively, the technique of linking
$N_\text{steps}$ to $\ln(1/\delta)$ by comparing (\ref{eq:N:steps2})
to (\ref{eq:evolve:to}) should work the same. We did not carry out the
details as this was not an issue over the range of $\lambda$
considered here.

Next, to leading order in $\nu$, the computational cost of the steps
in (\ref{eq:cost2}) is $C_\text{tot}= C(\delta)\times
\frac{1}{3}(2\nu)^3\times8\times N_\text{tot}$, where $C(\delta)$ is
the cost of a floating-point operation with roundoff threshold
$\delta$, $\frac{1}{3}(2\nu)^3$ is the number of multiplications
required to solve the $2\nu\times2\nu$ linear system associated with a
$\nu$-stage fully-implicit Runge-Kutta step for a linear ODE with two
components (see \cite{wilkening:irk}), and $8$ accounts for the four
multiplications and four additions required to perform one complex
multiplication and one complex addition. Optimizing $C_\text{tot}$
boils down to minimizing $\nu^3\delta^{-1/(2\nu-1)}$ with $\delta$
fixed. Minimizing $(2\nu-1)^3\delta^{-1/(2\nu-1)}$ instead gives
$2\nu-1=\frac{1}{3}\ln(1/\delta)$, for which
$K_\nu(\delta)=\frac{4}{9}\exp(3)\approx8.927$. Using $2\nu\le
\frac{10}{9}(2\nu-1)$ for $\nu\ge5$, we obtain
\begin{equation}\label{eq:ctot}
  C_\text{tot} \le 3.76\,C(\delta) \ln^6(1/\delta)\frac{\lambda}{\seglen}.
\end{equation}
Since $C(\delta)$ grows like $\ln^2(1/\delta)$ or $\ln(1/\delta)\ln[
\ln(1/\delta)]$, depending on the arbitrary precision implementation,
we conclude that the cost of computing $\rho'(\lambda)$ by this method
with accuracy $\delta$ grows slower than $\delta^{-\gamma}$ for any
$\gamma>0$ as $\delta\rightarrow0$.  This translates into significant
performance gains over other methods for computing spectral density
functions with high accuracy, as shown in \S\ref{sec:compare} and
\cite{wilkening:irk}.

We find that the estimate $n=0.6226\ln(1/\delta)$ for the optimal
number of collocation points agrees closely with comparisons of actual
running times in practice --- the same $n$ is close to
optimal. However, the bounds (\ref{eq:N:steps2}) and (\ref{eq:cost2})
overpredict the number of steps taken using adaptive stepsize control
\cite{wilkening:irk} by a factor of roughly $3(2\nu-1)$, and the
optimal choice of $(2\nu-1)$ is a few times larger than
$\frac{1}{3}\ln(1/\delta)$. This factor of $(2\nu-1)$ can be obtained
in the analysis of Appendix~\ref{sec:step:bound} if we set
$rM=(2\nu-1)$ instead of $5/3$ in (\ref{eq:trunc:error}).  ($M$ is a
type of Lipschitz constant for the ODE, $r$ is the radius of a disk
centered at $x$ in the complex plane, and $h$ is the
stepsize). However, for technical reasons explained in
\cite{wilkening:irk}, our current analysis requires $hM\le 2/3$, which
breaks down if $r$ is increased in this way. We believe it should be
possible to remove this barrier using A-stability of the scheme rather
than a Neumann series to bound the condition number of the implicit
Runge-Kutta equations, but we do not know how to do this. If it is
indeed possible to increase $rM$ to $(2\nu-1)$ without losing control
of the constants in the formulas, then $K_\nu(\delta)$ would be of the
form $\frac{\text{const}}{2\nu-1} \delta^{-1/(2\nu-1)}$ and the result
(\ref{eq:ctot}) could be improved to contain $\ln^5(1/\delta)$ instead
of $\ln^6(1/\delta)$.

We remark that these optimization problems are only intended to serve
as a guideline for choosing $n$ and $\nu$ and deriving a rough
complexity estimate. The correctness of the algorithm does not depend
on choosing $n$ or $\nu$ optimally, and does not rely on asymptotics.

\section{Numerical examples}
\label{sec:num}

We now consider two examples illustrating the use of (\ref{eq:fhat2}),
(\ref{eq:fhat3}) to solve $u_t=-Lu$ with initial conditions
$u(x,0)=f_j(x)$, namely
\begin{equation}\label{eq:ex12}
  \text{Example 1:}\quad f_1(x)=x, \qquad\quad
  \text{Example 2:}\quad f_2(x)=x^2.
\end{equation}
Example 1 is harder to compute since $Lf_1$ has a singularity at $x=0$
that leads to an infinite initial speed $u_t$ there. This causes $\hat
f_1(\lambda)$ to decay slowly (like $\lambda^{-2}$), just as the
Fourier transform of a function with a slope-discontinuity decays
slowly.  Nevertheless, this singular example is relevant to the problem
of resistive damping in a plasma; see Section~3 of \cite{landreman1}.

\begin{figure}
  \begin{center}
    \includegraphics[width=.99\linewidth]{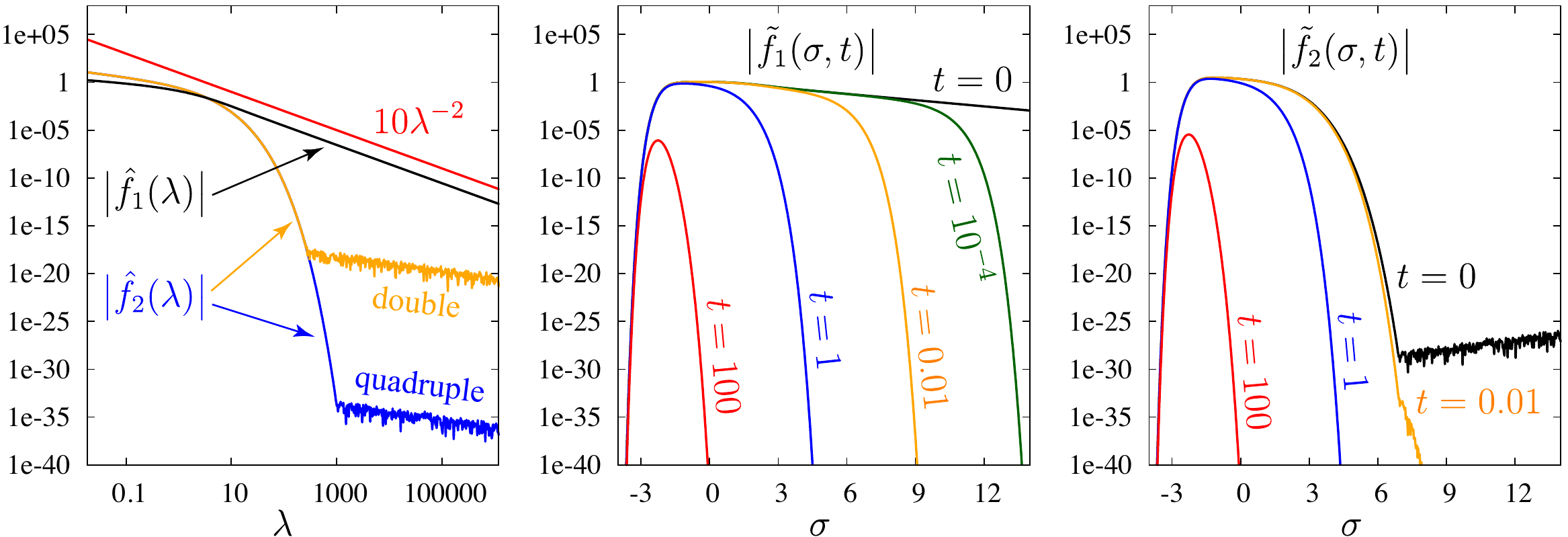}
  \end{center}
  \caption{\label{fig:fhat} Plots of $\big|\hat{f}(\lambda)\big|$ and
    $\big|\tilde f(\sigma,t)\big|$ at various times for the two
    examples in (\ref{eq:ex12}).  The effective support (exceeding
    $10^{-30}$) of $\tilde f_2(\sigma,t)$ lies between $-4<\sigma<7$
    for $t\ge0$ while that of $\tilde f_1(\sigma,t)$ extends beyond
    $\sigma=14$ until $t=10^{-4}$ due to the slow $\lambda^{-2}$ decay
    rate of $\hat f_1(\lambda)$. Both $\tilde f_1(\sigma,t)$ and
    $\tilde f_2(\sigma,t)$ turn out to be negative everywhere they can
    be distinguished from roundoff error.  }
\end{figure}

We computed $\hat f_1(\lambda)$ and $\hat f_2(\lambda)$ at
the grid points $\lambda_j$ in (\ref{eq:grid}) using the method
explained above, in which a third component is added to $\vec r$ in
(\ref{eq:y:ode}) to represent $F(x)=\int_0^x
se^{-s^2/2}f(s)y_1(s;\lambda)\,ds$, and the solution is evolved until
$F(x)$ reaches its limiting value.  The results are shown in
Figure~\ref{fig:fhat}, where we have adopted the notation
\begin{equation}\label{eq:ftil}
  \begin{aligned}
  \tilde f(\sigma,t) &:= \hat f(e^\sigma)\exp(-e^\sigma t)Y(e^\sigma)\rho'(e^\sigma)e^\sigma, \\
  \tilde u(x,t) &:= \big[ u(x,t) - 4\hat{f}(0)/\sqrt{\pi} \big]e^{-x^2/2},
  \end{aligned}
  \quad
  \tilde u(x,t) = \int_{-\infty}^\infty
  \frac{y_1(x;e^\sigma)}{xY(e^\sigma)}\tilde f(\sigma,t)\,d\sigma.
\end{equation}
The tilde here is not related to the one used in
(\ref{eq:xstar})--(\ref{eq:m:clever}) to denote solutions of
(\ref{eq:y:ode}) starting at $x=x^*$.  Note that
$\tilde f(\sigma,t) = \hat u(\lambda,t)Y(\lambda)\rho'(\lambda)\lambda$
with $\hat u(\lambda,t)=\hat f(\lambda)e^{-\lambda t}$ and
$\lambda=e^\sigma$, and the extra factor of $\lambda$ accounts
for $d\lambda = e^\sigma d\sigma$ in
the change of variables.  Also, if $\tilde f(\sigma,t^*)$
can be represented efficiently as a function of $\sigma$ for some
fixed $t^*\ge0$, then for $t\ge t^*$ we have
\begin{equation}\label{eq:ftil:eval}
  \tilde f(\sigma,t)=\tilde f(\sigma,t^*)e^{-e^\sigma(t-t^*)},
\end{equation}
which is easy to evaluate.
We will represent $\tilde f(\sigma,t^*)$ using a Fourier series on
$-4=\sigma_L \le \sigma \le \sigma_R = 14$, where $t^*=10^{-4}$ in
Example 1 and $t^*=0$ in Example 2.

\begin{figure}
  \begin{center}
    \includegraphics[width=\linewidth]{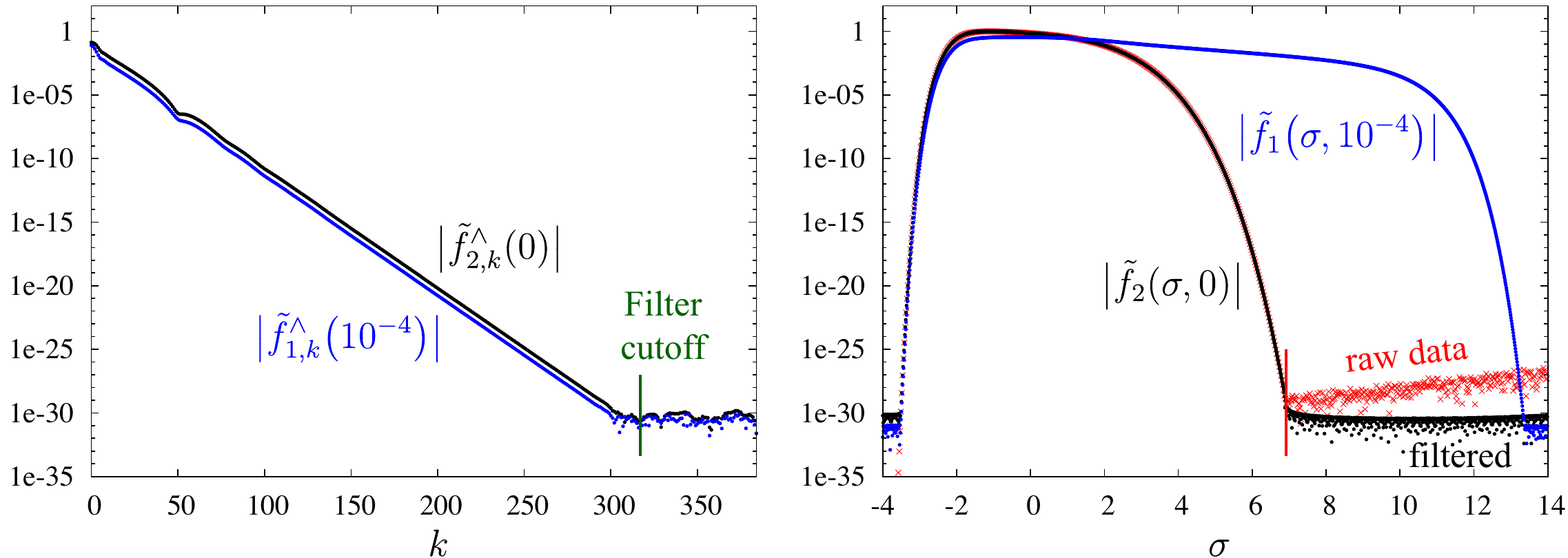}
  \end{center}
  \caption{\label{fig:filter} Sampling $\tilde f_1(\sigma,0.0001)$ and
    $\tilde f_2(\sigma,0)$ at the 768 grid points in (\ref{eq:grid})
    is sufficient to reduce the highest frequency Fourier modes to
    roundoff error in quadruple-precision arithmetic. (right) The raw
    data (red) is evaluated on the original grid while the filtered
    data (black and blue) is evaluated on a finer mesh with 1536 grid
    points.}
\end{figure}

Figure~\ref{fig:filter} shows that the grid spacing in (\ref{eq:grid})
is sufficient to represent $\tilde f_1(\sigma,t)$ over $-4\le\sigma\le
14$ for $t\ge10^{-4}$ and $\tilde f_2(\sigma,t)$ for $t\ge0$ to
quadruple-precision accuracy using the FFT. The red $\times$ markers
were obtained by computing $\hat f_2(\lambda)$ and $\rho'(\lambda)$
directly. The black markers were obtained from the red by truncating
the data at $\sigma=6.921875$ (the vertical red line), applying the
FFT (to all 768 points), truncating the Fourier series at $k=317$
(the vertical green line), and transforming back.  The blue markers
show a filtered version of the raw data labeled $t=10^{-4}$ in
Figure~\ref{fig:fhat}. In this case, the Fourier series was truncated
at $k=317$ with no initial filter in $\sigma$.  Note that roundoff
error causes $|\tilde f(\sigma,t)|$ to grow to around $10^{-30}$ near
$\sigma=-4$ in both examples.  This is not a problem since
$Y(\lambda)$ was included in the raw data before the FFT was computed.
For larger values of $t$ in Example 1 and all $t\ge0$ in Example 2,
the domain can be reduced to $-4\le\sigma\le8$ so that only 512 raw
data points are needed.  However, we will continue to work with the
grid (\ref{eq:grid}) for illustration.

To perform the integral in (\ref{eq:ftil}), we use the trapezoidal
rule over the interval $-4\le\sigma\le 14$. We increase the number of
collocation points as $x$ increases in order to resolve the
increasingly oscillatory integrals involved.  Much of this work can be
done once and for all, independent of the initial condition $f(x)$.
To this end, we pre-compute
\begin{equation}\label{eq:g:def}
  g(x;\lambda) = \frac{y_1(x;\lambda)}{xY(\lambda)}, \qquad\quad
  \lambda = e^\sigma
\end{equation}
at selected $x$ locations on a nested hierarchy of grids
\begin{equation}\label{eq:nested}
  \sigma^{(p)}_j = -4 + \frac{3j}{128\times 2^{p}}, \qquad
  0\le j < 768\times 2^{p}, \qquad 0\le p\le 10.
\end{equation}
For $p\ge1$, $g(x,\exp(\sigma^{(p)}_j))$ only has to be computed for
odd indices $j$ since it is already known for even indices from the
previous level.  Note that $\sigma^{(0)}_j$ coincides with $\sigma_j$
in (\ref{eq:grid}).  This makes it easy to interpolate the values of
$\tilde f(\sigma,t)$ rather than computing $m(\lambda^+)$ and $\hat
f(\lambda)$ at the new grid points directly. To increase the size of
$\tilde f$ by a factor of $2^{p}$, we simply zero-pad its FFT and
compute the inverse transform. For example, the black and blue markers
in the right panel of Figure~\ref{fig:filter} were computed in
this way on the $p=1$ grid with 1536 nodes.

\begin{figure}[t]
  \begin{center}
    \includegraphics[width=\linewidth]{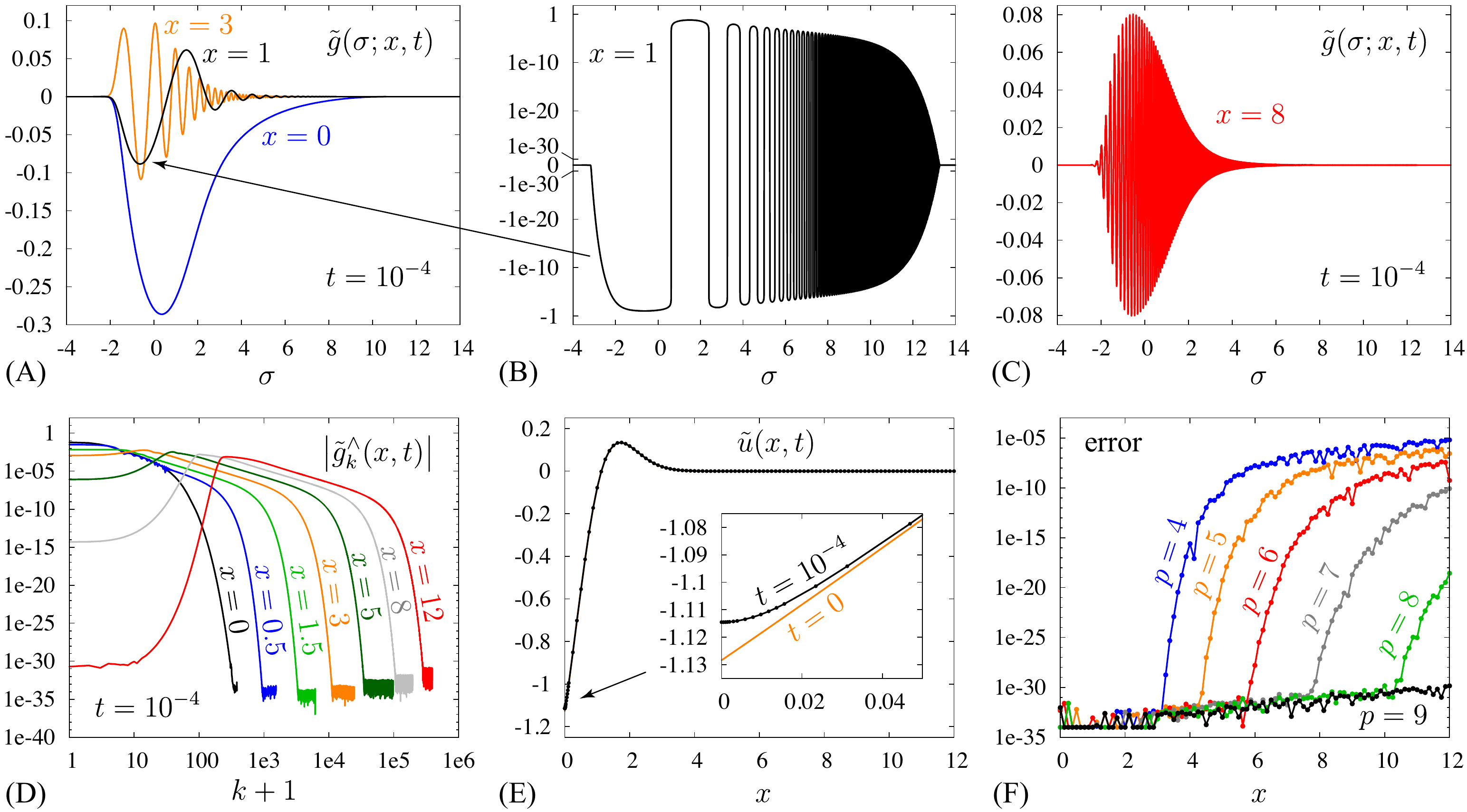}
  \end{center}
  \caption{\label{fig:fourier1} Plots of $\tilde g(\sigma;x,t)$, its
    Fourier modes $\tilde g^{\wedge}_k(x,t)$, the solution $\tilde
    u(x,t)$, and the error $|\tilde u^{(p)}(x,t) - \tilde
    u^{(10)}(x,t)|$ for Example 1 at $t=10^{-4}$. As $x$ increases,
    $\tilde g(\sigma;x,t)$ becomes more oscillatory and more grid
    points are needed to avoid aliasing errors in the trapezoidal rule
    integration scheme.
  }
\end{figure}

The results of this calculation for Example 1 are given in
Figure~\ref{fig:fourier1}.  Panels (A)--(C) show the integrand
$\tilde g(\sigma;x,t)$ in the reconstruction formula 
\begin{equation}\label{eq:g:til}
  \tilde u(x,t) = \int_{-\infty}^\infty \tilde g(\sigma;x,t)\,d\sigma, \qquad\quad
  \tilde g(\sigma;x,t) = g(x,e^\sigma)\tilde f(\sigma,t),
\end{equation}
as a function of $\sigma$ for $x=0$, $x=1$, $x=3$ and $x=8$ with
$t=0.0001$ fixed.  Note that as $x$ increases, $\tilde g$ becomes more
oscillatory as a function of $\sigma$.  Panel (B) shows the $x=1$
solution stretched vertically to a signed logarithmic scale.  This was
done by plotting $\opn{arcsinh}\big(\frac{1}{2}10^{32}\tilde g\big)$
on the $y$-axis and placing tick marks where $\tilde g=\pm 10^{-10k}$.
Note that $\tilde g$ becomes highly oscillatory as it decays.  On a
stretched scale, solutions with other values of $x$ have a similar
envelope to the $x=1$ solution shown here, and range from having no
oscillations ($x=0$) to very rapid oscillations ($x=12$).  Panel (D)
gives the magnitudes of the Fourier modes of $\tilde g(\sigma;x,t)$
with $x$ and $t$ held fixed.  The mode amplitudes of the FFT are
normalized by $1/N$, where $N=768\times2^{p}$ is the number of grid
points.  With this scaling, the $N$-point trapezoidal rule gives
$\tilde u(x,t)=18\tilde g^{\wedge}_0(x,t)$, where
$18=\sigma_\text{max}-\sigma_\text{min}$.  The curves labeled $x=0$,
$x=0.5$, etc., were computed with $p=0$, $2$, $4$, $6$, $8$, $9$ and
$10$, respectively.  These levels were chosen so that $\tilde
g^{\wedge}_k$ decays to roundoff error before $k$ reaches the Nyquist
frequency $k=N/2$, which is the largest mode shown for each curve.
Panel (E) shows the solution $\tilde u(x,t)$ at $t=10^{-4}$, obtained
by integrating $\tilde g(\sigma;x,t)$.  For comparison, we also plot
\begin{equation}
  \tilde u(x,0) = \left( x - \frac{2}{\sqrt{\pi}}\right)e^{-x^2/2},
\end{equation}
which agrees closely with $\tilde u(x,10^{-4})$ except near $x=0$,
where $u_t$ is initially infinite.  Panel (F) gives the error in the
reconstructed solution using $p=10$ as the exact solution.  Higher
values of $x$ require finer grids to resolve the oscillations in
$\tilde g(\sigma;x,t)$.  For this example, the $p=9$ and 10 solutions
are identical to 30 digits of accuracy.

\begin{figure}[t]
  \begin{center}
    \includegraphics[width=\linewidth]{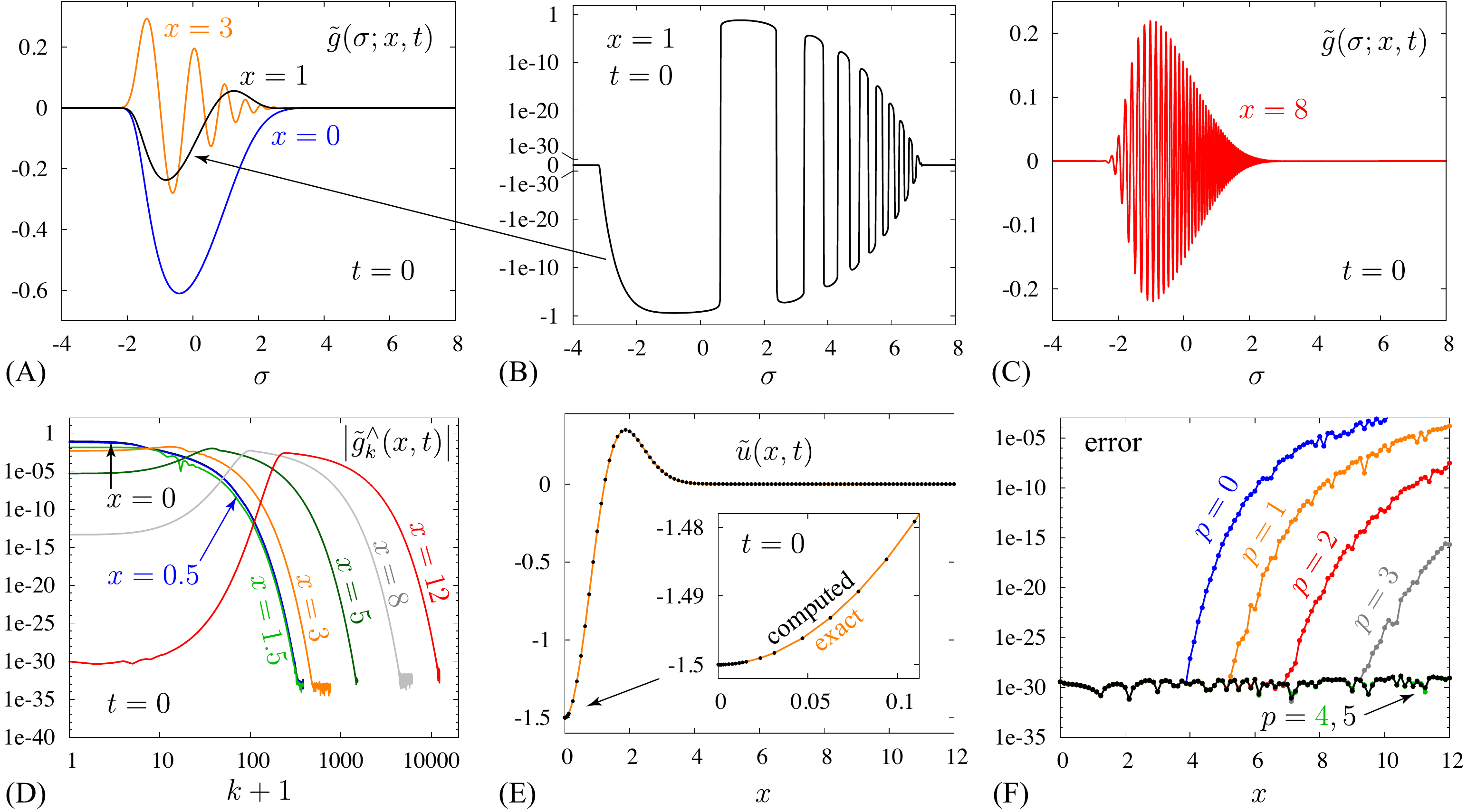}
  \end{center}
  \caption{\label{fig:fourier2} Analogous results for Example 2. Note
    that the axes are scaled differently in
    Figure~\ref{fig:fourier1}. }
\end{figure}

In Figure~\ref{fig:fourier2}, we present analogous results for
Example~2.  Since $L$ maps $f(x)=x^2$ to a smooth function that is
well-behaved at the origin, the
reconstruction can be done at $t=0$ to recover $u(x,0)=f(x)$.  By
contrast, we needed $t\ge10^{-4}$ to overcome the slow decay of $\hat
f(\lambda)$ in Example~1.  The main difference between
Figures~\ref{fig:fourier1} and~\ref{fig:fourier2} is that the
effective support of $\tilde g$ is smaller in the latter case, leading
to less oscillatory integrals. This is evident on comparing panel (B)
in both figures, and also on observing that fewer Fourier modes are
needed in panel (D) to reach machine-precision in Example 2.  More
specifically, the curves labeled $x=0$, $x=0.5$, etc., in (D) were
computed with $p=0$, $0$, $0$, $1$, $2$, $4$ and $5$,
respectively.  In panel (E), the orange curve gives the exact initial
condition
\begin{equation}\label{eq:exact}
  \tilde u(x,0) = \left(x - \frac{3}{2}\right)e^{-x^2/2}
\end{equation}
while the black markers are computed using the trapezoidal rule on the
numerically computed $\tilde g(\sigma;x,t)$.  Panel (F) gives the true
errors relative to the exact initial condition (\ref{eq:exact}) rather
than using the solution on the finest grid as the benchmark.  In
Example 2, the solution reaches roundoff error already at $p=4$, which
corresponds to a grid 32 times coarser than the $p=9$ mesh needed in
Example 1.

The accelerating frequency of oscillation that occurs in $\tilde
g(\sigma;x,t)$ as $\sigma$ increases is partly due to our choice of
$\lambda=e^\sigma$ in the change of variables (\ref{eq:ftil}).  When
$x$ is small, this choice is very good for representing $\tilde
g(\sigma;x,t)$ with a small number of Fourier modes. However, we can
do better for larger $x$.  From the asymptotic analysis in
Appendix~\ref{sec:asym}, we expect $x^{-1}y_1(x;\lambda)$ in
(\ref{eq:g:def}) to oscillate like
$x^{-1/4}\cos(\sqrt{8}x^{5/2}\sqrt{\lambda})$ at leading order.  Thus,
to achieve a nearly constant number of grid points per cycle with
respect to $\lambda$ holding $x$ fixed, we should change variables so
that $\lambda\sim\xi^2$ for large $\xi$.  We also want
$\lambda\rightarrow0$ as $\xi\rightarrow-\infty$.  We tried functions
of the form
\begin{equation}\label{eq:chgvar}
  \lambda = A(\sqrt{1+\xi^2}+\xi)^2 = A(\sqrt{1+\xi^2}-\xi)^{-2}
\end{equation}
and found that $A=25$ works nicely. The first formula is used for
positive $\xi$ and the second for negative $\xi$.

\begin{figure}[t]
  \begin{center}
    \includegraphics[width=\linewidth]{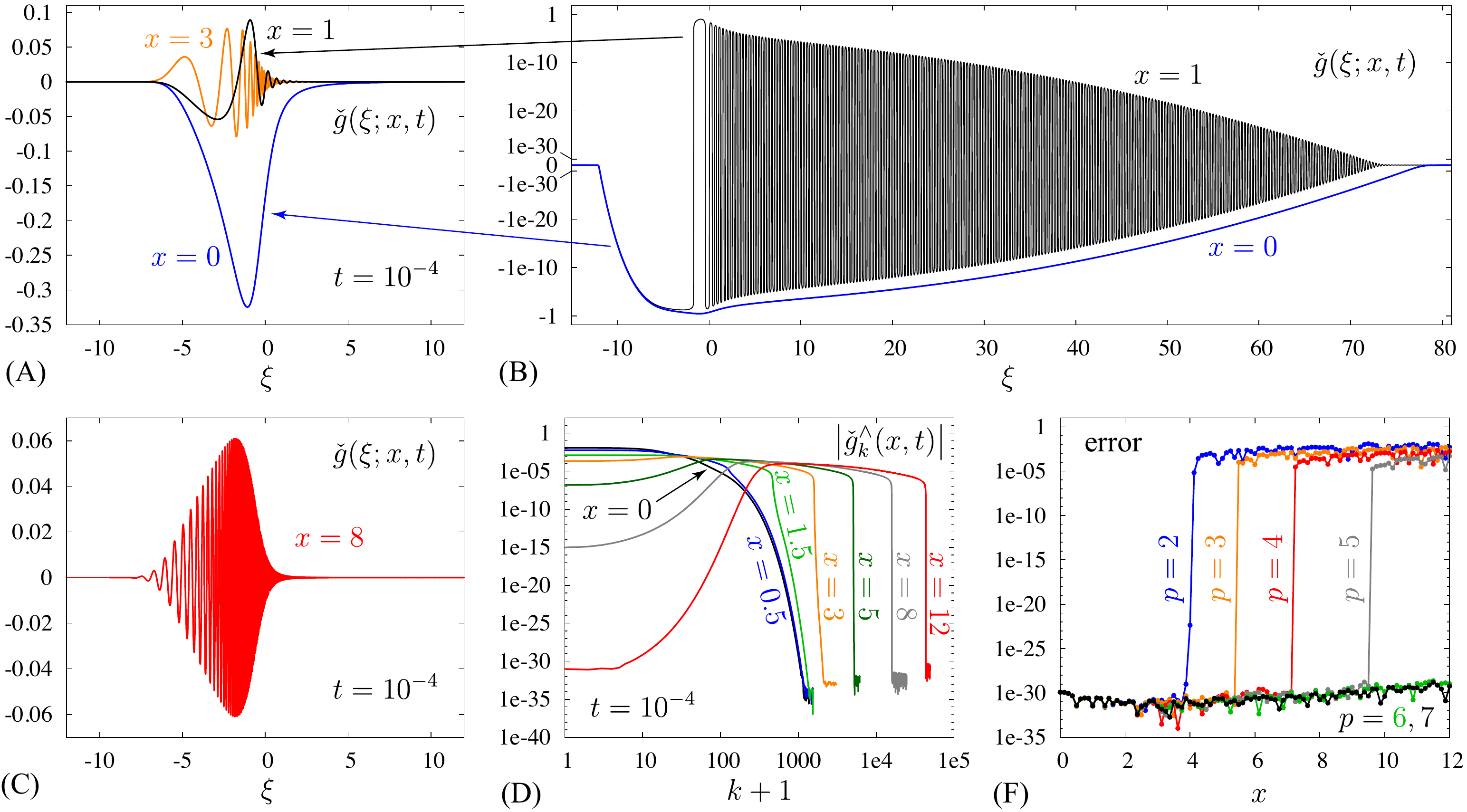}
  \end{center}
  \caption{\label{fig:fourier3} The alternative change of variables
    (\ref{eq:chgvar}) leads to a more efficient representation of the
    oscillatory integrand $\check g(\xi;x,t)$ for $x\ge5/8$.
    Panel (E) is omitted as it looks identical to
    Figure~\ref{fig:fourier1}(E).  }
\end{figure}

Figure~\ref{fig:fourier3} shows the results for Example 1 at
$t=0.0001$ with the alternative integration variable.  Plots (A)--(C)
show the integrand $\check g(\xi;x,t)$ in the reconstruction formula
\begin{equation}\label{eq:g:check}
  \tilde u(x,t) = \int \check g(\xi;x,t)\,d\xi, \qquad
  \check g(\xi;x,t) = g(x;\lambda)\check f(\xi,t),
\end{equation}
where $\lambda$ is related to $\xi$ via (\ref{eq:chgvar}),
$\check f(\xi,t) = 2\tilde f(\sigma,t)/\sqrt{1+\xi^2}$,
and $\sigma=\ln\lambda$.
Note that the oscillation frequency is nearly uniform over
$0\le\xi\le81$ in (B), unlike the result in Figure~\ref{fig:fourier1}.
To compute the integrals, we pre-compute $g(x;\lambda)$ on a nested
grid similar to (\ref{eq:nested}) but over $-15\le\xi\le81$, namely
\begin{equation}\label{eq:nested:xi}
  \xi^{(p)}_j = -15 + \frac{j}{8\times 2^{p}}, \qquad
  0\le j < 768\times 2^{p}, \qquad 2\le p\le 7.
\end{equation}
The lowest level is labeled 2 so that the number of grid points at a
given level is the same in
Figures~\ref{fig:fourier1}--\ref{fig:fourier3}.  The drawback of using
$\xi$ is that more grid points are needed at $x=0$ to represent
$\check g(\xi;x,t)$ than $\tilde g(\sigma;x,t)$.  The benefit is that
fewer grid points are needed for larger $x$.  The following table
gives the index $N$ at which the Fourier modes of $\tilde
g(\sigma;x,t)$ and $\check g(\xi;x,t)$ reach roundoff error in
quadruple-precision with $t=0.0001$
\begin{equation*}
  \begin{array}{c||c|c|c|c|c|c}
    x        &    0 &  5/8 &    1 &     3 &      8 &     12 \\ \hline
    N_\sigma &  320 & 1200 & 2000 & 11420 & 109000 & 290000 \\
    N_\xi    & 1160 & 1220 & 1300 &  2120 &  16000 &  44000
  \end{array}
\end{equation*}
We use these numbers as a guideline for the optimal number of
collocation points to use in the trapezoidal rule.  As $x$ increases,
the cost of reconstructing the solution via (\ref{eq:g:til}) or
(\ref{eq:g:check}) increases
due to more collocation points being needed to resolve the
oscillations in $\tilde g$ or $\check g$, and more timesteps being
needed to evolve to $x$.  Note that the $\sigma$ variable is better
for small $x$ because $\tilde g$ and $\check g$ grow at similar rates
near $\sigma=-2$ and $\xi=-5$, respectively, but the domain for $\xi$
is several times larger than for $\sigma$. However, for larger $x$,
the oscillations dominate the smoothness properties of $\tilde g$ and
$\check g$, and are spread out more uniformly in the $\xi$
variable. Thus, fewer collocation points are wasted in less
oscillatory regions. In our code, we use the Fourier representation of
$\tilde f(\sigma,t)$ in Figure~\ref{fig:filter} to evaluate $\check
f(\xi,t)$ on the $p=2$ grid since this only requires 768 evaluations
of $m(\lambda^+)$ and $\hat f(\lambda)$ rather than 3072.  Once
$\check f(\xi,t)$ is known at $p=2$, we use the FFT to interpolate to
higher levels.

\section{Summary of the algorithm and comparison with other methods}
\label{sec:compare}

We give here a brief summary of the method presented in detail in the
sections above, which may help readers interested in implementing
it. The key steps of the algorithm may be summarized as
follows:

\begin{enumerate}[leftmargin=.2in,rightmargin=.2in,topsep=.05in,itemsep=.05in]

\item \label{stepMESH} Pick a ``coarse'' mesh in $\lambda$-space,
  e.g.~$\lambda_j=e^{\sigma_j}$ or $25[(1+\xi_j^2)^{1/2} -
  \xi_j]^{-2}$, where $\sigma_j$ (or $\xi_j$) is uniformly spaced
  from $\sigma_L$ to $\sigma_R$ (or $\xi_L$ to $\xi_R$). The number of
  grid points and choice of $\sigma_L$ and $\sigma_R$ may be adjusted
  later, in steps~\ref{stepFTIL} and~\ref{stepFFT}.

\item \label{stepINIT} For each $\lambda$ on the mesh, evolve
  $y_1(x;\lambda)$ forward to $x^*$, the location of its first
  negative extremum, $y_\text{max}(\lambda)$. If the ODE is singular
  at the origin, i.e.~(\ref{eq:y:ode}) is used rather than
  (\ref{eq:y:ode:intro}), use the series solution (\ref{eq:taylor}) to
  initialize the ODE to the right of the singularity, e.g.~at
  $x=10^{-6}$.  Use enough terms of the series to achieve
  roundoff-level accuracy.  Also define $Y(\lambda)$ as in
  (\ref{eq:Y:def}) and evolve $y_0(x;\lambda)$ to $x^*$ to obtain
  $\Phi(x^*;\lambda)$ in (\ref{eq:xstar}). We use an arbitrary
  (e.g.~50th) order, fully implicit Runge-Kutta collocation (IRK)
  method to advance the solutions in $x$.

\item \label{stepCX} Let $\lambda_k=\lambda+i\seglen\theta_k$, with
  $\lambda$ a grid point from steps 1 and 2 and
  $\theta_1$,\dots,$\theta_\nn$ the Chebyshev
  points in (\ref{eq:cheb:pts}).  Evolve $\tilde \Phi(x;\lambda_k)$
  forward from the identity at $x=x^*$ until its columns become
  linearly dependent to machine precision (determined by monitoring
  the Wronskian). Record $\tilde m(\lambda_k) = - \tilde
  y_0(x;\lambda_k)/\tilde y_1(x;\lambda_k)$.
  % FIXME   , which stops changing once $x$ is large enough.
  Extrapolate to obtain $\tilde m(\lambda^+)$ as
  well as $m(\lambda^+)$ and $\rho'(\lambda)$.  Adjust $\seglen$, if
  necessary, to obtain appropriate Chebyshev mode decay rates
  (Fig.~\ref{fig:cheb}).

\item Compute the transform of the initial condition at the grid
  points via $\hat f(\lambda)=\lim_{x\rightarrow\infty} F(x)$,
  $F(x)=\int_0^x [w^{1/2}(s)f(s)]y_1(s;\lambda)\,ds$. $F$ is evolved
  via an ODE, along with $y_1$ and $z_1$ in (\ref{eq:yz:def}), until
  $F$ stops changing, which happens rapidly due to
  $w^{1/2}(x)=xe^{-x^2/2}$.  If the ODE is singular at the origin,
  initialize $F$ at $x=10^{-6}$ using the series solution of the ODE.

\item \label{stepFTIL}
  Evaluate $\tilde f(\sigma,t^*) = \hat f(e^\sigma)\exp(-e^\sigma
  t^*)Y(e^\sigma)\rho'(e^\sigma)e^\sigma$ on the grid at $t^*=0$.  If
  necessary, go back to step 1 and adjust the mesh endpoints so that
  $\tilde f$ decays rapidly to zero as $\sigma\rightarrow \sigma_L^+$
  and $\sigma\rightarrow\sigma_R^-$.  It may be necessary to increase
  $t^*$ to achieve sufficient decay at the right endpoint
  (Fig.~\ref{fig:fhat}).

\item \label{stepFFT} Compute the FFT of $\tilde f(\sigma,t^*)$ for
  the purpose of interpolation.  Adjust the number of mesh points (in
  step 1) as necessary so the Fourier modes of $\tilde f(\sigma,t^*)$
  decay to roundoff accuracy.  Optionally, filter the modes
  (Fig.~\ref{fig:filter}).

\item Compute the basis functions $g(x;\lambda) =
  y_1(x;\lambda)/[xY(\lambda)]$ by solving the ODE for $y_1$ on a
  nested hierarchy of grids, $\lambda=\exp\big\{\sigma_j^{(p)}\big\}$,
  as in (\ref{eq:nested}).  Record $g(x;\lambda)$ at the
  points $x$ where $u(x,t)$ is to be evaluated.  Start with level
  $p=0$ and add levels as needed in step \ref{stepTRAP}.

\item \label{stepTRAP} Compute $\tilde
  g(\sigma;x,t)=g(x,e^\sigma)\tilde f(\sigma,t)$ on successive levels
  of the grid hierarchy and evaluate\,
  $u(x,t)=\int_{\sigma_L}^{\sigma_R} \tilde g(\sigma;x,t)\,d\sigma$\,
  via the trapezoidal rule.  Since the trapezoidal rule is the
  constant mode of the FFT, stop at level $p$ of the grid hierarchy
  when the FFT of $\tilde g(\sigma;x,t)$ decays to roundoff accuracy.
  Evaluation of $\tilde f(\sigma,t)$ is done via (\ref{eq:ftil:eval}),
  with $\tilde f(\sigma,t^*)$ evaluated via the FFT from step
  \ref{stepFFT}.

\end{enumerate}
Steps \ref{stepFTIL}--\ref{stepTRAP} should be modified via
(\ref{eq:g:check}) and (\ref{eq:nested:xi}) if $\xi$ is used instead
of $\sigma$ in step~\ref{stepMESH}. We note that the procedure is very
general: $L$ can be replaced by any singular Sturm-Liouville operator
that is regular (or of limit circle type, with Taylor or Frobenius
series solutions) at the origin and of limit point type at $\infty$,
provided the spectral density function $\rho'(\lambda)$ is smooth
enough. Smoothness of $\rho$ affects the smoothness of $m(\lambda)$ in
(\ref{eq:m:diff}) and (\ref{eq:m:cauchy}), which determines how
effectively $m(\lambda)$ can be extrapolated to the real axis from the
upper half-plane via Chebyshev polynomials. Steps
\ref{stepFTIL}--\ref{stepTRAP} also rely on $\rho$ being sufficiently
smooth to ensure that the Fourier modes of $\tilde f$ and $\tilde g$
decay rapidly.  We note that each step of the algorithm provides
\emph{a-posteriori} error estimates based on the decay of Chebyshev
and Fourier modes. This allows rapid selection of mesh parameters
through a few iterations of steps~\ref{stepMESH}, \ref{stepFTIL}
and~\ref{stepFFT}.

Since the new algorithm is based on an exact mathematical formula
expressing the solution $u(x,t)$ at a later time in terms of the
initial condition $f(x)$, it will be as accurate as the quantities
$\hat f(\lambda)$, $y_1(x;\lambda)/[xY(\lambda)]$, and
$[Y(\lambda)\rho'(\lambda)]$ that appear in (\ref{eq:fhat2}), up to
quadrature error in the trapezoidal rule, which is controlled by
checking that the Fourier modes decay to machine precision.  Computing
$\hat f(\lambda)$, $y_1(x;\lambda)$ and $Y(\lambda)$ involves solving
linear ODEs, which are easily solved to machine precision using
high-order Runge-Kutta methods. Achieving high relative accuracy in
$y_1(x;\lambda)$ requires that it not be a multiple of the recessive
solution, but this is guaranteed for $\lambda\ne0$ since the point
spectrum of $L$ is $\{0\}$.  Computing
$m(\lambda)=-\lim_{x\rightarrow\infty}\frac{y_0(x;\lambda)}{y_1(x;\lambda)}$
can also be done to machine precision using high-order Runge-Kutta
methods. The extrapolation procedure to compute
$m(\tau^+)=\lim_{\veps\rightarrow0^+}m(\tau+i\veps)$ is numerically
stable in finite precision arithmetic, as shown in
\S\ref{sec:roundoff}. However, high relative accuracy of the complex
number $m(\tau^+)$ does not imply high relative accuracy of the
imaginary part, and digits can be lost when evaluating $\rho'(\tau) =
\frac{1}{\pi}\im\{m(\tau^+)\}$ when $\tau\ll1$. We showed how to avoid
this loss of accuracy in the imaginary part by factoring the
fundamental matrix as in (\ref{eq:xstar}) and using the Wronskian
identity to correct the smaller singular value of one of the factors.

To confirm that the new algorithm is spectrally accurate, we have checked in joint work with
Landreman \cite{vsck2} that the solution $u(x,t)$ computed as above agrees to
roundoff accuracy (14 or 29 digits) with the projected dynamics of
(\ref{eq:U:intro}) in spaces of orthogonal polynomials. Such agreement
provides strong evidence that monitoring Chebyshev and Fourier mode
decay rates provides accurate \emph{a-posteriori} error estimates.
The current algorithm may be viewed as a method of approximating
exact integral formulas for the solution,  while that in \cite{vsck2}
may be regarded as a (nearly) exact evolution of a finite-dimensional
approximation of the PDE.

To gauge the performance of the new algorithm,
in follow-up work \cite{wilkening:irk}, we
compare our method of computing the spectral density function
(steps~\ref{stepINIT} and~\ref{stepCX}) to the popular software
package SLEDGE \cite{pruess1993,fulton1994,fulton1998}, and to the
algorithm of Fulton, Pearson and Pruess
\cite{fulton2008c,fulton2008b}.  The algorithm in SLEDGE is based on
the Levitan-Levinson formula, $\rho(\lambda)=
\lim_{b\rightarrow\infty}\rho_b(\lambda)$, where $\rho_b$ is the
spectral function associated with the regular problem on the interval
$0<x<b$.  SLEDGE could not handle the weight function
$w(x)=x^2e^{-x^2}$ due to underflow for $x>26$ in double-precision. To
convert to a constant weight function, we made the change of variables
$u=w^{-1/2}y$ to obtain (\ref{eq:y:ode:intro}).  In this form, the
problem is regular at the origin, as discussed in
\S\ref{sec:rescaled}, and we were able to compute the spectral
function $\rho(\lambda)$ using SLEDGE to around 4 digits of
accuracy. Further refinement of the mesh led to failure of the
algorithm, apparently due to overflow when scaling the eigenfunctions
with small eigenvalues to satisfy $\Psi y_1'=1$ at $x=0$.  
Computing $\rho'(\lambda)$ from $\rho(\lambda)$ would lead to
additional loss of accuracy. Thus, SLEDGE was not found to be suitable
for our purposes.

The algorithm of Fulton, Pearson and Pruess
\cite{fulton2008c,fulton2008b} (FPP) proceeds by defining a sequence of
auxiliary functions $P_n$, $Z_n$, $R_n$ such that
\begin{equation}\label{eq:fpp:fn}
  f_n(s;\lambda) = \frac{1/\pi}{
    P_n(s;\lambda)U_1(s;\lambda)^2 +
    Z_n(s;\lambda)U_1(s;\lambda)U_1'(s;\lambda) +
    R_n(s;\lambda)U_1'(s;\lambda)^2}
\end{equation}
converges to $\rho'(\lambda)$ as $s\rightarrow\infty$, with improved
convergence rates as $n$ increases.  Here $U_1(s;\lambda)$ is related
to $u_1(x;\lambda)$ or $y_1(x;\lambda)$ in previous sections by the
Liouville transformation \cite{birkhoff,wilkening:irk} to
Schr\"odinger form, $-U''+Q(s)U=\lambda U$.  For $n=1$ and $n=2$, the
FPP procedure gives
\begin{gather}\label{eq:PZR}
  P_1=\sqrt\lambda, \qquad Z_1 = 0, \qquad
  R_1=\lambda^{-1/2}, \\
  P_2 = \sqrt{\lambda-Q}, \qquad
  Z_2 = -\frac{1}{2}Q'/(\lambda-Q)^{3/2}, \qquad
  R_2 = 1/\sqrt{\lambda-Q}.
\end{gather}
The formulas for $n=3$ are given in \cite{fulton2008c} and
\cite{wilkening:irk} along with explicit error estimates.  These
estimates imply that, in the present case, the error in approximating
$\rho'(\lambda)$ by $f_n(s;\lambda)$ decays like $s^{-(10n-8)/5}$,
where the Liouville transformation relates $x$ to $s$ via the ODE
$dx/ds=\sqrt{\Psi(x)}$; see Lemma~\ref{lem:x}, Appendix~\ref{sec:asym}
and \cite{wilkening:irk}.

The following table gives running times (in seconds) for
computing the spectral density function with our method and the FPP
method at the 768 grid points $\lambda_j = e^{\sigma_j}$ in
(\ref{eq:grid}) using a 3.33 GHz Intel Xeon X5680 system with 12
cores:
\begin{equation*}
   \begin{array}{r|c|c|c||c|c|c}
     & \text{FPP4d3} & \text{FPP15d3} & \text{WC50d25} & 
    \text{FPP31q3} & \text{WC76q25} & \text{WC76q48} \\ \hline
    \text{time} & 10637 & 9.75 & 16.2 &  21284 & 2021 & 1279
  \end{array}
\end{equation*}
Here the letters d and q stand for double- and quadruple-precision
computations with tolerances set to $10^{-15}$ and $10^{-30}$,
respectively. FPP4d3 employs the 4th order timestepper described in
\cite{fulton2008c}, which uses Richardson extrapolation to improve the
accuracy of a 2nd order frozen coefficient method, with $n=3$ in
(\ref{eq:fpp:fn}). FPP15d3 and FPP31q3 also use $n=3$ in
(\ref{eq:fpp:fn}), but with a 15th or 31st order timestepper. The letters
WC refer to the new algorithm presented in this article, and WC76q25 means
that we use a 38-stage, 76th order fully
implicit Runge-Kutta method with 25 Chebyshev extrapolation points.
Comparing WC76q25 to WC76q48 confirms the prediction in \S\ref{sec:complexity}
that $n=46$ extrapolation points should be close to optimal in quadruple-precision.
(Our FFT implementation required adjusting to $n=48$.)
The transformation to Liouville normal form for the FPP method
involves solving the \emph{nonlinear} system
$dx/ds=\sqrt{\Psi(x(s))}$, $dU/ds = W$, $dW/ds = [Q(s) -
\lambda]U$, where $x(s)$ is needed to compute $Q(s)$. Since fully
implicit methods are difficult to implement for nonlinear
equations, we used a spectral deferred correction scheme \cite{dutt}
in the FPP15d3 and FPP31q3 cases. The scheme orders in the table were
chosen optimally by trial and error in all but the FPP4d3 cases.

The running times above show that the choice of timestepper is
critical for efficiently computing spectral density functions with
high accuracy.  In double-precision, our method is 650 times faster
than the FPP algorithm described in \cite{fulton2008c}. However, we
were able to improve their algorithm to be 1.66 times faster than ours
by implementing a better timestepper.  In quadruple-precision, our
method is 16.6 times faster than their method using the best
timestepper available for each method.  To explain this, we recall
\cite{fulton2008c,wilkening:irk} that the error
$|f_n(s;\lambda)-\rho'(\lambda)|$ in their method decays like
$O(s^{-22/5})$ when $n=3$. Thus, to reduce the error to $O(\delta)$,
the solution must be evolved to $s_\text{max} =
O(1/\delta)^{5/22}$. Empirically, $s_\text{max} \approx
0.5(1/\delta)^{5/22}\lambda^{-7/10}$ works well in double and
quadruple-precision arithmetic over $e^{-4}\le\lambda\le e^{14}$,
though establishing precise dependence on $\lambda$ is difficult.
Indeed, for smaller values of $\lambda$ this will not be adequate to
traverse the growth region of Figure~\ref{fig:qplot}.  (The potential
$Q(s)$ must drop below $\lambda$ for the convergence theory of
\cite{fulton2008c} to be valid).  From (\ref{eq:N:steps2}) and
Lemma~\ref{lem:x}, this translates into $N_\text{steps} \approx 6
K_\nu(\delta)\lambda^{-1/5}(1/\delta)^{5/22}$. By contrast, our method
requires $N_\text{tot}\approx 0.24
K_\nu(\delta)(11+\lambda^{-11/8})\ln^3(1/\delta)$ steps over this
range of $\lambda$, where again the $\lambda$-dependence is
empirical and not precisely known.  Since
$(1/\delta)^{5/22}$ exceeds $\ln^3(1/\delta)$ once $\delta<10^{-23}$,
this is roughly the transition point where our method should become
more efficient. In practice, as shown in the table above, the methods
are already comparable in double-precision ($\delta=10^{-15}$) since
converting to Liouville normal form in the FPP method introduces the
nonlinear equation $dx/ds = \sqrt{\Psi(x)}$.  A more detailed
comparison of the convergence rates of the two methods will be given
in \cite{wilkening:irk}.

An advantage of the FPP method is that $\lambda$ is real, so the
solution remains real and complex linear algebra is avoided.  This
comes at the cost of having to evaluate higher derivatives of $Q(s)$
that appear in the formulas for $P_n$, $Z_n$ and $R_n$ in
(\ref{eq:PZR}) for $n\ge 3$; see \cite{fulton2008c,wilkening:irk} for
details.  This is a significant issue when $Q(s)$ is complicated, and
could lead to loss of accuracy if the derivatives are computed
numerically.  In our case, $Q(s)$ is the term in parentheses in
(\ref{eq:F}) below, with $x$ replaced by $x(s)$. In
\cite{wilkening:irk}, to make this practical, we had to resort to
asymptotics, making use of Lemmas~\ref{lem:psi} and~\ref{lem:x}
from Appendix~\ref{sec:lem} below to evaluate $Q'''(s)$ in $Z_3$.  By
contrast, our method involves only the original ODE. No derivatives of
the potential need to be computed, and the equation need not even be
converted to Liouville normal form.

\section{Conclusion}

We have studied the dynamics of a model partial differential equation
that is used in plasma physics to compare the merits of different
discretization schemes for the speed variable in numerical solvers. To
do so, we used the spectral transform associated with a singular
Sturm-Liouville operator $L$ to represent the solution in such a way
that the dynamics becomes trivial, through multiplication by
$e^{-\lambda t}$.  Our algorithm relies on expressing the spectral
density function of $L$ in terms of the Titchmarsh-Weyl $m$-function,
evaluating the $m$-function along a line segment in the complex plane,
and extrapolating it to the real axis using Chebyshev polynomials.
Our method is very general, and will work for any singular
Sturm-Liouville problem that is of limit point type at infinity and of
limit circle type (or regular) at the origin, and if the $m$-function
has enough smoothness near the real $\lambda$-axis to be well
approximated by polynomials along line segments in the transverse
direction. Furthermore, the complexity of the new method
for computing the $m$-function has been analyzed and shown to be
comparable to other methods in double-precision, and faster in
quadruple-precision, when optimal timestepping algorithms are used for
all methods.

The solution of the PDE computed in this way can be compared with
approximate solutions obtained with the discretization methods
traditionally used in plasma physics, and provides a basis for error
quantification. Our construction of the solution and its behavior for
certain initial conditions are also very helpful in explaining the
behavior of approximate solutions. For example, we found that for
singular initial conditions the solution of the PDE often cannot be
resolved to the desired level of accuracy until $t$ surpasses a
critical value $t^*$, because the decay rate of the spectral transform
of the solution is only algebraic at $t=0$, and slow (but at least
exponential) in the moments that follow. This has strong implications
for the projected dynamics of this equation in finite-dimensional
spaces of orthogonal polynomials, and for the choice of these
polynomials, as we present elsewhere in subsequent work.

\appendix

\section{Bound on the number of steps}
\label{sec:step:bound}
A detailed analysis of arbitrary-order $\nu$-stage Runge-Kutta
collocation methods of order $2\nu$, presented elsewhere
\cite{wilkening:irk}, reveals that if such a scheme is used to evolve
the fundamental matrix $\Phi_1$ for
\begin{equation}\label{eq:ode:balanced}
  \begin{pmatrix} y \\ z \end{pmatrix}' =
  A(x)\begin{pmatrix} y \\ z \end{pmatrix},
  \qquad
  A(x) =
  \begin{pmatrix} 0 & \sqrt{\lambda/\Psi(x)} \\
          (V(x)-\lambda)/\sqrt{\lambda\Psi(x)} &
        -\Psi'(x)/(2\Psi(x))
\end{pmatrix} 
\end{equation}
from $x$ to $x+h$, with $\Phi_1(x)=I$ and $\lambda\in\mathbb{C}$, the
local truncation error is bounded by
\begin{equation}\label{eq:trunc:error}
  \|\Phi_{1,\text{numerical}}(x+h) - \Phi_{1,\text{exact}}(x+h)\| \le
  \big(1.15\alpha^{-1}\big)^{2\nu-1}e^{2rM(x)},
\end{equation}
where $M(x)=5|\lambda|^{1/2}\la x\ra^{3/2} + 3.5 |\lambda|^{-1/2}\la
x\ra^{1/2}$ is a bound on $\max_{|z-x|\le 1/5} \|A(z)\|$,
$(z\in\mathbb{C}, x\in\mathbb{R})$, $\la x\ra=\sqrt{1+x^2}$, $\alpha$
controls the size of a Bernstein ellipse \cite{trefethen:book} with
foci at $x$ and $x+h$ and semi-major and semi-minor axes of length
$h(\alpha+\alpha^{-1})/4$ and $h(\alpha-\alpha^{-1})/4$, respectively,
and $r$ is the radius of a disk centered at $x$ in the complex plane
containing the Bernstein ellipse, which requires $h/r\le
4\alpha/(\alpha+1)^2$. To obtain (\ref{eq:trunc:error}), it is assumed
in \cite{wilkening:irk} that $\nu\ge5$, $\alpha\ge3$, $r\le1/5$,
$hM\le 2/3$.  The first-order system (\ref{eq:ode:balanced}) is
equivalent to (\ref{eq:y:ode}) and (\ref{eq:y:ode:intro}) when $z$ is
defined by $z(x)=\sqrt{\Psi(x)/\lambda}\,y'(x)$, which is scaled so
that $z$ grows at the same rate as $y$ (like $x^{3/4}$) as
$x\rightarrow\infty$ when $\lambda$ is real. In practice, all three
systems (\ref{eq:y:ode}), (\ref{eq:y:ode:intro}) and
(\ref{eq:ode:balanced}) perform similarly, but the error analysis is
simplest for (\ref{eq:ode:balanced}).

We interpret (\ref{eq:trunc:error}) as a relative error in advancing
any fundamental matrix $\Phi(x)$ from $x$ to $x+h$ since
$\Phi(x+h)=\Phi_1(x+h)\Phi(x)$. In floating point arithmetic with
roundoff threshold $\delta$, the accuracy of the result will cease to
improve when this relative error reaches $O(\delta)$.  The right-hand
side of (\ref{eq:trunc:error}) will be less than $\delta$ if we assume
$\nu\ge5$, choose $r$ so $rM=5/3$, and require $\alpha\ge
(5/3)(1/\delta)^{1/(2\nu-1)}$. Increasing $\alpha$ to 3 if necessary,
the condition $h=9r/(4\alpha)$ implies $h/r\le4\alpha/(\alpha+1)^2$,
and we are led to the stepsize constraint
\begin{equation}\label{eq:Knu}
  h(x)M(x)\le\frac{1}{K_\nu(\delta)}, \qquad
  K_\nu(\delta) =
  \max\left(\frac{3}{2},\frac{4}{9}(1/\delta)^{1/{(2\nu-1)}}\right).
\end{equation}
The intermediate assumptions that $r\le1/5$, $hM\le2/3$, and
$\alpha\ge3$ are ensured by (\ref{eq:Knu}) since
$rM=5/3$, $M(x)\ge2\sqrt{5\times 3.5}\la x\ra$,
and $\alpha=9rM/(4hM)\ge45/8$.

To derive (\ref{eq:N:steps}), we note that counting steps using the
largest stepsize allowed by (\ref{eq:Knu}) will give the points
$x_{n+1} = x_n + h(x_n)$. This is Euler's method for the ODE
$dx/ds=h(x)$ with steps of size $\Delta s=1$. The change in $s$ after
$N$ steps is then $N$.  Since $h(x)=[K_\nu M(x)]^{-1}$ is a positive,
decreasing function, solutions of the continuous problem are
increasing and concave down. Hence, Euler's method will overpredict
the solution $x(s)$ of this ODE.  As a result, solving the ODE will
overpredict the change in $s$ needed to achieve a specified change in
$x$ using Euler's method. Using separation of variables, we conclude
that $\int_{x_1}^{x_2} K_\nu M(x)\,dx$ is an upper bound on the number
of steps required to advance the solution of (\ref{eq:ode:balanced})
from $x_1$ to $x_2$ with the maximum stepsize allowed by
(\ref{eq:Knu}).

\section{Technical lemmas}
\label{sec:lem}
In this section we present four technical lemmas needed in
Appendix~\ref{sec:asym} to establish the asymptotic behavior of the
solutions $u$ and $y$ of (\ref{eq:eval}) and
(\ref{eq:y:ode:intro}), respectively.

%%%%%%%% LEMMAS %%%%%%%%%%

\vspace*{1ex}
\begin{lemma} \label{lem:psi} The function $\Psi(x)$ in (\ref{eq:psi})
  is real analytic, even, positive and satisfies
  $\Psi(x)\doteq 1/(2x^3)$, where
  $f(x)\doteq g(x)$ means that $f^{(n)}(x)-g^{(n)}(x)=o(x^{-k})$
  as $x\rightarrow\infty$ for all integers $n\ge0$, $k\ge0$.
\end{lemma}

\vspace*{1ex}
\begin{proof}
  Taylor expansion shows that $\Psi(z)=(2z^3)^{-1}\big[\erf(z) -
    (2/\sqrt\pi)ze^{-z^2}\big]$ has a removable singularity at $z=0$
  with limiting value $\Psi(0)=2/(3\sqrt\pi)\approx0.3761$; thus,
  $\Psi(z)$ is entire. It is even since $2z^3$, $\erf(z)$ and
  $ze^{-z^2}$ are odd.  The formula
\begin{equation*}
  \Psi'(x) = -4\pi^{-1/2}\int_0^x (s/x)^4 e^{-s^2}ds < 0 \qquad (x>0)
\end{equation*}
shows that $\Psi(x)$ is decreasing on $(0,\infty)$.  Since
$\lim_{x\rightarrow\infty} \Psi(x)=0$, it follows that $\Psi(x)>0$ for
$x\ge0$.  Since $\Psi(x)$ is even, it is positive for $x<0$ as well.
A straightforward induction argument shows that $g(x) = (2x^3)^{-1}
- \Psi(x)$ has derivatives of the form
\begin{equation}\label{eq:dn:psi}
  g^{(n)}(x) = \frac{(-1)^n}{\sqrt{\pi}}\left[
    \frac{(n+2)!}{2x^{n+3}}\int_x^\infty e^{-s^2}\,dx +
    \sum_{j=0}^n c_{nj} x^{2j-n-2}e^{-x^2}\right],
\end{equation}
where $c_{00}=1$, $c_{n0}=(n+1)(n!/2 + c_{n-1,0})$, $c_{nn} =
2c_{n-1,n-1}$ and
\begin{equation*}
  c_{nj} =  (n+1-2j)c_{n-1,j} + 2c_{n-1,j-1}\qquad (1\le j\le n-1).
\end{equation*}
Since $\int_x^\infty e^{-s^2}\,dx\le (2x)^{-1}e^{-x^2}$, there is
a polynomial $p_n(x)$ of degree $n+1$ such that
$|g^{(n)}(x)|x^k\le p_n(x^2)x^{k-n-4}e^{-x^2}$, which converges to
$0$ as $x\rightarrow\infty$, as claimed.
\end{proof}

%%%%%%%%%%%%%%%%%%%%%%%%

\vspace*{1ex}
\begin{lemma}\label{lem:pows}
  Let $q\in\mathbb{R}$ and define $h(z)=z^q$.  Suppose $g(x)$ and its
  derivatives grow slowly as $x\rightarrow\infty$, i.e.~there exist
  integers $k_n\ge0$ such that $g^{(n)}(x) = O(x^{k_n})$ for
  $n\ge0$. Suppose also (increasing $k_0$ if necessary) that
  $g(x)^{-1}=O(x^{k_0})$. Then $f\doteq g$ implies $h\circ f\doteq
  h\circ g$, provided one of the following is true: $q$ is an integer;
  $g$ and $f$ are real-valued; or the inverse of the distance from
  $g(x)$ to the negative real axis is $O(x^{k_0})$.
\end{lemma}

\vspace*{1ex}
\begin{proof}
  By hypothesis, there exist $c_0\ge1$, $x_0\ge1$ such that
\begin{equation}\label{eq:g:bnd}
  c_0^{-1} x^{-k_0} \le \big|g(x)\big| \le c_0 x^{k_0}, \qquad
  (x\ge x_0).
\end{equation}
Increasing $x_0$ if necessary, we may assume $\big| g(x) - f(x) \big|
\le (2c_0)^{-1}x^{-k_0}\le c_0 x^{k_0}$ for $x\ge x_0$. Any
point $\zeta$ on the line segment $\gamma(x)$ joining $g(x)$ to $f(x)$
in the complex plane satisfies $|g(x)-\zeta|\le|g(x)-f(x)|$, and hence
\begin{equation}\label{eq:line:seg}
  (2c_0)^{-1}x^{-k_0} \le |\zeta| \le (2c_0)x^{k_0}, \qquad
  (x\ge x_0, \; \zeta\in\gamma(x)).
\end{equation}
In the complex case, if $q$ is not an integer, we also have $|\im
g(x)|\ge c_0^{-1}x^{-k_0}$ whenever $x\ge x_0$ and $\re g(x)\le0$;
thus, $g(x)$ is closer to each $\zeta\in\gamma(x)$ than to the
negative real axis and $\gamma(x)$ does not cross the branch cut of
$h(z)$. If $q$ is an integer, there is no branch cut. If $f$ and $g$
are real-valued, then $f(x)$ and $g(x)$ have the same sign for $x>x_0$
(since $g(x)$ is closer to $f(x)$ than to the origin), so $\gamma(x)$
is either a subset of the positive real axis or lies along the ray
from the origin through $(-1)^q$. Either way, $\gamma(x)$ does not
cross the branch cut.  Next, for any $n\ge0$ and $x\ge x_0$, we have
\begin{equation*}
  \Big| h^{(n)}(f(x)) - h^{(n)}(g(x)) \Big| =
  \bigg| \int_{\gamma(x)} h^{(n+1)}(\zeta)\,d\zeta\bigg| \le
  \max_{\zeta\in\gamma(x)}\big| h^{(n+1)}(\zeta)\big|\,\big|f(x) - g(x)\big|.
\end{equation*}
Multiplying by $x^k$, where $k\ge0$, and using (\ref{eq:line:seg}),
we obtain
\begin{equation*}
  | \cdots |x^k \;\; \le \;\;
  q(q-1)\cdots(q-n)(2c_0 x^{k_0})^{|q-n-1|}x^k|f(x)-g(x)|
  \;\; \rightarrow \;\; 0,
  \quad (x\rightarrow\infty).
\end{equation*}
Since $h^{(n)}(g(x))$ grows slowly (i.e.~polynomially) in $x$, the
$k=0$ case also implies that $h^{(n)}(f(x))$ grows slowly.  Finally,
we use Fa\`a-di Bruno's formula \cite{faa:dibruno}
\begin{equation*}
  \frac{d^n}{dx^n} h(f(x)) = \sum_{\pi\in P_n} h^{(|\pi|)}(f(x))\prod_{b\in\pi}
  f^{(|b|)}(x),
\end{equation*}
where $P_n$ is the set of partitions of $\{1,2,\dots,n\}$, $|\pi|$ is
the number of blocks in the partition $\pi$, and $|b|$ is the number
of integers in block $b$.  Thus, for any integer $k\ge0$,
\begin{equation}\label{eq:faa:bnd}
\begin{aligned}
  &\big| (h\circ f)^{(n)}(x) - (h\circ g)^{(n)}(x) \big|x^k
  \le \\
  &\hspace*{0.5in} \sum_{\pi\in P_n}
  \bigg| h^{(|\pi|)}(f(x))\prod_{b\in\pi} f^{(|b|)}(x) - 
    h^{(|\pi|)}(g(x))\prod_{b\in\pi} g^{(|b|)}(x) \bigg|x^k.
\end{aligned}
\end{equation}
Subtracting and adding telescoping terms, e.g.
\begin{equation*}
  \big|Aa_1a_2-Bb_1b_2\big|\le\big|Aa_1(a_2-b_2)\big| +
  \big|A(a_1-b_1)b_2\big| + \big|(A-B)b_1b_2\big|,
\end{equation*}
the right-hand side of (\ref{eq:faa:bnd}) is bounded by a finite sum
of terms in which one factor is a difference, either
$\big|f^{(|b|)}(x)-g^{(|b|)}(x)\big|$ or
$\big|h^{(|\pi|)}(f(x))-h^{(|\pi|)}(g(x))\big|$, and the rest grow
slowly in $x$.  Since the difference converges to zero faster than any
polynomial, the right-hand side of (\ref{eq:faa:bnd}) converges to
zero as $x\rightarrow\infty$, as claimed.
\end{proof}

%%%%%%%%%%%%%%%%%%%%%%%%%%%%%%

\vspace*{1ex}
\begin{lemma}\label{lem:prods}
  If $h(x)$ and its derivatives grow slowly,
  i.e.~there exist integers $k_n\ge0$ such that $h^{(n)}(x) = O(x^{k_n})$
  for $n\ge0$, then $f\doteq g$ implies $hf\doteq hg$.  It suffices to
  check that $h\doteq H$ with $H^{(n)}(x)=O(x^{k_n})$.
\end{lemma}

\vspace*{1ex}
\begin{proof}
  We see that $\left|\left( [hf]^{(n)}(x) - [hg]^{(n)}(x) \right)x^k\right|$
  may be bounded by
\begin{equation}\label{eq:prod:bnd}
  \sum_{j=0}^n {n \choose j} C_j x^{k_j+k} \left| f^{(n-j)}(x)
    - g^{(n-j)}(x)\right|
\end{equation}
for large $x$, where $C_j$ is a bound on $|h^{(j)}(x)|x^{-k_j}$ for
large $x$.  But (\ref{eq:prod:bnd}) converges to zero as
$x\rightarrow\infty$ due to $f^{(n-j)}(x) = g^{(n-j)}(x) +
o(x^{k_j+k})$.  Finally, we note that if $h\doteq H$, then
$\big|h^{(n)}(x) - H^{(n)}(x)\big|$ can be made smaller than
any multiple of $x^{k_n}$ for large $x$, so if one is
$O(x^{k_n})$, so is the other.
\end{proof}

\vspace*{1ex}
\begin{lemma} \label{lem:x} The solution $x(s)$ of
  $x'=\sqrt{\Psi(x)}$, $x(0)=0$ exists for all $s\in\mathbb{R}$ and is
  an increasing, real analytic, odd function of $s$.  There is a
  constant $c$ such that $x(s)\doteq \sqrt[5]{25/8}(s-c)^{2/5}$.
\end{lemma}

\vspace*{1ex}
\begin{proof}
  By Lemma~\ref{lem:psi}, $\Psi(x)$ is real analytic, even, and
  positive.  It follows that $x(s)$ is increasing, odd, and real
  analytic for as long as the solution exists (see \cite{coddington}
  regarding analyticity).  Moreover, $s=\int_0^x \Psi(r)^{-1/2}\,dr$,
  which gives
\begin{equation}\label{eq:sx}
  s - \frac{2\sqrt{2}}{5} x^{5/2} =
  c - \int_x^\infty \big[\Psi(r)^{-1/2} - \sqrt 2 r^{3/2}\big]\,dr,
\end{equation}
where $c=\int_0^\infty [\Psi(x)^{-1/2} - \sqrt{2}x^{3/2}]\,dx \approx
1.6247$.  Since the integrand in (\ref{eq:sx}) and each of its derivatives
is $o(r^{-k})$ for all $k\ge0$, we may apply Lemma~\ref{lem:pows} to
conclude
\begin{equation}\label{eq:xs:err}
  x = \sqrt[5]{25/8}(s-c)^{2/5} + \veps_0(x), \qquad \veps_0(x)\doteq 0.
\end{equation}
Since $dx/ds=\sqrt{\Psi(x)}$, we see that
\begin{equation*}
  \frac{d^nx}{ds^n} = \frac{d^n}{ds^n}\left[\sqrt[5]{25/8}(s-c)^{2/5}\right]
  + \veps_n(x), \quad
  \veps_n(x) = \sqrt{\Psi(x)}\,\frac{d}{dx}\veps_{n-1}(x),
  \quad (n\ge1).
\end{equation*}
By Lemmas~\ref{lem:psi} and~\ref{lem:pows},
$\sqrt{\Psi(x)}\doteq(2x^3)^{-1/2}$, which has derivatives that grow
slowly (in fact decay).  Thus, we may apply Lemma~\ref{lem:prods}
inductively to conclude that $\veps_n(x)\doteq 0$ for $n\ge0$. Since
$\veps_n(x)=o(x^{-k})$ for any $k\ge0$, it follows from
(\ref{eq:xs:err}) that $\veps_n(x(s))=o(s^{-k})$ for $k\ge0$, as
claimed.
\end{proof}

%%%%%%%%%%%%%%%%%%%%%%%%%

\section{Asymptotics of the ODE}
\label{sec:asym}
In this section, we study the asymptotic behavior of solutions of
\begin{equation}\label{eq:a1}
  -(\Psi w u')' = \lambda w u
\end{equation}
for large $x$, where $\Psi(x)=[\erf(x)-x\erf'(x)]/(2x^3)$ and
$w(x)=x^2 e^{-x^2}$. The case $\lambda=0$ does not require asymptotic
arguments as the general solution (\ref{eq:u:lam0}) can be written
down in closed form.  We will show that if $\lambda\ne0$, two linearly
independent solutions of (\ref{eq:a1}) exist of the form
\begin{equation}\label{eq:asym:u}
  u_{\pm}(x) = x^{-1/4} e^{x^2/2} P_0(x) \exp\left\{\pm i
    \sqrt{\frac{8\lambda}{25}}P_1(x)x^{5/2}\right\}
  \big[1 + O\big(x^{-7/2}\big)\big]
\end{equation}
for $x\gg1$, where
\begin{equation}\label{eq:pq:def}
\begin{aligned}
  P_0(x) &= 1 + \frac{1}{8x\lambda} + \frac{5}{128x^2\lambda^2} +
  \frac{15}{1024x^3\lambda^3}, \\
  P_1(x) &= 
    1 - \frac{5}{12x\lambda} - \frac{5}{32x^2\lambda^2} +
    \frac{5}{128x^3\lambda^3} + \frac{25}{6144x^4\lambda^4} +
    \frac{7-1152\lambda^4}{8192x^5\lambda^5}.
\end{aligned}
\end{equation}
When $\lambda$ is real and positive,
any real-valued solution of (\ref{eq:a1}) may still be written as a
linear combination $u(x)=Au_+(x)+Bu_-(x)$. Reality requires
$B=\overline{A}$, which implies
\begin{equation}
  u(x) = C x^{-1/4} e^{x^2/2} P_0(x) \cos\left\{
    \sqrt{8\lambda/25}\,P_1(x)x^{5/2} - \theta \right\}
  \big[1 + O\big(x^{-7/2}\big)\big]
\end{equation}
for some $C,\theta\in\mathbb{R}$, which yields (\ref{eq:asym}).  For
all other values of $\lambda\in\mathbb{C}\setminus\{0\}$, one of the
modes $u_\pm(x)$ in (\ref{eq:asym:u}) grows super-exponentially as
$x\rightarrow\infty$ while the other decays.  The recessive (decaying)
mode is uniquely determined by (\ref{eq:asym:u}), but the dominant
(growing) mode is not. When $\lambda>0$, neither mode dominates the
other, so both are determined uniquely by (\ref{eq:asym:u}).

While it is possible to derive (\ref{eq:asym:u}) by guessing its form
and computing successive terms of $P_0(x)$ and $P_1(x)$ from
(\ref{eq:a1}) iteratively, it is difficult to prove error bounds for
the resulting series.  Instead, we will use a Liouville transformation
\cite{birkhoff} to convert the ODE to normal form
and use WKB theory \cite{bender,olver} to study the asymptotics.

To convert the general second order self-adjoint equation $-(pu')' +
qu = \lambda w u$ to Liouville normal form \cite{birkhoff}, $-U''(s) +
Q(s)U(s) = \lambda U(s)$, one solves the ODE $dx/ds=\sqrt{p(x)/w(x)}$
to obtain $x(s)$, and then defines
\begin{equation*}
  U(s) = \gamma(s)u(x(s)), \qquad
  \gamma(s)=\sqrt[4]{p(x(s))w(x(s))}, \qquad
  Q(s)=\frac{q(x(s))}{w(x(s))} + \frac{\gamma''(s)}{\gamma(s)}.
\end{equation*}
To fit in the framework of WKB theory, it is convenient to change the
sign of $Q$ and absorb $\lambda$ into the potential.  Thus, we convert
(\ref{eq:a1}) to the form
\begin{equation}\label{eq:can_form}
  -U''=Q(s)U
\end{equation} 
by the change of variables $dx/ds = \sqrt{\Psi(x)}$,
$U(s)=\Psi(x(s))^{1/4}y(x(s))$, $y(x)=w(x)^{1/2}u(x)$,
$Q(s)=F(x(s))$, and
\begin{equation}\label{eq:F}
  F(x) = \lambda - \left(\frac{1}{4}\Psi''(x) - \frac{1}{16}\Psi(x)^{-1}\Psi'(x)^2 +
  (1-x^2)\frac{\Psi'(x)}{x} + (x^2 - 3)\Psi(x)\right).
\end{equation}
The same result is obtained if we start from (\ref{eq:y:ode:intro})
instead of (\ref{eq:a1}).  Since $\Psi$ is even and entire,
$\Psi'(z)/z$ has a removable singularity at $z=0$.  By
Lemma~\ref{lem:psi}, $0<\Psi(x)\le2/(3\sqrt\pi)$ on the real
axis. Thus, $F(x)-\lambda$ and $\sqrt{\Psi(x)}$ are real analytic on
all of $\mathbb{R}$. By Lemma~\ref{lem:x}, $Q(s)-\lambda$ is also real
analytic. (We subtract $\lambda$ to make $F(x)$ and $Q(s)$ real valued
for real arguments).  By Lemmas~\ref{lem:psi}, \ref{lem:pows}
and~\ref{lem:prods}, we may replace $\Psi(x)$ by $1/(2x^3)$ in
(\ref{eq:F}) to conclude
\begin{equation}\label{eq:Ftil}
  F(x)\doteq
  \lambda - \frac{1}{2}x^{-1} + \frac{9}{32}x^{-5},
\end{equation}
where $f(x)\doteq g(x)$ means that $f^{(n)}(x)-g^{(n)}(x)=o(x^{-k})$
as $s\rightarrow\infty$ for all integers
$n\geq0,\;k\geq0$. Lemmas~\ref{lem:pows} and~\ref{lem:prods} were both
used (the latter twice) to convert $\Psi(x)^{-1}\Psi'(x)^2$ into
$(9/2)x^{-5}$.

The WKB approximation \cite{bender,olver} of the solution of
(\ref{eq:can_form}) is
\begin{equation}\label{eq:WKB}
  U_\pm(s) \sim
  Q(s)^{-1/4}\exp\left\{\pm i\int^s\sqrt{Q(r)}\,dr\right\},
\end{equation}
where $f(s)\sim g(s)$ means $f(s)/g(s)\rightarrow1$ as
$s\rightarrow\infty$.  When $Q(s)$ is real and positive,
(\ref{eq:WKB}) can be derived by performing another Liouville
transformation, namely $d\xi/ds = \sqrt{Q(s)}$,
$W(\xi(s))=Q(s)^{1/4}U(s)$, to convert the ODE to $d^2W/d\xi^2 =
[-1+\phi(\xi)]W$. Neglecting $\phi$ gives $W\sim e^{\pm i\xi}$. If
$\phi(\xi)$ is small, error estimates can be derived in the $\xi$
coordinate system \cite{olver}.

In our case, $Q(s)$ is complex-valued, so the change of variables
$\xi=\int^s \sqrt{Q(r)}\,dr$ requires that we complexify the dependent
variable \cite{olver}.  We prefer to work with complex functions of
the real variable $x$, and only use $s$ and $\xi$ as intermediate
steps to finding a representation of (\ref{eq:a1}) that is suitable
for perturbation analysis.  To this end, we still define $W=Q^{1/4}U$,
but treat it as a function of~$s$ rather than~$\xi$.  Substitution of
the identity
\begin{equation}\label{eq:dt2}
  \frac{d^2U}{ds^2} = 
  Q^{1/4} \frac{d}{ds}\left[Q^{-1/2}\frac{d}{ds}
    \left(Q^{1/4}U\right)\right] +
  Q^{1/4} \left(\frac{d^2}{ds^2}Q^{-1/4}\right)U
\end{equation}
in the equation $d^2U/ds^2 = -QU$ gives
\begin{equation}\label{eq:W2}
  Q^{-1/2} \frac{d}{ds}Q^{-1/2}\frac{dW}{ds}
  = [-1+\phi(s)]W, \qquad
  \phi(s) = -Q^{-3/4}\frac{d^2}{ds^2}Q^{-1/4},
\end{equation}
where we have adopted the convention that differential operators act
on all products that follow them unless otherwise indicated by
parentheses.  The left-hand side of (\ref{eq:W2}) plays the role of
$d^2W/d\xi^2$ in \cite{olver}.  Changing back to the $x$-coordinate
system and writing $W(s)=v(x(s))$,
we find that if $v(x)$ and $u(x)$ are related by
\begin{equation}\label{eq:uv:relation}
  u(x) = w(x)^{-1/2}[F(x)\Psi(x)]^{-1/4}v(x),
\end{equation}
then $u(x)$ satisfies (\ref{eq:a1}) iff $v(x)$ satisfies
\begin{equation}\label{eq:v:psi}
  \sqrt{\frac{\Psi}{F}}\frac{d}{dx}\sqrt{\frac{\Psi}{F}}
    \frac{dv}{dx}
  = \left(-1+\sqrt{\frac{\Psi}{F}}\,\psi\right)v, \qquad
  \psi = -F^{-1/4}\frac{d}{dx}\sqrt{\Psi}\frac{d}{dx}
    F^{-1/4}.
\end{equation}
We again recognize the left-hand side as $d^2W/d\xi^2$.

Two technical issues concern zeros of $F(x)$ and the branch cut of the
square root and fourth root functions along the negative real axis.
We claim there is an $R_1\ge0$ such that $F(x)$ is bounded away from
zero and neither $F(x)$ nor $\Psi(x)/F(x)$ crosses the negative real
axis for $x\ge R_1$.  Since we are interested in the asymptotics of
$u(x)$ for large $x$, we only need to solve (\ref{eq:v:psi}) for $x\ge
R_1$.  If $\lambda\not\in\mathbb{R}$, then $|F(x)|\ge|\im\lambda|$ since
$\im\{F(x)\}=\im\lambda$, and $R_1=0$ suffices.  If $\lambda$ is real
(and non-zero), then $F(x)$ may have zeros, but by (\ref{eq:Ftil})
there is an $R_1\ge 0$ such that for $x\ge R_1$, $F(x)$ and $\lambda$ have
the same sign, and $|F(x)|\ge |\lambda/2|$.  If the sign is positive,
the branch cut is avoided, and if the sign is negative, we treat
$\im\lambda=\im\{F(x)\}=0^+$ in all formulas involving fractional
powers of $\lambda$ or $F(x)$, e.g.~$\sqrt{\Psi(x)/F(x)}=
-i\sqrt{\Psi(x)/|F(x)|}$.

If $\psi$ were zero in (\ref{eq:v:psi}), $v(x)=e^{\pm i\xi(x)}$ would
be independent solutions, where
\begin{equation}\label{eq:xi:def}
  \xi(x) = \int^x \sqrt{\frac{F(r)}{\Psi(r)}}\,dr.
\end{equation}
Lemmas~\ref{lem:pows} and~\ref{lem:prods} justify replacing
$\Psi^{-1}$ by $2r^3$ and $F$ by (\ref{eq:Ftil}) to obtain an
asymptotic formula for $\sqrt{F/\Psi}$:
\begin{equation}\label{eq:psi_F_asym}
  \sqrt{\frac{F(r)}{\Psi(r)}}\doteq
  \sqrt{2\lambda} \left(1 - \frac{1}{2r\lambda} + \frac{9}{32r^5\lambda}
    \right)^{1/2}r^{3/2}.
\end{equation}
We choose the integration constant in (\ref{eq:xi:def}) so that
\begin{equation}\label{eq:xi:asym}
  \xi(x)=(8\lambda/25)^{1/2}P_1(x)x^{5/2}+O(x^{-7/2}), \qquad
  (x\rightarrow\infty),
\end{equation}
where $P_1(x)$ was defined in (\ref{eq:pq:def}). This result is
obtained by expanding the asymptotic formula for $\sqrt{F/\Psi}$ in a
binomial series and integrating term by term.  We now look for
solutions of (\ref{eq:v:psi}) that are perturbations of $e^{\pm
  i\xi(x)}$:
\begin{equation}\label{eq:vpm}
  v_+(x)=e^{i\xi(x)}(1+h_+(x)), \qquad
  v_-(x)=e^{-i\xi(x)}(1+h_-(x)).
\end{equation}
The functions $h_+(x)$ and $h_-(x)$ must satisfy
\begin{equation}\label{eq:h:ODE}
  \frac{d}{dx}\sqrt{\frac{\Psi}{F}}\frac{dh}{dx} \pm 2i\frac{dh}{dx} =
  \psi(x)(1+h(x)).
\end{equation}
Solutions of the homogeneous problem (with $\psi\equiv0$) are
$h_\pm\equiv1$ and $h_\pm(x) = e^{\mp 2i\xi(x)}$.  Using variation of
parameters \cite{coddington} to solve the non-homogeneous problem
yields an integral equation for the solution with initial conditions
$h(x_0)=h_0$, $h'(x_0)=0$:
\begin{equation}\label{h:int:eq}
  h_\pm(x) = h_0 \pm \frac{1}{2i}\int_{x_0}^x \left(
    1 - e^{\pm 2i\{\xi(y) - \xi(x)\}} \right)
  \psi(y) [1 + h_\pm(y)]\,dy.
\end{equation}
The location of $x_0$ will be chosen below, and depends on $\lambda$
and the $\pm$ sign.  Simplifying the formula for $\psi$ in
(\ref{eq:v:psi}) gives its behavior to leading order as
$x\rightarrow\infty$:
\begin{equation}\label{eq:psi:asym}
  \psi(x)
  = \frac{2\Psi'FF' - 5\Psi(F')^2+4\Psi FF''}{16\Psi^{1/2}F^{5/2}}
  \sim \frac{-7\sqrt{2}}{32\lambda^{3/2}x^{9/2}}, \qquad
  (x\rightarrow\infty).
\end{equation}
Our goal is to use this in (\ref{h:int:eq}) to produce solutions
$h_\pm(x)$ of (\ref{eq:h:ODE}) that decay like $O(x^{-7/2})$.
If we are successful, then (\ref{eq:asym:u}) will follow from
(\ref{eq:vpm}), (\ref{eq:uv:relation}), (\ref{eq:xi:asym}) and
\begin{equation}\label{eq:p:derive}
  w(x)^{-1/2}[F(x)\Psi(x)]^{-1/4} = 
  \frac{
  x^{-1/4}e^{x^2/2}P_0(x)\big[1 + O(x^{-4})\big]}{(\lambda/2)^{1/4}},
  \qquad (x\rightarrow\infty),
\end{equation}
with $P_0(x)$ as in (\ref{eq:pq:def}).  Indeed, the factor of
$(\lambda/2)^{-1/4}$ can be dropped by linearity, and the factors of
$1+O(x^{-4})$, $\exp\{O(x^{-7/2})\}$ and $1+h_\pm(x)$ combine to make
$1+O(x^{-7/2})$ in (\ref{eq:asym:u}), as claimed.  The coefficients
of $P_0(x)$ were obtained from (\ref{eq:p:derive}) using
(\ref{eq:Ftil}), $\Psi(x)\doteq(2x^3)^{-1}$, Lemmas~\ref{lem:pows}
and~\ref{lem:prods}, and the binomial series.

Let us therefore study solutions of the integral equation
(\ref{h:int:eq}).  Let $\alpha=\im\sqrt{\lambda}$ and $\sgn=\pm1$,
depending on the case considered in (\ref{eq:vpm}).  We drop the $\pm$
subscript on $h$ and define
\begin{equation*}
  K(x,y) = \frac{\beta}{2i}\big[ 1 - E(x,y) \big]\psi(y), \qquad
  E(x,y) = e^{2i\beta\{\xi(y)-\xi(x)\}}, \qquad (\beta=\pm1)
\end{equation*}
which appear in (\ref{h:int:eq}).
We claim that there is an $R\ge R_1$ (depending on $\lambda$ and
$\sgn$) such that one of the following holds:
\begin{align*}
    &\text{case 1: } & y\ge x\ge R \;\;
    &\Rightarrow \;\; |K(x,y)|\le|\psi(y)|, \\
    &\text{case 2: } &  |\alpha|>0 \;\;\text{and}\;\; \big(
    x\ge y\ge R \;\; &\Rightarrow \;\;
    |K(x,y)|\le|\psi(y)|\,, \, |E(x,y)|\le\eta(x,y)\big),
\end{align*}
  where $\eta(x,y)=\exp\{-2|\alpha|y^{3/2}(x-y)\}$. If $\alpha=0$,
  then $\lambda>0$ and setting $R=R_{1}$ suffices to establish case
  1. Indeed, since $F(x)>0$ for $x\ge R_1$, it follows that $\xi(x)$
  in (\ref{eq:xi:def}) is real, $E(x,y)$ is on the unit circle, and
  $|1-E|\le 2$ for $x\ge R$, $y\ge R$.  If $\alpha>0$, we see from
  (\ref{eq:psi_F_asym}) that
  $\im\sqrt{F(r)/\Psi(r)}\,r^{-3/2}\rightarrow (\sqrt2\alpha)$ as
  $r\rightarrow\infty$.  Since $\sqrt{2}>1$, there exists $R\ge R_1$
  such that
\begin{equation*}
  \im\sqrt{F(r)/\Psi(r)}\ge\alpha r^{3/2}, \qquad (r\ge R).
\end{equation*}
It then follows from (\ref{eq:xi:def}) that
\begin{equation}\label{eq:xi:order}
  \im\{\xi(y)-\xi(x)\} \ge \int_x^y \alpha r^{3/2}\,dr \ge
  \alpha x^{3/2}(y-x)\ge 0, \qquad (y\ge x\ge R).
\end{equation}
Since $E(x,y)=e^{2i\sgn\{\xi(y)-\xi(x)\}}$, we see that $E(x,y)$ lies
inside or on the unit circle when $\sgn=1$ and $y\ge x\ge R$. Thus,
case 1 holds when $\sgn=1$.  Interchanging $x$ and $y$ in (\ref{eq:xi:order})
and evaluating $|E(x,y)|$ shows that case 2 holds when $\sgn=-1$. Similar
arguments show that if $\alpha<0$, cases 1 and 2 hold when $\sgn=-1$ and
$\sgn=1$, respectively.

In case 1, we set $h_0=0$ in (\ref{h:int:eq}) and send $x_0$ to
infinity. This yields the equation
\begin{equation*}
  h(x) = \mathbb{K}[1+h](x), \qquad
  \mathbb{K}[f](x) = -\int_x^\infty K(x,y)f(y)\,dy.
\end{equation*}
Increasing $R$ if necessary, we may assume
$\int_{R}^\infty|\psi(y)|\,dy\le1/2$.  Then $\|\mathbb{K}\|\le
1/2$, where $\mathbb{K}$ is regarded as an operator on
$BC\big([R,\infty)\big)$, the Banach space of bounded, continuous
functions in the uniform norm. Thus, $h = (\mathbb{K} + \mathbb{K}^2
+\mathbb{K}^3 + \cdots)1$ is the unique bounded, continuous function
that satisfies the integral equation. By the dominated convergence
theorem and Leibniz integral rule, such a solution of the integral
equation also satisfies the ODE (\ref{eq:h:ODE}). We note that
$h_n=(\mathbb{K} + \mathbb{K}^2 + \cdots + \mathbb{K}^n)1$ can be
computed via the Picard iteration $h_0(x)=0$,
$h_{n+1}=\mathbb{K}[1+h_n]$.  Standard estimates
\cite{coddington,olver} on the size of $|h_{n+1}(x)-h_n(x)|$ in the
Picard iteration scheme give the bound
\begin{equation}\label{eq:h:bnd1}
  |h(x)| \le \exp\left(\int_x^\infty |\psi(r)|\,dr\right) - 1, \qquad
  (x>R).
\end{equation}
By (\ref{eq:psi:asym}), $h(x)=O(x^{-7/2})$ as $x\rightarrow\infty$,
as required.

In case 2, we define $\mathbb{K}[f](x) = \int_{R}^x
K(x,y)f(y)\,dy$ and proceed in the same manner,
again assuming $\int_{R}^\infty |\psi(y)|\,dy\le1/2$.
This establishes existence and uniqueness of
a bounded, continuous solution of $h=\mathbb{K}[1+h]$,
along with the bound
\begin{equation}\label{eq:h:bound}
  |h(x)| \le \left(\exp\left\{\int_{R}^x |\psi(r)|\,dr\right\} - 1\right)
  \le \left(e^{1/2}-1\right) \le 2/3, \qquad (x\ge R).
\end{equation}
Although $h(x)$ will not in general approach 0 as
$x\rightarrow\infty$, we will show below that it approaches a limiting
value, $c$, with $|c|\le2/3$.  Defining
\begin{equation}
  \tilde h(x) = \frac{h(x)-c}{1+c}, \qquad h_0 = \frac{-c}{1+c},
\end{equation}
we find that $\tilde h$ satisfies $\tilde h = h_0 +
\mathbb{K}[1+\tilde h]$, and hence (\ref{eq:h:ODE}). Thus, adjusting
the initial condition from 0 to $h_0$ merely shifts and rescales the
solution. If we can show that $h(x)-c=O(x^{-7/2})$ as $x\rightarrow\infty$, then
$\tilde h(x)$ will be the desired solution of (\ref{eq:h:ODE}) that
decays as $O(x^{-7/2})$.  To prove that $c=\lim_{x\rightarrow\infty}
h(x)$ exists, we integrate (\ref{eq:h:ODE}) from $x_1$ to $x_2$,
assuming $x_2>x_1\ge R$:
\begin{equation}\label{eq:cauchy}
  h(x_2) - h(x_1) = \frac{\sgn}{2i}\int_{x_1}^{x_2} \psi(x)\big(1+h(x)\big)\,dx
  - \frac{\sgn}{2i}\sqrt{\frac{\Psi(x)}{F(x)}}\,h'(x)\bigg\vert_{x_1}^{x_2}.
\end{equation}
Recall that $\sgn=\pm1$ distinguishes the case in (\ref{eq:vpm}).  The
first term on the right is $O(x_1^{-7/2})$ due to (\ref{eq:h:bound})
and (\ref{eq:psi:asym}). If we can show that
$\sqrt{\Psi/F}h'(x)=O(x^{-7/2})$, then we are done: (\ref{eq:cauchy})
implies that for any sequence $x_1<x_2<\cdots$ with
$x_n\rightarrow\infty$, $h(x_n)$ is a Cauchy sequence.  So
$c=\lim_{x\rightarrow\infty}h(x)$ exists.  Sending $x_2$ to $\infty$
and replacing $x_1$ by $x$ in (\ref{eq:cauchy}) then gives
$c-h(x)=O(x^{-7/2})$, as required.

To show that $\sqrt{\Psi/F}h'(x)=O(x^{-7/2})$, we differentiate the
integral equation $h=\mathbb{K}[1+h]$.  This gives
\begin{equation}\label{eq:h:deriv}
  \sqrt{\frac{\Psi(x)}{F(x)}}\,h'(x) = \int_{R}^x E(x,y)\psi(y)[1+h(y)]\,dy,
  \qquad (x>R).
\end{equation}
Using the bound $|E(x,y)|\le\eta(x,y)=\exp\{-2|\alpha|y^{3/2}(x-y)\}$
for $x\ge y\ge R$, and breaking the integral in (\ref{eq:h:deriv})
into two segments of length $(x-R)/2$, we obtain
\begin{equation}\label{eq:h:deriv:bnd}
  \left| \sqrt{\frac{\Psi(x)}{F(x)}}\,h'(x) \right| \;\le\;
  \frac{5}{3}\eta\left(x,\jt \frac{x+R}{2}\right)
  \int_{R}^{\frac{x+R}{2}} |\psi(y)|\,dy +
  \frac{5}{3}\int_{\frac{x+R}{2}}^x |\psi(y)|\,dy,
\end{equation}
where $5/3$ is a bound on $|1+h(y)|$.  The first integral on the
right is bounded by $1/2$ while
\begin{equation*}
  \eta\left(x, \frac{x+R}{2}\right) =
  \exp\left\{ -|\alpha|\left( \frac{x+R}{2} \right)^{1/2}\left(
      \frac{x^2-R^2}{2} \right) \right\},
\end{equation*}
which decays super-exponentially as $x\rightarrow\infty$.  Since
$\int_x^\infty |\psi(y)|\,dy = O(x^{-7/2})$, the second
integral in (\ref{eq:h:deriv:bnd}) is bounded by
\begin{equation*}
  C \left( \frac{x+R}{2} \right) ^{-7/2} = \;
  2^{7/2} C \left(1 + \frac{R}{x}\right)^{-7/2}
  x^{-7/2} \;=\; O(x^{-7/2}).
\end{equation*}
Thus, $\sqrt{\Psi/F}\,h'(x) = O(x^{-7/2})$ as claimed.

\section{Analyticity of the spectral density function}
\label{sec:analytic}
A number of authors have proved that the spectral function $\rho(x)$
is absolutely continuous when the potential in the Schr\"odinger
equation is of bounded variation or decreases sufficiently rapidly at
infinity; see e.g.~\cite{pearson:singular, weidmann, lavine}.
However, we are not aware of any work that establishes conditions to
ensure that $\rho'(\lambda)$ will be real analytic for $\lambda>0$.
In this appendix, we prove this for the operator $L$ in
(\ref{eq:PDE}), studied throughout this paper.

Our proof will be to show that $u_\pm(x;\lambda)$ in (\ref{eq:asym:u})
can be made to depend analytically on $\lambda$ in a complex
neighborhood of each $\lambda_0>0$.  We can then construct a
fundamental matrix $\tilde\Phi(x;\lambda)$ for
(\ref{eq:1st:order:alt}) of the form
\begin{equation*}
  \tilde\Phi(x;\lambda) =
  \begin{pmatrix}
    y_+(x;\lambda) & y_-(x;\lambda) \\
    \Psi(x)y'_+(x;\lambda) & \Psi(x)y'_-(x;\lambda)
  \end{pmatrix},
\end{equation*}
where $y_\pm(x;\lambda) = xe^{-x^2/2}u_\pm(x;\lambda)$. The
construction of $u_\pm(x;\lambda)$ will involve a fixed-point
(i.e.~Picard) iteration, as in Appendix~\ref{sec:asym}, producing
solutions for sufficiently large $x$, say $x\ge x_0$, where $x_0$ does
not depend on $\lambda$ in the neighborhood.  Similar to what we did
in (\ref{eq:xstar}), we can express the fundamental matrix
$\Phi(x;\lambda)$ with correct boundary conditions at $x=0$ in the
form
\begin{equation}\label{eq:x0:def}
  \Phi(x;\lambda) = \tilde\Phi(x;\lambda)C(\lambda), \qquad\quad (x\ge x_0),
\end{equation}
where $C(\lambda) = \tilde\Phi(x_0;\lambda)^{-1}\Phi(x_0;\lambda)$
depends analytically on $\lambda$ near $\lambda_0$, since
$\Phi(x_0;\lambda)$ is an entire function of $\lambda$ (see
\cite{coddington}) while $\tilde\Phi(x_0;\lambda)^{-1}$ is analytic in
the neighborhood where $u_\pm(x;\lambda)$ are analytic.  The
$m$-function may then be written
\begin{equation*}
  m(\lambda) = -\lim_{x\rightarrow\infty}
  \frac{y_+(x;\lambda)c_{11}(\lambda) + y_-(x;\lambda)c_{21}(\lambda)}{
    y_+(x;\lambda)c_{12}(\lambda) + y_-(x;\lambda)c_{22}(\lambda)} =
  \begin{cases}
    -c_{21}(\lambda)/c_{22}(\lambda), & \im\lambda>0, \\
    -c_{11}(\lambda)/c_{12}(\lambda), & \im\lambda<0.
  \end{cases}
\end{equation*}
This equation clarifies the source of the discontinuity across the
continuous spectrum in the $m$-function: when $\im\lambda$ changes
sign, the growing and decaying solutions switch, i.e.~$y_+(x;\lambda)$
and $y_-(x;\lambda)$ reverse roles. From (\ref{eq:m:diff}), we know
$\overline{m(\lambda)}=m(\bar\lambda)$. Thus,
\begin{equation}\label{eq:rho:prime}
  \rho'(\lambda) = \frac{1}{\pi}\im\{m(\lambda^+)\} =
  \frac{1}{\pi}\frac{m(\lambda^+) - m(\lambda^-)}{2i} =
  \frac{1}{2\pi i}\,\frac{\det C(\lambda)}{c_{12}(\lambda)c_{22}(\lambda)},
\end{equation}
which is real-valued for real $\lambda$ and analytic in the
neighborhood on which $u_\pm(x;\lambda)$ can be constructed to depend
analytically on $\lambda$. Reality of the final formula in
(\ref{eq:rho:prime}) for $\lambda\in\mathbb{R}$ can be confirmed by
noting that $C(\lambda)=\tilde\Phi^{-1}\Phi$ with $\Phi$ real and
$\tilde\Phi$ having complex conjugate columns (since
$\overline{u_+(x;\lambda)} = u_-(x;\bar\lambda)$). We note that
(\ref{eq:rho:prime}) could potentially be used to compute
$\rho'(\lambda)$ without complexifying $\lambda$, but by solving a
terminal value problem as well as an initial value problem. This is
the key idea of the Fulton, Pearson, Pruess algorithm
\cite{fulton2008c}, though their derivation is different than ours and
leads to a representation, namely (\ref{eq:fpp:fn}), that does not
allow $\lambda$ to be complex.

The WKB analysis in Appendix~\ref{sec:asym} fails to produce analytic
functions $u_\pm(x;\lambda)$ because different formulas are used for
the cases when $\im\lambda$ is positive or negative. The problem
occurs in case 2, where we adjusted $h_0$ in (\ref{h:int:eq}) to
achieve $\lim_{x\rightarrow\infty} h_\pm(x)=0$. In fact,
(\ref{h:int:eq}) is not the most general solution of (\ref{eq:h:ODE});
one could add to it a term of the form $\pm\frac{h_1}{2i}\left( 1 -
  e^{\pm2i\{\xi(x_0)-\xi(x)\}}\right) \sqrt{\Psi(x_0)/F(x_0)}$, so
that $h_\pm(x_0) = h_0$ and $h'_\pm(x_0) = h_1$. Varying $h_1$ and
solving for $h_0$ (as we did with $h_1=0$) leads to a one-parameter
family of solutions of (\ref{eq:h:ODE}) that approach zero as
$x\rightarrow\infty$, consistent with the observation at the beginning
of Appendix~\ref{sec:asym} that the recessive solution is unique while
the dominant solution is not.  Our task now is to analytically
continue $u_+(x;\lambda)$ from the upper $\lambda$-half-plane into the
lower half-plane, and $u_-(x;\lambda)$ from the lower half-plane into
the upper half-plane, in order to select the ``right'' dominant
solutions.

It will be necessary to complexify $x$ in Appendix~\ref{sec:asym}, as
well as in Lemmas~\ref{lem:psi}--\ref{lem:prods}. We will use the
letter $z$ (to replace $x$) as there is no chance of confusion with
the second component of $\vec r$ in (\ref{eq:y:ode}). First, we modify
the notation $f(z)\doteq g(z)$ to mean that there is a number $R\ge0$
such that $f(z)$ and $g(z)$ are both defined and analytic on the
region $S_R=\{re^{i\theta}\,:\, r>R\,,\,|\theta|< \pi/6\}$, and, for
all non-negative integers $n$ and~$k$, $f^{(n)}(z) - g^{(n)}(z) =
o(|z|^{-k})$ as $(z\in S_R)\rightarrow\infty$, i.e.~as
$|z|\rightarrow\infty$ with $z$ remaining in $S_R$. We then have

\vspace*{1ex}
\begin{lemma}\label{lem:psi:z}
  $\Psi(z)$ satisfies $\Psi(z)\doteq 1/(2z^3)$.
\end{lemma}

\vspace*{1ex}
\begin{proof}
  Since $\Psi(z)$ is entire, both $\Psi(z)$ and $(2z)^{-3}$ are
  defined on $S_R$ with $R=0$.  Defining $g(z) = (2z)^{-1} - \Psi(z)$,
  we see that (\ref{eq:dn:psi}) holds with $x$ replaced by $z=x+iy$
  and the integral interpreted as $\int_x^\infty e^{-(s+iy)^2}ds$. 
  Note that $\big|\int_x^\infty e^{-(s+iy)^2}ds\big| \le
  e^{y^2}\int_x^\infty e^{-s^2}ds \le (2x)^{-1}e^{y^2-x^2}$. Since
  $z\in S$, we have $y^2\le x^2/3$, $|z|^2\le 4x^2/3$, and $x^2-y^2\ge
  \frac{2}{3}x^2 \ge \frac{1}{2}|z|^2$. It follows that
  $(2x)^{-1}e^{y^2-x^2}\le (\sqrt{3}|z|)^{-1}e^{-|z|^2/2}$.
  Similarly, $|e^{-z^2}|\le e^{-|z|^2/2}$. The rest of the proof of
  Lemma~\ref{lem:psi} works the same, i.e.~there is a polynomial
  $p_n(x)$ of degree $n+1$ such that $|g^{(n)}(z)||z|^k \le
  p_n(|z|^2)|z|^{k-n-4}e^{-|z|^2/2}$, which converges to zero
  as $(z\in S_R)\rightarrow\infty$.
\end{proof}

\vspace*{1ex} Since $g(z)=1/(2z^3)$ satisfies $\frac{1}{2}|z|^{-3}\le
|g(z)|\le 2|z|^3$, which is of the form (\ref{eq:g:bnd}), the proof of
Lemma~\ref{lem:pows} is easily modified to show $\Psi(z)^{-1}\doteq
2z^3$. In particular, there is an $R$ large enough that $\Psi(z)^{-1}$
has no poles for $z\in S_R$, which also follows from
$|\Psi(z)|\ge|2z^3|^{-1} - |\Psi(z)-(2z^3)^{-1}|$ and
Lemma~\ref{lem:psi:z}. Adapting Lemma~\ref{lem:prods} to the
case of analytic functions on $S_R$ is also straightforward, and
implies that
$F(z;\lambda)=\lambda-\big[\jt
\frac{1}{4}\Psi''(z) - \frac{1}{16}\Psi(z)^{-1}\Psi'(z)^2 +
  (1-z^2)z^{-1}\Psi'(z) + (z^2 - 3)\Psi(z)
\big]$
satisfies
\begin{equation}\label{eq:F:z:lam}
  F(z;\lambda)-\lambda \doteq -\frac{1}{2}z^{-1} + \frac{9}{32}z^{-5},
\end{equation}
where the left-hand side is independent of $\lambda$. Next we fix
$\lambda_0>0$ and consider $\lambda\in B$, where $B=B_a(\lambda_0)$ is
the open ball of radius $a=\lambda_0/3$ centered at $\lambda_0$. By
(\ref{eq:F:z:lam}), we may increase $R$ if necessary so that
$|F(z;\lambda)-\lambda|<\lambda_0/6$.  We then have
$|F(z;\lambda)-\lambda_0|<\lambda_0/6+\lambda_0/3=\lambda_0/2$, which
implies that $|\opn{Arg}F(z;\lambda)|<\pi/6$ for $\lambda\in B$ and
$z\in S_R$. Since $\Psi(z)^{-1}\doteq 2z^3$, we may increase $R$
further if necessary to conclude that $\Psi(z)^{-1}=2z^3(1+\veps(z))$
for $z\in S_R$, where $|\veps(z)|<\sin(\pi/12)$. As a result,
$|\opn{Arg}\Psi(z)|=\big|\opn{Arg}\big[\Psi(z)^{-1}\big]\big|<7\pi/12$
for $z\in S_R$.  Thus, $F(z;\lambda)$, $\Psi(z)^{-1}$ and
$F(z;\lambda)/\Psi(z)$ are bounded away from the branch cut of the
square root and fourth root functions for $\lambda\in B$ and $z\in
S_R$.  If we define
\begin{equation*}
  v(z;\lambda) = w(z)^{1/2}[F(z;\lambda)/\lambda]^{1/4}[2\Psi(z)]^{1/4}u(z;\lambda),
\end{equation*}
then $u(z;\lambda)$ will satisfy (\ref{eq:a1}) iff $v(z;\lambda)$ satisfies
\begin{equation}\label{eq:v:psi2}
  \sqrt{\frac{\Psi}{F}}\frac{d}{dz}\sqrt{\frac{\Psi}{F}}
    \frac{dv}{dz}
  = \left(-1+\sqrt{\frac{\Psi}{F}}\,\psi\right)v, \qquad
  \psi = -F^{-1/4}\frac{d}{dz}\sqrt{\Psi}\frac{d}{dz}
    F^{-1/4}.
\end{equation}
This may be derived analogously to (\ref{eq:uv:relation}) and
(\ref{eq:v:psi}).  We also define
\begin{equation}\label{eq:xi:def2}
  \xi(z;\lambda) = \int^z\sqrt{\frac{F(\zeta;\lambda)}
    {\Psi(\zeta)}}\,d\zeta, \qquad\quad (\lambda\in B, \, z\in S_R)
\end{equation}
and look for solutions of (\ref{eq:v:psi2}) of the
form $v_\pm(z;\lambda) = e^{\pm i\xi(z;\lambda)}[1 + h_\pm(z;\lambda)]$.
The equation for $h_\pm(z;\lambda)$ is
\begin{equation}\label{eq:h:ODE2}
  \frac{d}{dz}\sqrt{\frac{\Psi}{F}}\frac{dh}{dz} \pm 2i\frac{dh}{dz} =
  \psi(z;\lambda)(1+h), \qquad h = h_\pm(z;\lambda).
\end{equation}
We can solve this equation for $\lambda\in B$ and $z\in S_R^\pm$, where
\begin{equation}\label{eq:BR2}
    S_R^+ = \{re^{i\theta}\,:\, r>R, \; -\pi/24<\theta<\pi/6\}, \qquad
    S_R^- = \{ \bar z \,:\, z\in S_R^+\}.
\end{equation}
This is done by solving the integral equation
\begin{equation*}
  h(z;\lambda) = -\int_{\gamma(z)} K(z,\zeta;\lambda)
  [1+h(\zeta;\lambda)]\,d\zeta,
\end{equation*}
where $K(z,\zeta;\lambda) = \frac{\beta}{2i}\big[ 1 -
E(z,\zeta;\lambda) \big]\psi(\zeta;\lambda)$,
$E(z,\zeta;\lambda)=e^{2i\beta\{\xi(\zeta;\lambda)-\xi(z;\lambda)\}}$,
$\gamma(z)$ is the path from $z$ to $\infty$ in the direction
$e^{i\beta\pi/6}$, and $\beta=\pm1$ records the sign in (\ref{eq:h:ODE2}).
Thus, we seek a solution of $h=\mathbb{K}[1+h]$, where
\begin{equation}\label{eq:K2}
  \mathbb{K}f(z;\lambda) = -\int_0^\infty K(z,z+\hat\gamma s;\lambda)
  f(z+\hat\gamma s;\lambda)\hat\gamma\,ds,
  \qquad \hat\gamma = e^{i\beta\pi/6}.
\end{equation}
Focusing on the $\beta=+1$ case, we will show below that (increasing
$R$ if necessary) there is a positive, decreasing function $g(s)$
defined for $s\ge0$ such that $\int_0^\infty g(s)\,ds\le 1/2$ and
$|K(z,z+\hat\gamma s;\lambda)|\le g(s)$ for $z\in S_R^+$, $\lambda\in
B$ and $s\ge0$.

Let $\mc{B}$ be the Banach space of (jointly)
holomorphic functions $f(z;\lambda)$ that are bounded on $S_R^+\times
B$, with norm $\|f\|_\mc{B} = \sup_{z,\lambda}|f(z;\lambda)|$.  Then
since the integrand of (\ref{eq:K2}) is holomorphic in $z$ and
$\lambda$ for fixed $s$ and is uniformly dominated by $
g(s)\|f\|_\mc{B}$, $\mathbb{K}$ maps $\mc{B}$ to $\mc{B}$ and has norm
$\le 1/2$. Indeed, continuity of $\mbb{K}f(z;\lambda)$ follows from
the dominated convergence theorem, and analyticity follows from
Morera's theorem and Fubini's theorem. The fixed-point iteration
$h=[\mbb K + \mbb K^2 + \mbb K^3 + \cdots]1$ leads to a holomorphic
function $h\in\mc B$ that satisfies $h=\mbb K[1+h]$ and
$\|h\|_\mc{B}\le 1$. A change of variables shows that for any $t\ge0$
we have
\begin{equation*}
  h(z+\hat\gamma t;\lambda) = - \int_t^\infty K(z+\hat\gamma t,
  z+\hat\gamma s;\lambda)[1+h(z+\hat\gamma s;\lambda)]\hat\gamma\,ds.
\end{equation*}
Applying $\hat\gamma^{-1}\partial_t[\cdots]$ and
$\hat\gamma^{-1}\partial_t\big\{\sqrt{\Psi/F}\hat\gamma^{-1}\partial_t
[\cdots]\big\}$ to this equation and setting $t=0$ shows that
$h(z;\lambda)$ satisfies (\ref{eq:h:ODE2}).
Finally, we will see below that
\begin{equation}\label{eq:K:bnd}
  |K(z,z+\hat\gamma s;\lambda)|\le g\big(s+3(|z|-R)/4\big), \qquad
  (z\in S_R^+, \, \lambda\in B).
\end{equation}
It follows that $|h(z;\lambda)|\le 2\int_{3(|z|-R)/4}^\infty
g(s)\,ds$, which converges to zero as $|z|\rightarrow\infty$.  In
particular, for real $x$, $h(x;\lambda)\rightarrow0$ as
$x\rightarrow\infty$.  As a result, $v_+(x;\lambda)\sim e^{+
  i\xi(x;\lambda)}$ for large $x$, and $u_+(x;\lambda)$ has the
form (\ref{eq:asym:u}) and depends analytically on $\lambda\in B$ for
fixed~$x$. In the $\beta=-1$ case, the same construction works on
$S_R^-$, and in fact the partial sums $h_-^{(n)}=\sum_{k=1}^n\mbb{K}^k1$
are related to those above by $h_-^{(n)}(z;\lambda) =
\overline{h_+^{(n)}(\bar z;\bar \lambda)}$ for $z\in S_R^-$, $\lambda\in B$.
As a result, $u_-(x;\bar\lambda) = \overline{u_+(x;\lambda)}$
for $x>R$. In (\ref{eq:x0:def}), the point $x_0$ where
$C(\lambda)=\tilde\Phi(x_0;\lambda)^{-1}\Phi(x_0;\lambda)$
is defined can be any number greater than $R$.

It remains to construct $g(s)$ so that $\int_0^\infty g(s)\,ds\le 1/2$
and (\ref{eq:K:bnd}) holds.  First we claim that $|E(z,z+\hat\gamma
s;\lambda)|$ is a decreasing function of $s$ when $\lambda\in B$ and
$z\in S_R^\pm$ are fixed, and therefore remains bounded by 1 for
$s\ge0$.  Focusing on the $\beta=+1$ case, this is equivalent to
claiming that $\im\{\xi(z+\hat\gamma s;\lambda)\}$ is an increasing
function of $s$.  This follows from $(d/ds)\xi(z+\hat\gamma s;\lambda)
= \sqrt{F(\zeta;\lambda)/\Psi(\zeta)}\hat\gamma$, where
$\zeta=z+\hat\gamma s$. Indeed, since $\zeta\in S_R^+$ and $\lambda\in
B$, the arguments used above to bound $F(z;\lambda)$ and
$\Psi(z)^{-1}$ away from the branch cut also imply
\begin{align*}
  \opn{Arg}\left(\sqrt{\frac{F(\zeta;\lambda)}{\Psi(\zeta)}}\hat\gamma\right) &=
  \jt \frac{1}{2}\opn{Arg} F + \frac{3}{2}\opn{Arg} z + \frac{1}{2}\opn{Arg}(1+\veps(z))
  + \opn{Arg}\hat\gamma \\
   &\in \left(
    -\frac{\pi}{12} - \frac{\pi}{16} - \frac{\pi}{24} + \frac{\pi}{6}\,,\,
    \frac{\pi}{12}+\frac{\pi}{4} + \frac{\pi}{24} + \frac{\pi}{6}\right) =
  \left(\frac{\pi}{48}\,,\,\frac{13\pi}{24}\right).
\end{align*}
Therefore, the imaginary part of $(d/ds)\xi(z+\hat\gamma s;\lambda)$
is positive.  Note that $\opn{Arg}\hat\gamma=\pi/6$ offsets the
negative contributions from the other terms, which is why adjusting
the contour of integration to point in the $\hat\gamma$
direction allows us to analytically continue $\lambda$ across the real
axis.

Since $|E(z,z+\hat\gamma s;\lambda)|\le1$ for $s\ge0$, the left-hand
side of (\ref{eq:K:bnd}) is bounded by $|\psi(z+\hat\gamma
s;\lambda)|$.  Increasing $R$ if necessary, we claim there is a
non-negative function $g(s)$ such that $\int_0^\infty g(s)\,ds\le1/2$
and $|\psi(z+\hat\gamma s;\lambda)|\le g\big(s+3(|z|-R)/4\big)\le
g(s)$ for $s\ge0$, $z\in S_R^+$, $\lambda\in B$. (We continue to
assume $\beta=+1$). To this end, we note that
\begin{gather*}
  \jt \psi(\zeta;\lambda) =
  \frac{1}{8}\underbrace{\Psi^{-1/2}\Psi'F'}_{C_1}F^{-3/2} -
  \frac{5}{16}\underbrace{\Psi^{1/2}(F')^2}_{C_2}F^{-5/2} +
  \frac{1}{4}\underbrace{\Psi^{1/2}F''}_{C_3} F^{-3/2}, \\
  \jt C_1 \doteq -\frac3{\sqrt8}\zeta^{-9/2}[1-\frac{45}{16}\zeta^{-4}], \quad
  C_2 \doteq \frac{\zeta^{-11/2}}{\sqrt{32}}[1-\frac{45}{16}\zeta^{-4}]^2, \quad
  C_3 \doteq -\frac{\zeta^{-9/2}}{\sqrt2}[1-\frac{135}{16}\zeta^{-4}],
\end{gather*}
where we used (\ref{eq:F:z:lam}) to compute $F'$ and $F''$, which are
independent of $\lambda$.  Increasing $R$ if necessary, we may assume
$|C_1|\le (4/3)|\zeta|^{-9/2}$, $|C_2|\le (1/5)|\zeta|^{-11/2}$ and
$|C_3|\le(4/5)|\zeta|^{-9/2}$ for $\zeta\in S_R$. We already
established that $|F(\zeta;\lambda)-\lambda_0|<\lambda_0/2$ for
$\lambda\in B$ and $\zeta\in S_R$, so
$|F(\zeta;\lambda)|>\lambda_0/2$.  It follows that
\begin{equation}\label{eq:psi:bnd}
  |\psi(\zeta;\lambda)|\le \frac{11}{30}\left(\frac{2}{\lambda_0}\right)^{3/2}
  |\zeta|^{-9/2} + \frac{1}{16}\left(\frac2{\lambda_0}\right)^{5/2}|\zeta|^{-11/2},
  \qquad (\zeta\in S_R, \, \lambda\in B).
\end{equation}
Finally, for $z\in S_R^+$ and $\zeta=z+\hat\gamma s$ with $s\ge0$, we
have $|\zeta|^2=|z|^2 + s^2 - 2|z|s\cos[5\pi/6 + \opn{Arg}(z)]$, by
the law of cosines.  Since $-\pi/24<\opn{Arg(z)}<\pi/6$ and
$\cos(19\pi/24)<-3/4$, $|\zeta|^2\ge |z|^2+s^2+(3/2)|z|s$. It follows
that $|\zeta|\ge (s+3|z|/4)$.  Defining $g(s)$ to be the right-hand
side of (\ref{eq:psi:bnd}) with $|\zeta|$ replaced by $(s+3R/4)$, we
have that $g(s)$ is a positive, decreasing function for $s\ge0$ and
 $|\psi(z+\hat\gamma s;\lambda)| \le g(s+3(|z|-R)/4) \le g(s)$ for
$z\in S_R^+$ and $\lambda\in B$, as claimed. Increasing $R$ if
necessary, $\int_0^\infty g(s)\,ds\le 1/2$.

\bibliographystyle{siam}
%\bibliography{refs}

\end{document}